\documentclass{article}

\usepackage{graphicx,psfrag}
\usepackage{pstricks}
\usepackage{subfigure}

\usepackage{amsmath}
\usepackage{amssymb}
\usepackage{amsthm}

\numberwithin{equation}{section}

\title{Newton's aerodynamic problem in media\\
of chaotically moving particles}

\author{Alexander Yu. Plakhov\\
        \texttt{plakhov@mat.ua.pt}
        \and
        Delfim F.~M.~Torres\\
        \texttt{delfim@mat.ua.pt}}

\date{Department of Mathematics\\
      University of Aveiro\\
      3810-193 Aveiro, Portugal}

\newtheorem{lemma}{Lemma}

\theoremstyle{remark}
\newtheorem{remark}{Remark}

\theoremstyle{definition}

\theoremstyle{definition}
\newtheorem{example}{Example}

      \newcommand {\al}   {\alpha}          \newcommand {\bt}  {\beta}
                \newcommand {\Gam}  {\Gamma}
                
              \newcommand {\ve}   {\varepsilon}
                 \newcommand {\vphi} {\varphi}
      \newcommand {\lam}  {\lambda}         
      
      \newcommand {\om}   {\omega}          
      \newcommand {\pl}   {\partial}        \newcommand {\s}    {\sigma}
           \newcommand {\UUU}  {{\cal U}}
      \newcommand {\VVV}  {{\cal V}}
      \newcommand {\RRRR}  {{\cal R}}       
      
      \newcommand {\RRR}  {{\mathbb R}}     \newcommand {\OOO}  {{\cal O}}
           
      \newcommand {\III}  {{\cal I}}        
              
      \newcommand {\HHH}  {{\cal H}}        \newcommand {\MMM}  {{\cal M}}
      \newcommand {\LLL}  {{\cal L}}        \newcommand {\KKK}  {{\cal K}}
      \newcommand {\AAA}  {{\cal A}}

\begin{document}

\maketitle

\begin{abstract}
We study the problem of minimal resistance for a body moving with
constant velocity in a rarefied medium of chaotically moving point
particles, in Euclidean space $\RRR^d$. The particles
distribution over velocities is radially symmetric. Under some
additional assumptions on the distribution function, the complete
classification of bodies of least resistance is made. In the case
of three and more dimensions there are two kinds of solutions: a
body similar to the solution of classical Newton's problem and a
union of two such bodies ``glued together'' by rear parts of their
surfaces. In the two-dimensional case there are solutions of five
different types: (a) a trapezium; (b) an isosceles triangle; (c)
the union of a triangle and a trapezium with common base; (d) the
union of two isosceles triangles with common base; (e) the union
of two triangles and a trapezium. The cases (a)--(d) are realized
for any distribution of particles over velocities, and the case
(e) is only realized for some distributions. Two limit cases are
considered, where the average velocity of particles is big and
where it is small as compared to the body's velocity. Finally,
using the obtained analytical results, we study numerically a
particular case: the problem of body's motion in a rarefied
homogeneous monatomic ideal gas of positive temperature in
$\RRR^2$ and in $\RRR^3$.
\end{abstract}

\section{Introduction}

In 1686, in his \emph{Principia} \cite{N}, I. Newton considered
the problem of body's motion in a homogeneous medium of point
particles. He assumed that collisions of the particles with the
body are absolutely elastic, the medium is very rare, so that the
particles do not mutually interact, and that initially the
particles are immovable, \textrm{i.e.}, thermal motion of
particles is not taken into account. These assumptions are not
satisfied in the ordinary conditions ``on earth'', but can be
approximately valid when considering motion of high-speed and
high-altitude flying vehicles such as missiles and artificial
satellites.

Newton considered the problem of finding the shape of body
minimizing resistance of the medium to the body's motion. He
solved this problem in the class of convex axially symmetric
bodies with the axis parallel to the body's velocity, of fixed
length along this axis and with fixed projection on a plane
orthogonal to the axis. Due to convexity of the body, each
particle hits the body at most once, and this fact allows one to
write down an explicit analytical formula for resistance. The body
of least resistance found by Newton can be described as follows:
the rear part of its surface is a flat disk, which is at the same
time the maximal cross section of the body by a plane orthogonal
to the symmetry axis. The front part of the surface is composed of
a smaller disk in the middle and of a strictly convex lateral
surface.

Let us also mention the two-dimensional analogue of Newton's
problem. Consider a class of convex figures in $\RRR^2$ that are
symmetric with respect to some straight line and have fixed length
along this line and fixed width; it is required to find the figure
from this class such that resistance to the motion of the figure
along this line is minimal. If the length does not exceed the
half-width, the solution is a trapezium with the angle 45$^0$ at
the base; elsewhere, the solution is an isosceles triangle.

Since the early 1990th the interest to Newton's problem revived.
In particular, there were obtained interesting results related to
minimization problems in wider classes of bodies obtained by
withdrawing or relaxing the conditions initially imposed by
Newton: axial symmetry \cite{BK}, \cite{BrFK}, \cite{LP1},
\cite{LP2} and convexity \cite{BFK}, \cite{CL1}, \cite{CL2},
\cite{P1}, \cite{P2}.

On the other hand, the assumptions of absolutely elastic
collisions and of absence of thermal motion in the medium are, at
the best, true only approximately. (Note that ``absence of thermal
motion'' means that the mean velocity of thermal motion of
particles is negligible as compared to the body's velocity.) In
\cite{friction}, the problem was studied under the more realistic
hypothesis of presence of friction at the moment of collision (so
that collisions are not absolutely elastic). In the present paper,
we address the minimization problem in a medium with thermal noise
of particles.

A convex and axisymmetric body moves in $\RRR^d$,\, $d \ge 2$
along its symmetry axis, in a medium of chaotically moving
particles; the medium is homogeneous, and distribution of the
particles over velocities is the same at every point. The
magnitude $V$ of velocity is constant. The length $h$ of the body
along the axis is fixed, and the maximal cross section of the body
by a hyperplane orthogonal to the axis is a unit
$(d-1)$-dimensional ball. We consider the problem of finding the
shape of body minimizing resistance of the medium. The main
results of this paper are as follows.

\begin{figure}
%
\subfigure[$V = 1$, $h = 1.97$]{
  \label{fig:sol3DV1Up321}
  \begin{minipage}[b]{0.5\textwidth}
  \centering
  \includegraphics[scale=0.5]{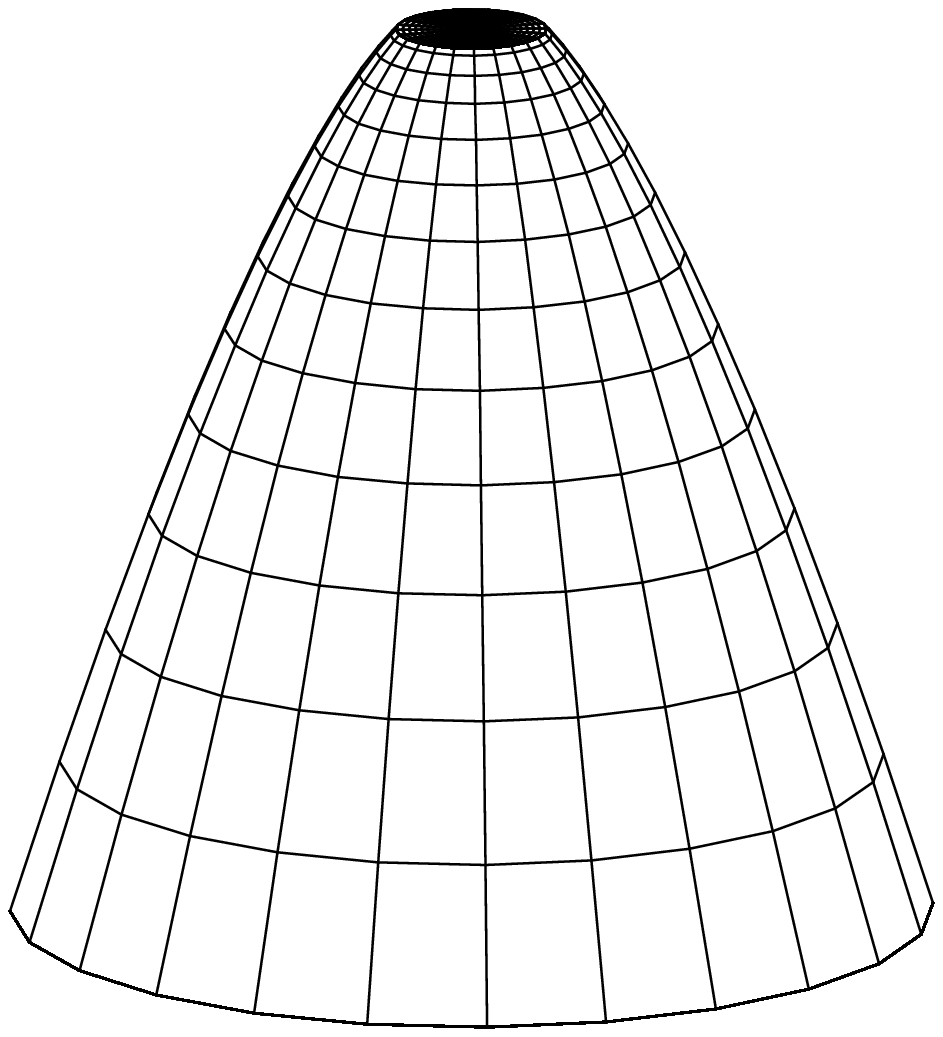}
  \end{minipage}}
\subfigure[$V = 1$, $h = 3.11$]{
  \label{fig:sol3DV1Up401}
  \begin{minipage}[b]{0.5\textwidth}
  \centering
  \includegraphics[scale=0.5]{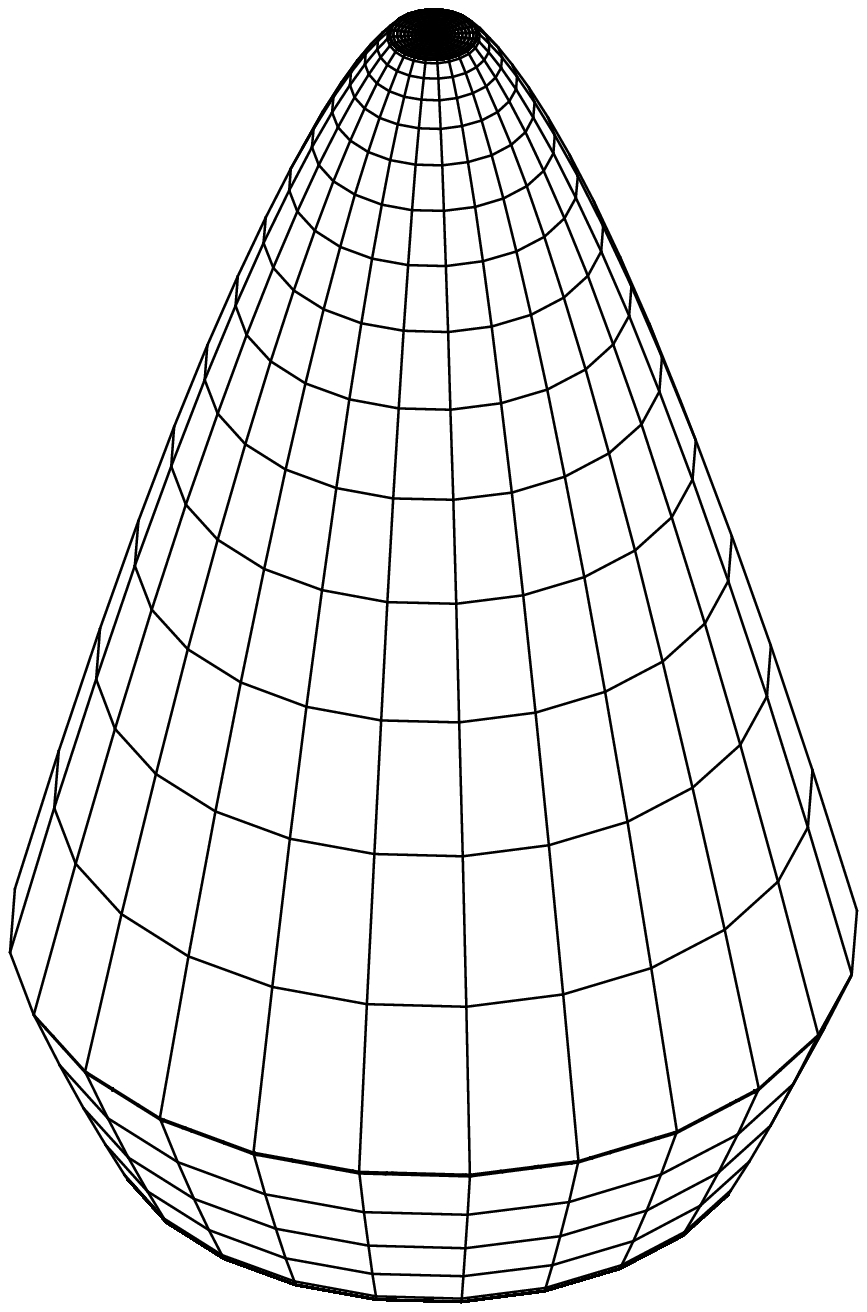}
  \end{minipage}}
\caption{Two solutions of the three-dimensional problem related to
motion in a rarefied monatomic homogeneous ideal gas. The mean
square velocity of gas molecules equals 1.} \label{fig:sol3DNLSd3}
\end{figure}

If $d \ge 3$, there are two different kinds of solutions. We shall
describe them in the case $d = 3$; if $d > 3$, the description is
quite similar. The solution of first kind is similar to the
solution of classical Newton's problem, that is, its surface can
be described in the same way as the surface of Newton's solution.
The solution of second kind is a union of two bodies similar to
Newton's solution ``glued together'' by rear parts of their
surfaces. The length (along the direction of motion) of the front
body is always more than the length of the rear body turned over.
The solution of first kind is realized for $h \le h_*$, and of
second kind, for $h > h_*$, where $h_* = h_*(V) > 0$ is a critical
value depending on $V$. The function $h_*(V)$ goes to infinity as
$V \to +\infty$ and to zero as $V \to 0^+$. The examples of
solutions of first and second kind are shown on
Fig.~\ref{fig:sol3DV1Up321} and on Fig~\ref{fig:sol3DV1Up401},
respectively. Here and in what follows the body is supposed to
move vertically upwards.

\begin{figure}
\vspace*{-1cm}
\begin{minipage}[t]{0.24\textwidth}
  \begin{center}
    \subfigure[$h = 0.7$]{
      \label{fig:2dh07u0p1068}
      \includegraphics[scale=0.2]{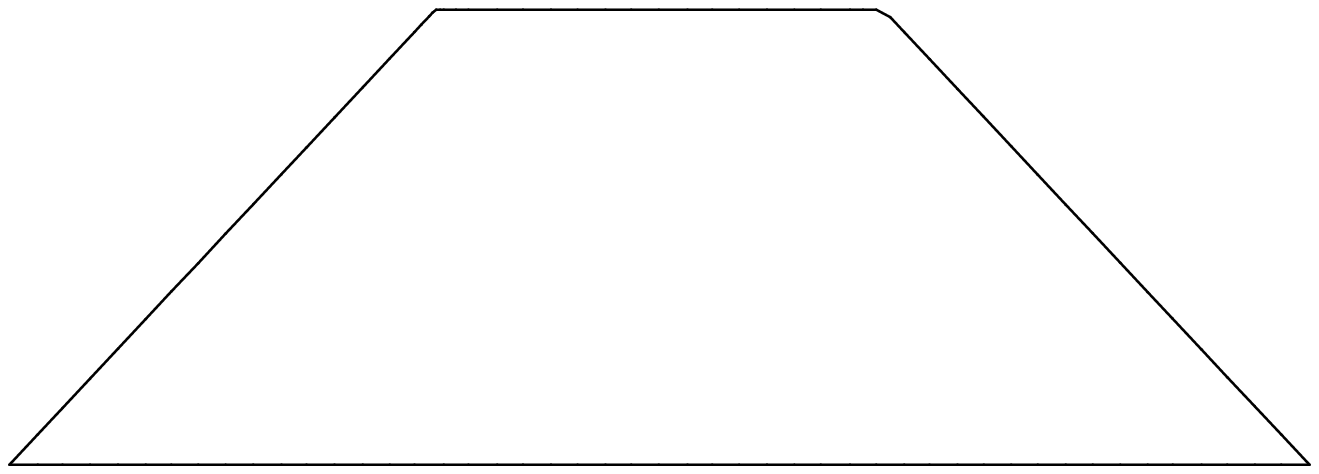}
    }\\[20pt]
    \subfigure[$h = 3$]{
      \label{fig:2d-tri-h3}
      \includegraphics[scale=0.2]{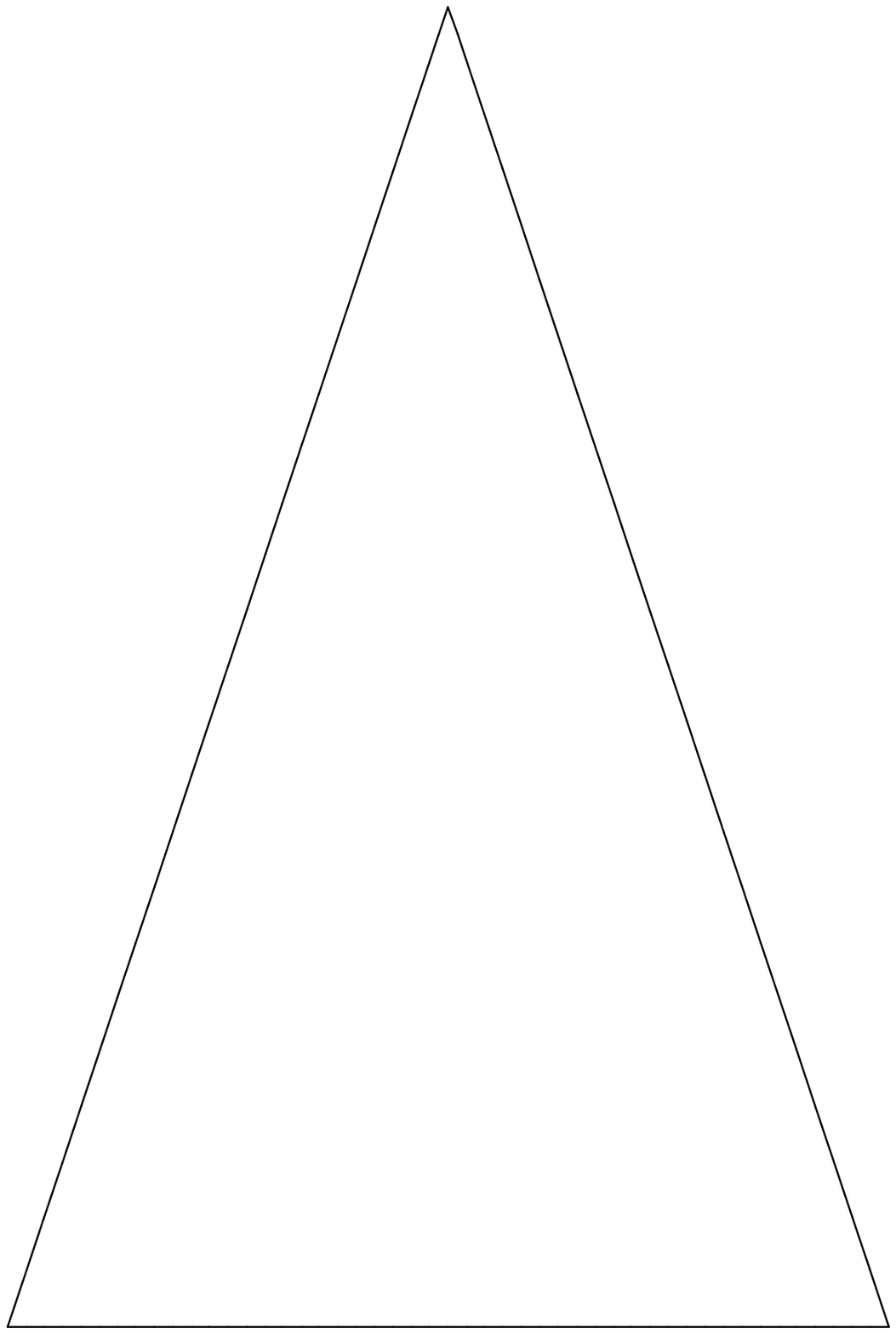}
    }
  \end{center}
\end{minipage}
  \begin{minipage}[t]{0.76\textwidth}
    \begin{minipage}[t]{0.5\textwidth}
      \begin{center}
        \subfigure[$h = 6$]{
          \label{fig:2d-tri-trap-hp479hm121}
          \includegraphics[scale=0.31]{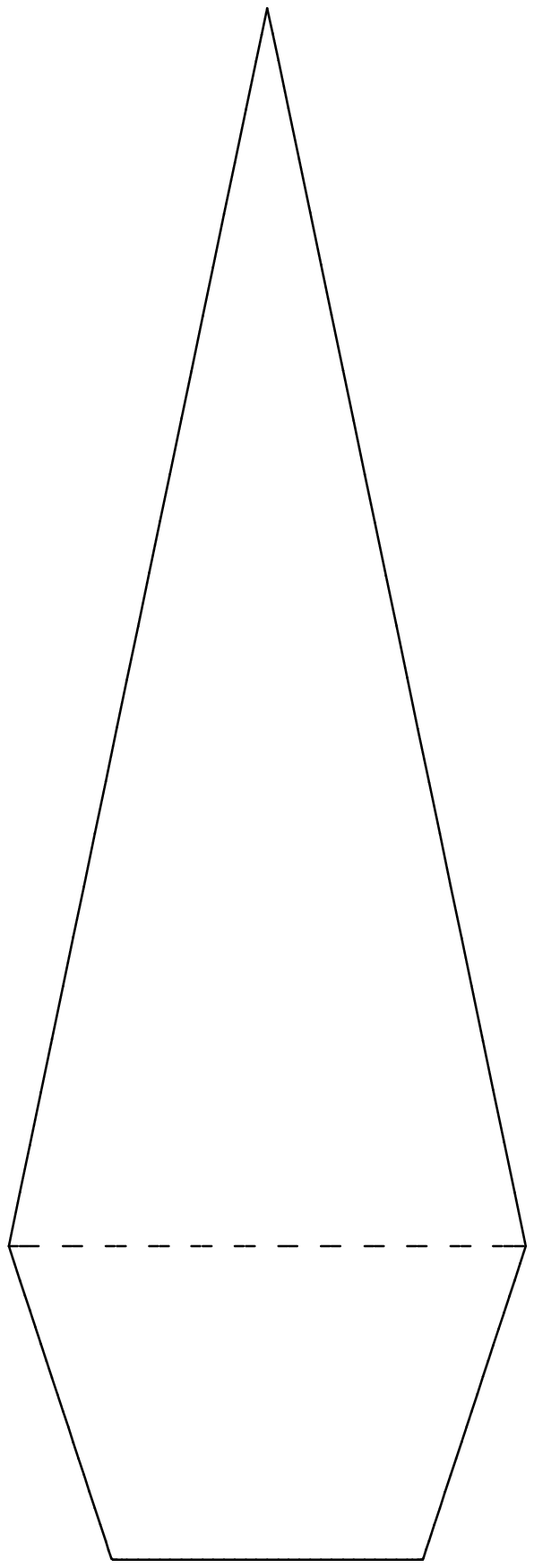}
        }\\[20pt]
        \subfigure[$h = 7.83$]{
          \label{fig:2d-tritri-hp479hm304}
          \includegraphics[scale=0.31]{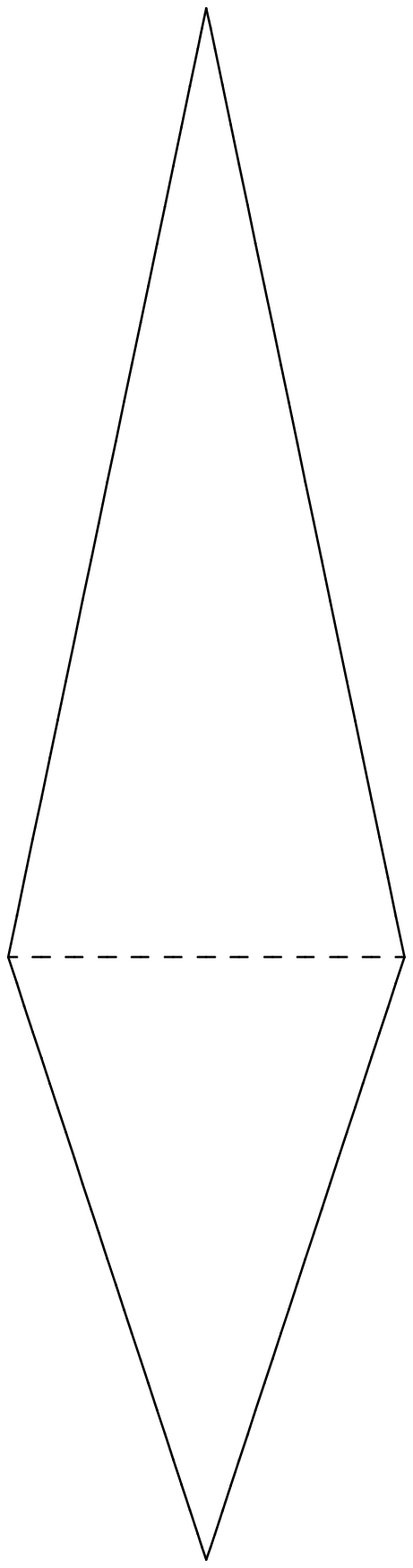}
        }
      \end{center}
    \end{minipage}
    \begin{minipage}[c]{0.5\textwidth}
      \begin{center}
        \subfigure[]{
          \label{fig:2dmixhp10-15hm1-23u1-072u2-2.24}
          \includegraphics[width=1.5cm,height=13cm]{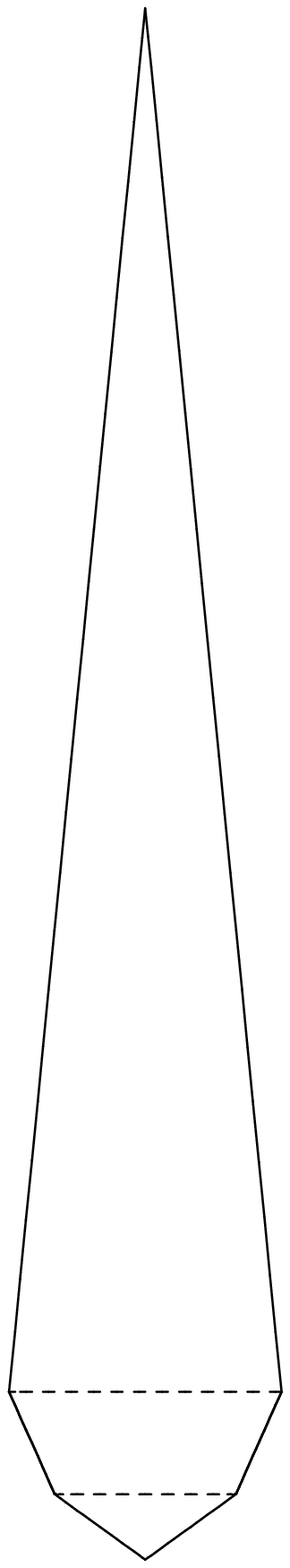}
        }
      \end{center}
    \end{minipage}
    \caption{The two-dimensional problem. The solutions corresponding to the cases (a) -- (d)
    are calculated numerically, for the motion with velocity $V = 1$ in a gas;
    the gas parameters are the same as on Fig.~1.}
    \label{fig:sol2D:fk:sk}
  \end{minipage}
\end{figure}

If $d = 2$, the classification of solutions is somewhat more
complicated. There are five different kinds of solutions: (a) a
trapezium, (b) an isosceles triangle, (c) the union of a triangle
and a trapezium, (d) the union of two isosceles triangles, (e) the
union of two triangles and a trapezium; see
Fig.~\ref{fig:2dh07u0p1068} --
Fig~\ref{fig:2dmixhp10-15hm1-23u1-072u2-2.24}. The solutions of
first kind are realized for $0 < h < u_+^0$, of second, for $u_+^0
\le h \le u_*$, of third, for $u_* < h < u_* + u_-^0$, of fourth
and fifth, for $h \ge u_* + u_-^0$. These values $u_+^0 =
u_+^0(V)$,\, $u_* = u_*(V)$,\, $u_-^0 = u_-^0(V)$ will be defined
in section 4.2; one has $0 < u_+^0(V) < u_*(V) < u_*(V) +
u_-^0(V)$. The solutions (a) -- (d) are realized for any
distribution of particles over velocities and for any positive
$V$; the solution (e) is realized only for some special
distributions and some values of $V$. The numerical computation of
a solution of kind (e) is a hard task, which is unsolved as yet.

%
%
%
%

In the limit cases, where the velocity of body is big and where it
is small as compared to the mean velocity of particles, the shape
of body of least resistance depends only on the length $h$, and
does not depend on the distribution of particles over velocities.
In the first limit case the optimal body coincides with the
solution of classical Newton's problem. In the second limit case,
for $d = 3$, the optimal body is a second kind solution symmetric
with respect to a plane perpendicular to the symmetry axis, the
inclination angle of the lateral surface at its upper and lower
points with respect to this plane being 51.8$^0$; and for $d = 2$,
the optimal body is one of the four figures:\, (a) a trapezium if
$0 < h < 1{.}272$;\, (b) an isosceles triangle if $h = 1{.}272$;\,
(c) the union of an isosceles triangle and a trapezium if $1{.}272
< h < 2{.}544$;\, (d) a rhombus if $h \ge 2{.}544$. In the cases
(a) -- (c) the inclination angle of lateral sides of these figures
with respect to the base equals 51.8$^0$, and in the case (d),
exceeds this value.

In a monatomic ideal gas the velocities of molecules are
distributed according to Gaussian law. Suppose that the mean
square velocity of molecules equals 1, then the kind of solution
is determined by two parameters: velocity of the body $V$ and its
length $h$. We define numerically the regions on the parameter
plane corresponding to different kinds of solutions;\, for some
special values of parameters we determine the shape of optimal
body and calculate the corresponding resistance. This work is made
in the two- and in the three-dimensional cases.

Imagine an observer travelling with the body through the medium;
he would detect a flux of particles falling on the immovable body.
In fact, this picture is more convenient for us, and will be taken
in the sequel.

The paper is organized as follows. In
section~\ref{Sec:CalPressResGC} the formulas for pressure of the
flux on the body's surface and for resistance force are derived;
two auxiliary lemmas of pressure distribution over the surface are
formulated; and the problem of minimal resistance is reduced to
the form more adapted for studying. In
section~\ref{sec:AuxMinPrb}, some auxiliary minimization problems
are solved. Using these results, in the next section we solve the
minimal resistance problem in general form. The solutions are
different for the cases $d = 2$ and $d \ge 3$. In
section~\ref{sec:GaussDV}, the obtained results are applied to the
flux corresponding to a rarefied monatomic homogeneous ideal gas
of positive temperature. In appendix~A, the auxiliary lemmas are
proved, and in appendix~B, asymptotic formulas for pressure
functions as $V \rightarrow 0^+$ are obtained.

\section{Calculation of pressure and resistance}
\label{Sec:CalPressResGC}

{\bf 2.1} \ Consider an immovable convex body $\cal B$ in
Euclidean space ${\mathbb R}^d$,\, $d \ge 2\,$ and a flux of
infinitesimal particles falling upon it. Velocities and masses of
the particles in general are different.  Let the function
$\rho(v)$ denote the distribution density over velocities of total
mass of particles in a unit volume, so that for any two
infinitesimal regions ${\cal V}$,\, $\mathcal{X} \subset \mathbb
R^d$ having $d$-dimensional volumes $|\mathcal{V}| = dv$,\,
$|\mathcal{X}| = dx$,\, the total mass of particles that are
contained in $\mathcal{X}$ and have velocities $v \in \cal V$
equals $\rho(v)\, dv\, dx$. Therefore, the particles' distribution
over velocities and density of the flux $\nu = \int_{\mathbb R^d}
\rho(v) dv$ are the same at each point and at each instant. It is
supposed that $\nu < \infty$.

An individual particle hitting the body at $x$ transmits the
impulse $m \cdot 2(v | n_x) n_x$ to the body, where $m$ and $v$
denote the particle's mass and its velocity before the collision,
$n_x$ means the outer unit normal to $\pl \cal B$ at $x$, and
$(\cdot\, |\, \cdot)$ means scalar product. Note that $(v | n_x) <
0$.

Let $\delta$ be an infinitesimal part of $\pl \cal B$ containing
$x$, and let $\mathcal{V}$ be an infinitesimal region containing
$v$, of volume $|\mathcal{V}| = dv$. The total mass of particles,
colliding with $\delta$ in a time interval $dt$ and having
velocity $v \in \mathcal{V}$, equals $dM = \rho(v)\, dv \cdot (v |
n_x)_{\!-}\, |\delta| dt$, where $|\delta|$ means
$(d-1)$-dimensional area of $\delta$, and $z_- := \max \{ -z,\, 0
\}$. The total impulse transmitted by these particles equals
$$
dM \cdot 2(v | n_x) n_x = -2 \rho(v)\, dv \cdot (v | n)_{\!-}^{\,\
2}\, |\delta| dt \cdot n_x.
$$
Integrating this value with respect to $v$, one obtains the total
impulse transmitted to $\delta$ per time $dt$,
$$
-2 \int_{\RRR^d} (v | n_x)_{\!-}^{\,\ 2} \rho(v)\, dv \cdot
|\delta| dt \cdot n_x.
$$
Dividing this value by $|\delta| dt$, one obtains that pressure of
the flux at $x$ equals $\pi(n_x)$, where
\begin{equation}\label{pi}
\pi(n) = - 2\int_{\RRR^d} (v | n)_{\!-}^{\,\ 2} \rho(v)\, dv \cdot
n.
\end{equation}

Integrating the pressure over $\partial\cal B$, one gets the total
force $R(\cal B)$ the flow is exerting on the body,
\begin{equation}\label{R}
R(\mathcal{B}) = \int_{\partial\mathcal{B}} \pi(n_x)\, d{\mathcal
H}^{d-1}(x),
\end{equation}
where ${\mathcal H}^{d-1}$ means $(d-1)$-dimensional Hausdorff
measure.
 \vspace{1mm}

{\bf 2.2} \ Denote $\RRR_+ := [0,\, +\infty)$, and denote by
$\AAA_d$ the set of functions $\sigma \in C^1(\RRR_+)$ such that
the function $\sigma'(r) /r$,\, $r > 0\,$ is negative, bounded
below, and monotone increasing, and
\begin{equation}\label{sigma}
\int_0^\infty r^2 \sigma(r)\, dr^{d} < \infty.
\end{equation}

\begin{remark}\label{r1}
Note that if $\s_1$,\, $\s_2 \in \AAA_d$ then $\s_1 + \s_2 \in \AAA_d$.
Also, if $\al$,\, $\bt > 0$,\, $\s \in \AAA_d$ and $\tilde \s(r) = \al \s(\bt r)$,
then $\tilde \s \in \AAA_d$.
\end{remark}

From now on, we suppose that the density function $\rho$ satisfies
the condition

\begin{quote}
{\bf A} \ $\rho(v) = \sigma(|v + V e_d|)$, where $\s \in
\AAA_d$,\, $V > 0$, and $e_d$ is the $d$th coordinate vector.
\end{quote}

Note that the relation (\ref{sigma}) implies that pressure
$\pi(n)$ (\ref{pi}) is always finite.

\begin{example}\label{e1}
Consider a rarefied homogeneous monatomic ideal gas in $\mathbb
R^3$ of absolute temperature $T > 0$. The distribution density of
molecules' mass over velocities equals $\sigma_{h} (|v|)$, where
$$
\sigma_h (r) = \nu \left( \frac{m}{2\pi kT} \right)^{3/2}\,
e^{-\frac{mr^2}{2kT}}
$$
(the Maxwell distribution); here $k$ is Boltzmann's constant and
$\nu$ is the gas density. Consider a body moving through the gas
with constant velocity of magnitude $V$ in the direction of third
coordinate vector $e_3$. In a frame of reference connected with
the body the distribution density over velocities equals
$\rho_h(v) = \sigma_h(|v + V e_3|)$. It is easy to check that
$\sigma_h \in \AAA_3$, so the condition A is fulfilled.
\end{example}

\begin{example}\label{e2}
Let, now, a rarefied ideal gas of temperature $T$ be a mixture of
$n$ homogeneous components, the $i$th component having density
$\nu_i$ and being composed of monatomic molecules of mass $m_i$.
Then the distribution density of molecules' mass over velocities
equals $\sigma_{nh} (|v|)$, where
$$
\sigma_{nh} (r) = \sum_{i=1}^n \nu_i\, \left( \frac{m_i}{2\pi kT}
\right)^{3/2}\, e^{-\frac{m_i r^2}{2kT}} \,.
$$
Taking account of remark 1, one concludes that $\sigma_{nh} \in
\AAA_3$.  As in the previous example, a body moves in the gas
along the third coordinate axis with velocity $V$. In a frame of
reference connected with the body the distribution density over
velocities equals $\rho_{nh}(v) = \sigma_{nh}(|v + V e_3|)$,
therefore, the condition A is also satisfied.
\end{example}

In the examples \ref{e1} and \ref{e2}, the force of resistance of
the gas to the body's motion is calculated according to (\ref{pi})
and (\ref{R}).
   \vspace{1mm}

In what follows, we shall also suppose the following condition to
be fulfilled:

\begin{quote}
{\bf B} \ The body $\mathcal{B}$ is convex, compact, and symmetric
with respect to the $d$th coordinate axis. Moreover, the maximal
cross section of the body by a hyperplane orthogonal to the
symmetry axis is a unit $(d-1)$-dimensional ball.
\end{quote}

By translation along the $d$th coordinate axis, the body can be
reduced to the form
\begin{equation*}
\mathcal{B} = \{ (x\text{'}, x_d) : |x\text{'}| \le 1,\
f_-(|x\text{'}|) \le x_d \le -f_+(|x\text{'}|)\}\, ,
\end{equation*}
where $x\text{'} = (x_1, \ldots, x_{d-1})$,\, $f_+$ and $f_-$ are
convex non-positive non-decreasing continuous functions defined on
$[0,\, 1]$. The length $h$ of body along the symmetry axis equals
\begin{equation*}
h = -f_+(0) - f_-(0).
\end{equation*}
 \vspace{1mm}

Now, let us specify the formulas for pressure $\pi(n)$ (\ref{pi})
and force $R(\cal B)$ (\ref{R}), using the conditions A and B.

At a regular point $x_+ = (x\text{'}, -f_+(|x\text{'}|))$ of the
upper part of the boundary $\partial \mathcal{B}$, the outer unit
normal vector is
\begin{equation}\label{nx+}
n_{x_+} = \frac{1}{\sqrt{f_+'(|x\text{'}|)^{\,2} + 1}}\
\left(f_+'(|x\text{'}|) \frac{x\text{'}}{|x\text{'}|},\ 1 \right),
\end{equation}
and using (\ref{pi}) and taking into account axial symmetry of
$\rho$ with respect to the $d$th coordinate axis, one finds that
pressure of the flux at this point equals
\begin{equation}\label{pinx+}
\pi(n_{x_+}) = -p_+ \left(f_+'(|x\text{'}|) \right) \cdot n_{x_+},
\end{equation}
where
\begin{equation}\label{p+}
p_+(u) := \left| \pi \left( \frac{1}{\sqrt{u^2+1}} \left( u,\, 0,
\ldots, 0,\, 1 \right)\right)\right|.
\end{equation}
Similarly, pressure of the flux at a regular point $x_- =
(x\text{'}, f_-(|x\text{'}|))$ of the lower part of $\partial
\mathcal{B}$ equals
\begin{equation}\label{pinx-}
\pi(n_{x_-}) = p_- \left(f_-'(|x\text{'}|) \right) \cdot n_{x_-},
\end{equation}
where
\begin{equation}\label{nx-}
n_{x_-} = \frac{1}{\sqrt{f_-'(|x\text{'}|)^{\,2} + 1}}\
\left(f_-'(|x\text{'}|) \frac{x\text{'}}{|x\text{'}|},\ -1 \right)
\end{equation}
and
\begin{equation}
\label{p-}
p_-(u) := -\left| \pi \left( \frac{1}{\sqrt{u^2+1}}
\left( u,\, 0, \ldots, 0,\, -1 \right)\right)\right|.
\end{equation}
From (\ref{p+}), (\ref{p-}), and (\ref{pi}) one obtains
\begin{equation}\label{formula}
p_\ve(u) = \ve \int_{\RRR^d} \frac{(v_1 u + \ve v_d)\!_-^{\,\
2}}{1 + u^2}\ \rho(v)\, dv, \ \ \ \text{where} \ \ \ \ve \in \{
-,\ + \}.
\end{equation}

Let us calculate $R(\mathcal{B})$. The integral in the right hand
side of (\ref{R}) is the sum of two integrals corresponding to the
upper and lower parts of $\partial \mathcal{B}$. Changing the
variable in each of these integrals and using the formulas
(\ref{nx+}), (\ref{pinx+}), (\ref{pinx-}), and (\ref{nx-}), one
obtains
$$
R(\mathcal{B}) = \int_{|x\text{'}|\le 1} p_+ (f_+'(|x\text{'}|))
\cdot \left( -f_+'(|x\text{'}|) \frac{x\text{'}}{|x\text{'}|},\ -1
\right)\, dx\text{'} + \ \ \ \ \ \ \ \ \ \ \ \
$$
$$
\ \ \ \ \ \ \ \ \ \ \ \ + \int_{|x\text{'}|\le 1} p_-
(f_-'(|x\text{'}|)) \cdot \left( f_-'(|x\text{'}|)
\frac{x\text{'}}{|x\text{'}|},\ -1 \right)\, dx\text{'}.
$$
Next, using that the functions $p_\ve (f_\ve'(|x\text{'}|))$ are
invariant, and the functions $\ve f_\ve'(|x\text{'}|)
\frac{x\text{'}}{|x\text{'}|}$ are anti-invariant with respect to
central symmetry $x\text{'} \to -x\text{'}$, one gets that
$$
\int_{|x\text{'}|\le 1} p_\ve (f_\ve'(|x\text{'}|)) \cdot \ve
f_\ve'(|x\text{'}|)\, \frac{x\text{'}}{|x\text{'}|}\ dx\text{'} =
0, \ \ \ \ \ve \in \{ -, + \},
$$
hence
$$
R(\mathcal{B}) = -\left( \int_{|x\text{'}|\le 1} p_+
(f_+'(|x\text{'}|))\, dx\text{'} + \int_{|x\text{'}|\le 1} p_-
(f_-'(|x\text{'}|))\, dx\text{'} \right)\ e_d.
$$
Therefore
\begin{equation*}
R(\mathcal{B}) = -a_{d-1}\, \left( \RRRR_+(f_+) + \RRRR_-(f_-)
\right) \cdot e_d,
\end{equation*}
where $a_{d-1}$ is the volume of a unit ball in $\RRR^{d-1}$, and
\begin{equation}\label{rpm}
\RRRR_\ve (f) = \int_0^1 \, p_\ve (f'(t))\, dt^{d-1}, \ \ \ \ \
\ve \in \{ -, + \}.
\end{equation}

Denote by ${\cal M}(h)$ the class of convex non-positive
non-decreasing continuous functions $f$ defined on $[0,\, 1]$ such
that $f(0) = -h$. Note that any function $f \in {\cal M}(h)$ is
differentiable everywhere except possibly on a finite or countable
set, and $f'$ is monotone, hence the integral (\ref{rpm}) is well
defined for any $f \in {\cal M}(h)$.
 \vspace{1mm}

{\bf 2.3} \ Thus, the problem of minimal resistance takes the
following form:
\begin{quote}
minimize $\RRRR_+(f_+) + \RRRR_-(f_-)$\\
provided that $f_+$ and $f_-$ are convex non-positive\\
non-decreasing functions satisfying the relation\\
$-f_+(0) - f_-(0) = h$.
\end{quote}
It will be solved in two steps. First, given $h_- \ge 0$,\, $h_+
\ge 0$, find
\begin{equation}\label{prs}
\inf_{f\in{\cal M}(h_-)} {\cal R}_-(f) \ \ \  \text{ and } \ \ \
\inf_{f\in{\cal M}(h_+)} {\cal R}_+(f).
\end{equation}
Second, given solutions $f_{h_-}^-$,\, $f_{h_+}^+$ of the problems
(\ref{prs}), find
$$
\mathrm{R}(h) := \inf_{h_+ + h_- = h} \left( {\mathcal
R}_+(f_{h_+}^+) + {\mathcal R}_-(f_{h_-}^-) \right).
$$
 \vspace{1mm}

{\bf 2.4} \ Let us formulate two auxiliary lemmas. Their proofs
are rather bulky, and are given in Appendix A.

Lemma~\ref{l1} states some properties of the functions $p_+$,\,
$p_-$, which will be needed in the subsequent sections.

\begin{lemma}\label{l1}
Let $\rho$ satisfy the condition A. Then

(a) there exist the limits $\lim_{u\to +\infty} p_\ve(u) =:
p_\ve(+\infty)$, besides\\
\hspace*{8mm} $p_+(+\infty) + p_-(+\infty) = 0$;

(b) $p_\ve \in C^1(\RRR_+)$,~ and  $p_\ve'(0) = \lim_{u\to
+\infty} p_\ve'(u) = 0$,~ $\ve \in \{ -,+ \}$;

(c) for $u > 0$,\, $p_+'(u) < p_-'(u)$;

(d) for $u > 0$,\, $p_+'(u) < 0$, and for any $u \ge 0$,\, $p_-(u)
> p_-(+\infty)$.
\end{lemma}

Lemma 2 specifies the form of functions $p_+$,\, $p_-$ for $d =
2$:~ the function $p_+$ has, in a sense, a simple behavior, and
the behavior of $p_-$ may be complicated. This specification will
be used in section 4.1 when constructing the body of least
resistance in two dimensions.

Designate $\s^{\al,\bt}(r) = \s(r) + \al\, \s(\bt r)$, where $\al
\ge 0$,\, $\bt > 0$. By virtue of remark~\ref{r1}, if $\s \in
\AAA_2$ then $\s^{\al,\bt} \in \AAA_2$, hence the function
$\rho^{\al,\bt}(v) = \s^{\al,\bt} (|v + V e_2|)$,\, $V
> 0$,\, defined on $\RRR^2$, satisfies the condition A. Denote by
$p_\ve^{\al,\bt}$,\, $\ve \in \{ -,+ \}$ the function
corresponding to the density $\rho^{\al,\bt}$, according to the
formula (\ref{formula}), and denote by $\bar p_\ve^{\al,\bt}$ the
maximal convex function defined on $\RRR_+$ that does not exceed
$p_\ve^{\al,\bt}$.

\begin{figure}
\subfigure[]{
  \label{fig:lemma2a}
  \begin{minipage}[b]{0.5\textwidth}
  \centering
  \psfrag{u}{$u$}
  \psfrag{ub}{$\bar{u}$}
  \psfrag{pp}{$p_+$}
  \includegraphics[scale=0.4]{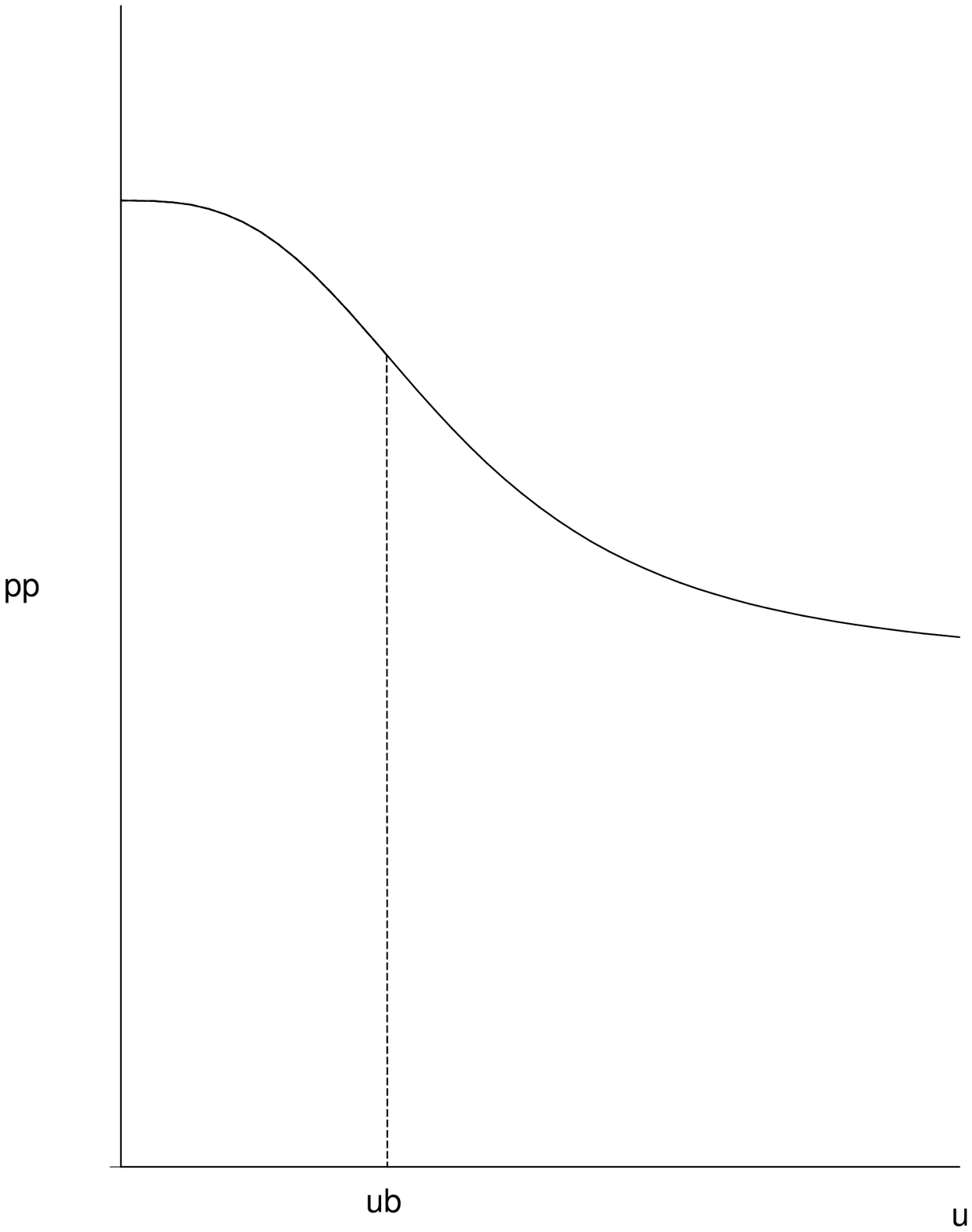}
  \end{minipage}}
\subfigure[]{
  \label{fig:lemma2b}
  \begin{minipage}[b]{0.5\textwidth}
  \centering
  \psfrag{u}{$u$}
  \psfrag{Oab}{$\cal{O}^{\al,\bt}$}
  \psfrag{pm}{$p_-^{\al,\bt}$}
  \includegraphics[scale=0.4]{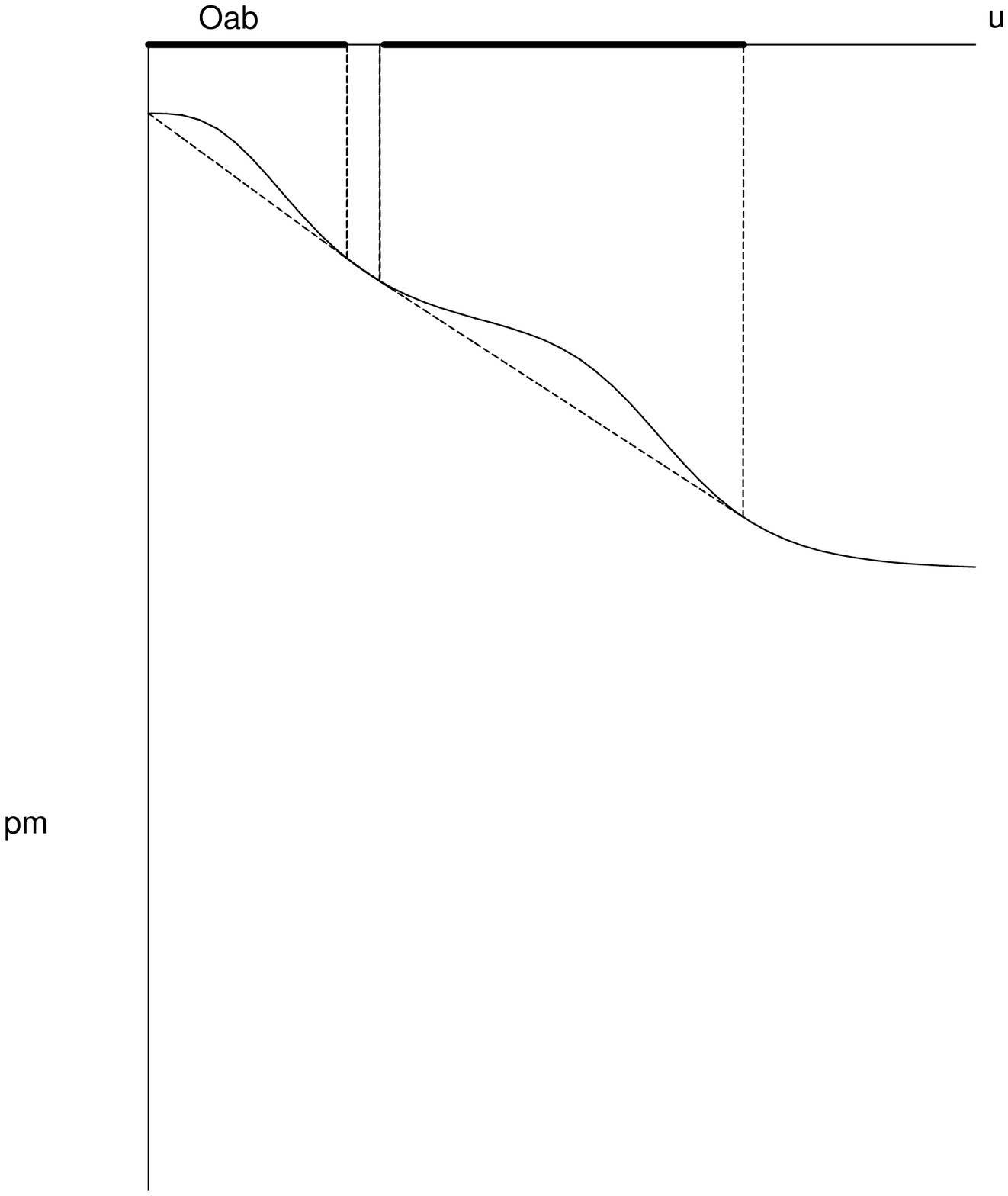}
  \end{minipage}}
\caption{}
\end{figure}

\begin{lemma}\label{l2} Let $d = 2$.

(a) If $\rho$ satisfies the condition A then for some $\bar{u} >
0$,\, $p_+'$ is monotone decreasing on $[0,\, \bar{u}]$ and
monotone increasing on $[\bar{u},\, +\infty)$ (see
Fig.~\ref{fig:lemma2a}).

(b) Suppose that $\s \in \AAA_2$,\, $V > 0$, and for any $n
> 0$ the function $r^n \s(r)$ monotonically decreases, for $r$
large enough. Then there exist $\al \ge 0$,\, $\bt > 0$ such that
the set $\OOO^{\al,\bt} = \{ u :\, p_-^{\al,\bt}(u) > \bar
p_-^{\al,\bt}(u) \}$ has at least two connected components (see
Fig.~\ref{fig:lemma2b}; $\OOO^{\al,\bt}$ is shown bold-faced on
the $x$-axis).
\end{lemma}

\section{Auxiliary minimization problems}
\label{sec:AuxMinPrb}

{\bf 3.1} \ The following lemma reduces the minimization problems
(\ref{prs}) to a simpler problem of minimization for a function
depending on a parameter.

\begin{lemma}\label{l3}
Let $p \in C(\RRR_+)$,\, $d \ge 2$,\, $\lam > 0$, and let a
function $f_h \in \MMM(h)$ satisfy the condition

\begin{quote} {\bf (C$_\lam$)} \ $f_h(1) = 0$,~ and for almost every $t$,\
$u = f_h'(t)$ is a solution of the problem
\begin{equation}\label{min}
t^{d-2}\, p(u) + \lam\, u \to \min.
\end{equation}
\end{quote}

Then $f_h$ is a solution of the minimization problem
\begin{equation}\label{minim}
\inf_{f\in\MMM(h)} \RRRR(f), \ \ \ \ \ \RRRR (f) = \int_0^1 \, p
(f'(t))\, dt^{d-1}.
\end{equation}
Moreover, any other solution of (\ref{minim}) satisfies the
condition (C$_\lam$) with the same $\lam$.
\end{lemma}

\begin{proof} In fact, the problem (\ref{minim}) can be considered
to be a degenerated case of the classical problem of optimal control
\cite{MR29:3316b}, and the statement of lemma is a consequence
of the Pontryagin maximum principle (see \cite{arXiv0404237}).
The proof we give here, however, is quite elementary
and does not appeal to the maximum principle (cf. \cite{T}).

For any $f \in \MMM(h)$ one has
\begin{equation}\label{ineq}
t^{d-2}\, p(f'(t)) + \lam\, f'(t) \ge t^{d-2}\, p(f_h'(t)) +
\lam\, f_h'(t)
\end{equation}
at almost every $t$. Integrating both sides of (\ref{ineq}) over
$t \in [0,\, 1]$, one gets
$$
\frac{1}{d-1}\, \int_0^1 \, p (f'(t))\, dt^{d-1} + \lam\, (f(1) -
f(0)) \ge
$$
\begin{equation}\label{ineqq}
\ge \frac{1}{d-1}\, \int_0^1 \, p (f_h'(t))\, dt^{d-1} + \lam\,
(f_h(1) - f_h(0)),
\end{equation}
and using that $f(1) \le 0 = f_h(1)$ and $f(0) = f_h(0) = -h$, one
obtains that $\RRRR (f) \ge \RRRR (f_h)$.

Next, suppose that $f \in \MMM(h)$ and $\RRRR (f) = \RRRR (f_h)$,
then, using the relation (\ref{ineqq}) and the equality $f(0) =
f_h(0)$, one gets that $f(1) \ge f_h(1) = 0$, hence $f(1) = 0$.
Therefore the inequality in (\ref{ineqq}) becomes equality, which,
in view of (\ref{ineq}), implies that
$$
t^{d-2}\, p(f'(t)) + \lam\, f'(t) = t^{d-2}\, p(f_h'(t)) + \lam\,
f_h'(t)
$$
for almost every $t$, hence $u = f'(t)$ is also a solution of
(\ref{min}), on a set of full measure. Thus, $f$ satisfies the
condition (C$_\lam$).
\end{proof}

{\bf 3.2} \ Assume, additionally, that $p \in C^1(\RRR_+)$
and that $p$ is bounded below. Denote by $\bar p$ the maximal
convex function defined on $\RRR_+$ that does not exceed
$p$. The function $\bar{p}$ is also continuously
differentiable, and for any $h \ge 0$,\, $u \ge 0$ one has
\begin{equation}\label{300-1}
p(u) \ge \bar p(u) \ge \bar p(h) + \bar p'(h) \cdot (u - h).
\end{equation}
Define the set $\OOO_p := \{ u:\, p(u) > \bar p(u) \}$. Obviously,
$\OOO_p$ is open, and hence is a union of a finite or countable
(maybe empty) set of disjoint open intervals.

The following lemma specifies the solution $f_h$ of the
minimization problem (\ref{minim}) in the case $d = 2$.

\begin{lemma}\label{l3a}
Suppose that the function $p$ is bounded below, $p \in
C^1(\RRR_+)$,\, $h \ge 0$,\, $\bar{p}'(h) < 0$, and denote by
$(h^{(-)},\, h^{(+)})$ the maximal interval contained in $\OOO_p$
such that $h^{(-)} \le h \le h^{(+)}$ (it may happen that
$h^{(-)} = h = h^{(+)}$, i.e., the interval is empty). \\
Then $h^{(+)} < +\infty$, and the following holds true:

(a) The function $f_h$ defined by
\begin{equation*}
f_h(t) = -h + h t, \ \ \ \ \ \ \ \ \ \ \ \ \ \ \ \ \ \ \ \ \ \ \ \ \ \ \ \ \
\end{equation*}
if $h^{(-)} = h = h^{(+)}$,~ and by
\begin{equation*}
f_h(t) = \left\{ \begin{array}{ll} -h + h^{(-)} t & \text{if } \ t \le t_0\,,\\
-h + h^{(-)} t_0 + h^{(+)} (t - t_0)  & \text{if } \ t \ge t_0\,,
\end{array} \right.
\ \ \text{where} \ t_0 = \frac{h^{(+)} - h}{h^{(+)} - h^{(-)}},
\end{equation*}
if $h^{(-)} < h^{(+)}$,~ is a solution of the minimization problem
\begin{equation}\label{300minimd=2}
\inf_{f\in\MMM(h)} \mathcal{R}(f), \ \ \ \ \ \mathcal{R}(f) =
\int_0^1 \, p (f'(t))\, dt\,,
\end{equation}
besides
\begin{equation}\label{300-infR}
 \hspace*{-0.5mm} \text{(b)} \hspace{26mm} \inf_{f\in\MMM(h)} \mathcal{R}(f) =
{\mathcal R}(f_h) = \bar p(h). \hspace*{27mm}
\end{equation}

(c) If $f$ is a solution of (\ref{300minimd=2}) then at almost
every $t$, the value $u = f'(t)$ satisfies the relations
$p(u) = \bar{p}(h) + \bar p'(h) \cdot (u - h)$ and $u \not\in \OOO_p$.
\end{lemma}

\begin{proof}
$p$ is bounded below and $\bar p'(h) < 0$, therefore there exists
a value $u > h$ such that $\bar p(u) = p(u)$, hence $h^{(+)} \le u
< +\infty$.

One has $\bar p(h^{(-)}) = p(h^{(-)})$,\, $\bar p(h^{(+)}) =
p(h^{(+)})$, and
\begin{equation}\label{300-2}
p(h^{(\pm)}) = \bar p(h) + \bar p'(h) \cdot (h^{(\pm)} - h).
\end{equation}
From (\ref{300-1}) and (\ref{300-2}) it follows that for any $u$
$$
p(u) - p(h^{(\pm)}) \ge \bar p'(h) \cdot (u - h^{(\pm)}),
$$
and designating $\lam = -\bar p'(h)$, one obtains
$$
p(u) + \lam u \ge p(h^{(\pm)}) + \lam h^{(\pm)}.
$$
This means that both $h^{(-)}$ and $h^{(+)}$ minimize the function
$p(u) + \lam u$.

Further, one easily sees that $f_h \in \mathcal{M}(h)$,\, $f_h(1)
= 0$, and the function $f_h'$ takes the values $h^{(-)}$ and
$h^{(+)}$ (which may coincide). Applying lemma \ref{l3}, one
obtains that $f_h$ is a solution of the problem
(\ref{300minimd=2}).

If $h^{(-)} \neq h^{(+)}$ then
\begin{eqnarray*}\label{300-3}
{\mathcal R}(f_h) & = & \int_0^{t_0} p(h^{(-)})\, dt +
\int_{t_0}^1 p(h^{(+)})\, dt = \nonumber\\
& = & \frac{h^{(+)} - h}{h^{(+)} - h^{(-)}}\ p(h^{(-)}) + \frac{h
- h^{(-)}}{h^{(+)} - h^{(-)}}\ p(h^{(+)}).
\end{eqnarray*}
On the other hand, excluding $\bar p'(h)$ from the relation
(\ref{300-2}), one obtains
\begin{equation*}
\frac{h^{(+)} - h}{h^{(+)} - h^{(-)}}\ p(h^{(-)}) + \frac{h -
h^{(-)}}{h^{(+)} - h^{(-)}}\ p(h^{(+)}) = \bar p(h),
\end{equation*}
hence ${\mathcal R}(f_h) = \bar p(h)$. If $h^{(-)} = h^{(+)} = h$
then ${\mathcal R}(f_h) = p(h) = \bar p(h)$,  so the formula
(\ref{300-infR}) is true.

Let, now, $f$ be a solution of (\ref{300minimd=2}). By lemma
\ref{l3}, for almost every $t$ the value $\hat{u} = f'(t)$
minimizes the function $p(u) + \lam u$, hence $p(\hat{u}) + \lam
\hat{u} = p(h^{(+)}) + \lam h^{(+)}$, and substituting $\lam =
-\bar p'(h)$ and using that $p(h^{(+)}) = \bar p(h) + \bar p'(h)
\cdot (h^{(+)} - h)$, one obtains
$$
p(\hat{u}) = \bar p(h) + \bar p'(h) \cdot (\hat{u} - h).
$$
Taking into account (\ref{300-1}), one gets that
$$
{p}(\hat{u}) = \bar p(\hat{u}) = \bar p(h) + \bar p'(h) \cdot
(\hat{u} - h),
$$
hence $\hat{u} \not\in \OOO_p$. Lemma \ref{l3a} is proved.
\end{proof}

{\bf 3.3} \ Suppose, in addition to the previous assumptions, that
there exists the limit $\lim_{u\to +\infty} p(u) =: p(+\infty)$
and that for any $u \in \mathbb{R}_+$,\, $p(u) > p(+\infty)$. Then
$\bar{p}'(u) < 0$ for any $u$. Denote $B = -\bar{p}'(0)$;\, one
has $B > 0$;\, the function $\bar{p}'$ is continuous and monotone
non-decreasing from $B$ to 0.

Suppose that $d \ge 3$, and introduce an auxiliary notation:~ $\om
= \frac{1}{d - 2}$,\, $q(u) = |\bar{p}'(u)|^{-\om}$,\, $Q(u) =
\int_0^u q(\nu)\, d\nu$. Both function $q$ and $Q$ are continuous
and monotone non-decreasing on $\RRR_+$;\, $q$ changes from
$B^{-\om}$ to $+\infty$, and $Q$ changes from 0 to $+\infty$.

\begin{lemma}\label{l3b}
Let $d \ge 3$,\, $h \ge 0$,\, $p \in C^1(\RRR_+)$,\, and let there
exist the limit $\lim_{u\to +\infty} p(u) = p(+\infty)$ and for
any $u \in \mathbb{R}_+$,\, $p(u) > p(+\infty)$. Then

(a) the set of values $U \ge 0$ satisfying the equation
\begin{equation}\label{eql3b}
U - \frac{Q(U)}{q(U)} = h
\end{equation}
is non-empty.

(b) Let $U$ be a solution of (\ref{eql3b}) and $t_0 =
\frac{q(0)}{q(U)}$. Introduce the function $f_h(t)$,\, $t \in
[0,\, 1]$ as follows:
\begin{quote}
for $t \in [0,\, t_0]$,\, \ $f_h(t) = -h$;

for $t \in [t_0,\, 1]$,\, \ $f_h$ is defined parametrically,
\begin{equation}\label{def1l3b}
 f_h = -h + \frac{u\, q(u) - Q(u)}{q(U)}, \hspace*{3mm}
\end{equation}
\begin{equation}\label{def2l3b}
t = \frac{q(u)}{q(U)}, \ \ \ \ \ \ \ \ u \in [0,\, U]. \hspace*{3mm}
\end{equation}
\end{quote}
The function $f_h$ is defined correctly, is strictly convex on $[t_0,\, 1]$,
and is a unique solution
of the minimization problem (\ref{minim}).

(c) The minimal value of $\mathcal{R}$ equals
\begin{equation}\label{infR-l3b}
\inf_{f\in\MMM(h)} \mathcal{R}(f) = {\mathcal R}(f_h) = \bar p(U)
+ \frac{Q(U)}{q(U)^{d-1}}\,.
\end{equation}
\end{lemma}

\begin{proof} (a) For arbitrary $c > 0$, one has
$$
U - \frac{Q(U)}{q(U)} = \int_0^U \left(1 - \frac{q(\nu)}{q(U)}
\right) d\nu \ge \int_0^c \left(1 - \frac{q(c)}{q(U)} \right)
d\nu,
$$
and taking into account that $\lim_{U\to +\infty} q(U) = +\infty$,
one concludes that for $U$ sufficiently large, $U -
\frac{Q(U)}{q(U)} > c/2$. This implies that the continuous
function $U - \frac{Q(U)}{q(U)}$ goes to $+\infty$ as $U \to
+\infty$;~ besides it vanishes at $U = 0$, hence the set of
solutions of (\ref{eql3b}) is non-empty. Therefore, (a) is proved.

Denote by $S(h)$ the set of points $u$ such that $\bar p(u) = \bar
p(h) + \bar p'(h) \cdot (u - h)$. Note that $S(h)$ coincides with
the connected component of $\bar\OOO_p$ containing $h$, if $h \in
\bar\OOO_p$, and $S(h) = \{ h \}$ otherwise. Thus, $S(h)$ is a
closed segment containing $h$; two segments $S(h_1)$,\, $S(h_2)$
either coincide or are disjoint. The condition $p(u) > p(+\infty)$
implies that all segments $S(u)$ are bounded. Obviously, the
family of non-degenerated segments (i.e., of those that are not
singletons) is at most countable. Denote by $\mathcal{S}$ the
union of non-degenerated segments. The function $\bar{p}'(u)$ is
monotone increasing on $\RRR_+ \setminus \mathcal{S}$ and is
constant on any non-degenerated segment $S(u) \subset
\mathcal{S}$, hence the set $\{ \bar{p}'(u), \ u \in \mathcal{S}
\}$ is at most countable.

Let $0 \le u_1 \le u_2$. After simple algebra one obtains
$$
\left[u_2 - \frac{Q(u_2)}{q(u_2)} \right] - \left[u_1 -
\frac{Q(u_1)}{q(u_1)} \right] =
$$
\begin{equation}\label{proofl3b-001}
= Q(u_1) \left[ \frac{1}{q(u_1)} - \frac{1}{q(u_2)} \right] +
\frac{1}{q(u_1)} \int_{u_1}^{u_2} \left( q(u_2) - q(\nu) \right)
d\nu .
\end{equation}
Both terms in the right hand side of (\ref{proofl3b-001}) are
non-negative, hence the function $U - \frac{Q(U)}{q(U)}$ is
monotone non-decreasing. If both $u_1$ and $u_2$ are solutions of
(\ref{eql3b}) then both terms in (\ref{proofl3b-001}) are equal to
zero, which implies that $q$ is constant on $[u_1,\, u_2]$; that
is, $\bar p'$ is constant on $[u_1,\, u_2]$;~ or, equivalently,
$u_1 \in S(u_2)$. This implies that the solution set of
(\ref{eql3b}) coincides with a segment $S(u)$.

(b) The relations (\ref{def1l3b}) and (\ref{def2l3b}) define
continuous functions $f_h$ and $t$ of $u$, varying from $-h$ to 0
and from $t_0$ to 1, respectively, when $u$ passes the interval
$[0,\, U]$. Moreover, the function $t = q(u)/q(U)$ is monotone
non-decreasing, each set $\{ u : q(u)/q(U) = t \}$ coincides with
some segment $S(u)$, and $f_h$ is constant on any such segment.
This means that the function $f_h(t)$ is defined correctly.

Let us calculate the left-hand and right-hand derivatives
$f_h'(t^-)$,\, $f_h'(t^+)$. Designate by $[u^-(t),\, u^+(t)]$ the
interval $\{ u : q(u)/q(U) = t \}$ and put $u = u^+(t)$. Let the
values $f_h + \Delta f_h$,\, $t + \Delta t$ correspond to the
argument $u + \Delta u$, and $\Delta t > 0$. One has
$$
\Delta t = \frac{q(u + \Delta u) - q(u)}{q(U)}\,,
$$
$$
\Delta f_h = \frac{(u + \Delta u)\, q(u + \Delta u) - Q(u + \Delta u)}{q(U)} -
\frac{u\, q(u) - Q(u)}{q(U)} =
$$
$$
= \frac{u\, (q(u + \Delta u) - q(u)) + \int_u^{u+\Delta u} (q(u + \Delta u) - q(\nu))
d\nu}{q(U)}\,,
$$
hence
\begin{equation}\label{proof3b-002}
\frac{\Delta f_h}{\Delta t} = u + \int_u^{u+\Delta u}
\frac{q(u + \Delta u) - q(\nu)}{q(u + \Delta u) - q(u)}\, d\nu\,.
\end{equation}
The integrand in the right hand side of (\ref{proof3b-002}) is
less than 1;~ due to definition of $u = u^+(t)$, one has $\Delta u
\to 0^+$ as $\Delta t \to 0^+$, therefore
$$
f_h'(t^+) = \lim_{\Delta t\to 0^+} \frac{\Delta f_h}{\Delta t} = u^+(t).
$$
Similarly, one finds
$$
f_h'(t^-) = u^-(t).
$$
Both functions $u^-(t)$ and $u^+(t)$ are positive, and for any
$t_1 < t_2$ one has $u^-(t_1) \le u^+(t_1) < u^-(t_2) \le
u^+(t_2)$. Therefore, the function $f_h$ is monotone increasing
and strictly convex on $[t_0,\, 1]$;~ moreover, it is constant on
$[0,\, t_0]$,\, $f_h(0) = -h$, and $f(1) = 0$. Thus, it is proved
that $f_h \in \mathcal{M}(h)$.

For any $t \in [t_0,\, 1]$, except possibly a countable set of
values, one has $u^-(t) = u^+(t) := \tilde{u} \in \RRR_+ \setminus
\mathcal{S}$, hence there exists the derivative $f_h'(t) =
\tilde{u}$. For any $u \ne \tilde{u}$ one has
$$
\bar{p}(u) > \bar{p}(\tilde{u}) + \bar{p}'(\tilde{u}) \cdot (u - \tilde{u}),
$$
and using that $p(u) \ge \bar{p}(u)$,\, $p(\tilde{u}) = \bar{p}(\tilde{u})$,
one obtains
$$
p(u) > p(\tilde{u}) + \bar{p}'(\tilde{u}) \cdot (u - \tilde{u}),
$$
hence
\begin{equation}\label{proof3b-003}
p(u) - \bar{p}'(\tilde{u}) \cdot u > p(\tilde{u}) - \bar{p}'(\tilde{u}) \cdot
\tilde{u}.
\end{equation}
Recall that $t = \frac{q(\tilde{u})}{q(U)} =
\frac{|\bar{p}'(U)|^\om} {|\bar{p}'(\tilde{u})|^\om}$ and $\om =
\frac{1}{d - 2}$. One has $t^{d-2} =
\frac{\bar{p}'(U)}{\bar{p}'(\tilde{u})}$, and multiplying both
parts of (\ref{proof3b-003}) by $t^{d-2}$ and designating
$-\bar{p}'(U) = \lam$, one obtains that
$$
t^{d-2} p(u) + \lam u > t^{d-2} p(\tilde{u}) + \lam \tilde{u}
$$
for any $u \ne \tilde{u}$. Thus, $\tilde{u} = f_h'(t)$ is a unique
value minimizing the function $t^{d-2} p(u) + \lam u$.

Let, now, $t \in (0,\, t_0)$. For $u > 0$ one has
$$
\bar{p}(u) \ge \bar{p}(0) + \bar{p}'(0) \cdot u.
$$
Using that $p(u) \ge \bar{p}(u)$,\, $p(0) = \bar{p}(0)$,\, $t_0^{2-d} =
\bar{p}'(0)/\bar{p}'(U) = - \bar{p}'(0)/\lam$, one obtains
$$
p(u) \ge p(0) + \bar{p}'(0) \cdot u = p(0) - \lam t_0^{2-d} u,
$$
hence for any $t \in (0,\, t_0)$
$$
p(u) + \lam t^{2-d} u > p(0),
$$
therefore the value $f_h'(t) = 0$ is a unique minimum of the
function $t^{d-2} p(u) + \lam u$. Applying lemma \ref{l3}, one
concludes that $f_h$ is a unique solution of (\ref{minim}).

(c) One has
\begin{equation}\label{proof3b-c0}
\mathcal{R}(f_h) = \int_0^{t_0} p(0)\, dt^{d-1} + \int_{t_0}^1 p(f_h'(t))\, dt^{d-1}.
\end{equation}
The first integral in the right hand side of (\ref{proof3b-c0}) equals
$$
\int_0^{t_0}\!\! ...\, = p(0) \left( \frac{q(0)}{q(U)}
\right)^{d-1}.
$$
Denote $\tilde{U} = \inf S(U)$. Changing the variable in the
second integral $t = q(u)/q(U)$,\, $u \in [0,\, \tilde{U}]$ and
taking into account that for almost every $t$,\, $f_h'(t) = u$,
one obtains that the second integral equals
$$
\int_{t_0}^1\!\! ...\, = \int_0^{\tilde{U}} p(u)\, d\! \left(
\frac{q(u)}{q(U)}\right) ^{d-1} = p(u) \left(
\frac{q(u)}{q(U)}\right)^{d-1} \bigg|_0^{\tilde{U}} -
\int_0^{\tilde{U}} \left( \frac{q(u)}{q(U)}\right)^{d-1} dp(u).
$$
Summing the first and the second integrals and taking into account
that $q(\tilde{U}) = q(U)$ and $d - 1 = 1 + 1/\om$, one obtains
\begin{equation}\label{proof3b-c1}
\mathcal{R}(f_h) = p(\tilde{U}) - \int_0^{\tilde{U}}
\left( \frac{q(u)}{q(U)}\right)^{1+1/\om} dp(u).
\end{equation}
The integral in (\ref{proof3b-c1}) can be represented as the sum
$$
\int_0^{\tilde{U}} (...) = \int_{[0,\tilde{U}]\setminus\mathcal{S}} (...)\,
+\, \sum_i \int_{S_i} (...)\,,
$$
where $S_i = S(u_i)$ are non-degenerated segments whose union
gives $[0,\, \tilde{U}]\cap \mathcal{S}$. If $u \in [0,\,
\tilde{U}] \setminus \mathcal{S}$, one has $p'(u) = \bar{p}'(u) =
-q(u)^{-1/\om}$, hence
$$
\int_{[0,\tilde{U}]\setminus\mathcal{S}} (...)\, =\,
-\int_{[0,\tilde{U}]\setminus\mathcal{S}}
\frac{q(u)}{q(U)^{1+1/\om}}\, du.
$$
Next, taking into account that the function $q$ is constant on
$S_i$ and that at endpoints of $S_i$,\, $p$ and $\bar{p}$
coincide, one obtains
$$
\int_{S_i} (...)\, =\, \int_{S_i} \left( \frac{q(u)}{q(U)}\right)^{1+1/\om} d\bar{u}\,
= -\int_{S_i} \frac{q(u)}{q(U)^{1+1/\om}}\, du.
$$
Summing these integrals, one gets
\begin{equation}\label{proof3b-c2}
\int_0^{\tilde{U}} \left( \frac{q(u)}{q(U)}\right)^{1+1/\om} dp(u)\, =\,
-\int_0^{\tilde{U}} \frac{q(u)}{q(U)^{1+1/\om}}\, du.
\end{equation}
Further, the function $q(u) = |\bar{p}'(u)|^{-\om}$ is constant on
$[\tilde{U},\, U]$, hence
\begin{equation}\label{proof3b-c3}
-\int_{\tilde{U}}^U \frac{q(u)}{q(U)^{1+1/\om}}\, du =
-q(U)^{-1/\om} (U - \tilde{U}) = \bar{p}'(U)\, (U - \tilde{U}) =
\bar{p}(U) - \bar{p}(\tilde{U}).
\end{equation}
Using that $\bar{p}(\tilde{U}) = p(\tilde{U})$ and applying
(\ref{proof3b-c1}), (\ref{proof3b-c2}), and (\ref{proof3b-c3}),
one gets
$$
\mathcal{R}(f_h) = \bar{p}(U) + \int_0^U \frac{q(u)}{q(U)^{1+1/\om}}\, du,
$$
and recalling that $Q$ is the primitive of $q$ and $1 + 1/\om = d
- 1$, one comes to the formula (\ref{infR-l3b}).
\end{proof}

\section{Solution of the minimal resistance problem}
\label{sec:SolPMR}

\subsection{Two-dimensional problem}

\subsubsection{Minimization of $\mathcal{R}_+$}

From statement (a) of lemma \ref{l2} it follows that there exist
values $u_+^0 > 0$ and $B_+ > 0$ such that
\begin{equation*}
\frac{p_+(u_+^0) - p_+(0)}{u_+^0} = p_+'(u_+^0) = -B_+
\end{equation*}
and
$$
\bar{p}_+(u) = \left\{ \begin{array}{ll} p_+(0) - B_+ u & \ \ \text{
if } \ 0 \le u \le u_+^0,\\
p_+(u) & \ \ \text{ if } \ u \ge u_+^0.
\end{array} \right.
$$
This implies that $\OOO_{p_+} = (0,\, u_+^0)$. Applying lemma
\ref{l3a}, one obtains that there exists a unique solution $f_h^+$
of the minimization problem
\begin{equation*}
\inf_{f\in\MMM(h)} \RRRR_+(f), \ \ \ \ \ \RRRR_+ (f) = \int_0^1 \,
p_+ (f'(t))\, dt,
\end{equation*}
defined by the relations
\begin{eqnarray}\label{fh1d=2}
f_h^+(t) & = & \left\{ \begin{array}{ll} -h & \text{ for } t \in [0,\ t_0],\\
-h + u_+^0 \cdot (t - t_0)  & \text{ for } t \in [t_0,\, 1],
\end{array} \right.\\
t_0 & = & 1 - h/u_+^0, \nonumber
\end{eqnarray}
if $0 \le h < u_+^0$, and
\begin{equation*}
f_h^+(t) = -h + ht, \ \  \ \ \ \ \ \ \ \ \ \ \ \ \ \ \ \ \ \ \ \ \ \ \
\end{equation*}
if $h \ge u_+^0$. The minimal resistance equals
$$
\inf_{f\in\MMM(h)} \mathcal{R}_+(f) = {\mathcal R}_+(f_h^+) = \bar p_+(h).
$$

\subsubsection{Minimization of $\mathcal{R}_-$}

Note that $p_-'(0) = 0$ and $\bar p_-'(0) < 0$, hence $\OOO_{p_-}$
contains an interval $(0,\, u_-^0)$,\, $u_-^0 > 0$, besides
$p_-(u_-^0) = p_-(0) + \bar p_-'(0) \cdot u_-^0$. Denote $B_- =
-\bar p_-'(0)$ and represent the open set $\OOO_{p_-}$ as the
union of its connected components $\OOO_i = (u_i^-,\, u_i^+)$,\,
$\OOO_{p_-} = \cup_i \OOO_i$. We shall suppose that the set of
indices $\{ i \}$ contains 1 and that $\OOO_1 = (0,\, u_-^0)$.
Statement (b) of lemma \ref{l2} and example \ref{e2} imply that in
some cases (for example, when considering pressure distribution of
a mixture of two homogeneous rarefied gases on the rear part of
surface of a moving body) $\OOO_{p_-}$ has at least two connected
components.

Consider the minimization problem
\begin{equation}\label{minimd=2-}
\inf_{f\in\MMM(h)} \RRRR_-(f), \ \ \ \ \ \RRRR_- (f) = \int_0^1 \,
p_- (f'(t))\, dt.
\end{equation}
Applying lemma \ref{l3a}, one obtains that there exists a solution
$f_h^-$ of this problem, besides
$$
\inf_{f\in\MMM(h)} \mathcal{R}_-(f) = {\mathcal R}_-(f_h^-) = \bar p_-(h).
$$
For $0 \le h < u_-^0$ one has
\begin{eqnarray}\label{fh1d=2-}
f_h^-(t) & = & \left\{ \begin{array}{ll} -h, & \text{ if } t \le t_0\\
-h + u_-^0 \cdot (t - t_0),  & \text{ if } t \ge t_0,
\end{array} \right.\\
t_0 & = & 1 - h/u_-^0. \nonumber
\end{eqnarray}
For $h \in \RRR \setminus \OOO_{p_-}$ one has
\begin{equation*}\label{fh2d=2-}
f_h^-(t) = -h + ht.
\end{equation*}
Finally, for $h \in \OOO_i$,\, $i \ne 1$ one has
\begin{eqnarray}\label{fhd=2a}
f_h^-(t) & = & \left\{ \begin{array}{ll} -h + u_i^- t,  & \text{ if } t \le t_i\\
-h + u_i^- t_i + u_i^+ (t - t_i),   & \text{ if } t \ge t_i\,,
\end{array} \right.\\
t_i & = & \frac{u_i^+ - h}{u_i^+ - u_i^-}. \nonumber
\end{eqnarray}

Notice that $f_h$ needs not be the unique solution of
(\ref{minimd=2-}). In some degenerated cases it may happen that
the right endpoint of some interval $\OOO_i$ coincides with the
left endpoint of another interval, $u_i^+ = u_j^-$,\, $i \ne j$;\,
then there exists a continuous family of solutions of
(\ref{minimd=2-}); the derivative of any function from this family
takes the values $u_i^-$,\, $u_i^+$, and $u_j^+$.

\subsubsection{Solution of the two-dimensional problem}

Thus, the problem of finding
$$
\mathrm{R}(h) = \inf_{h_+ + h_- = h} \left( {\mathcal
R}_+(f_{h_+}^+) + {\mathcal R}_-(f_{h_-}^-) \right)
$$
amounts to the problem
\begin{equation}\label{minp}
\min_{0 \le z \le h} p_h(z), \ \ \ \text{where} \ \ \ p_h(z) =
\bar p_+(z) + \bar p_-(h - z).
\end{equation}

The functions $\bar p_-'(u)$,\, $\bar p_+'(u)$ are continuous and
monotone non-decreasing, hence the function $p_h'(z)$,\, $0 \le z
\le h$\, is also continuous and monotone non-decreasing.

Using statement (c) of lemma \ref{l1}, one concludes that $B_+
> B_-$. Indeed, if $u_-^0 \le u_+^0$ then
$$
-B_- = \frac{p_-(u_-^0) - p_-(0)}{u_-^0} > \frac{p_+(u_-^0) -
p_+(0)}{u_-^0} \ge \frac{\bar p_+(u_-^0) - p_+(0)}{u_-^0} = -B_+,
$$
and if $u_-^0 > u_+^0$ then
$$
-B_- = p_-'(u_-^0) > p_+'(u_-^0) > p_+'(u_+^0) = -B_+.
$$
Thus, there exists a unique value $u_* > u_+^0$ such that $\bar
p_+'(u_*) = p_+'(u_*) = -B_-$. Consider four cases:

1) $0 < h < u_+^0$;

2) $u_+^0 \le h \le u_*$;

3) $u_* < h < u_* + u_-^0$;

4) $h \ge u_* + u_-^0$.

In the cases 1) and 2), for $0 \le z < h$, one has $p_h'(z) < \bar
p_+'(u_*) + B_- = 0$, hence $z = h$ is a unique value of argument
minimizing $p_h$. Therefore, the optimal values of $h_+$ and $h_-$
are $h_+ = h$,\, $h_- = 0$, and $f_{h_- =0}^- \equiv 0$.

1) $0 < h < u_+^0$. The function $f_{h_+ = h}^+$ is given by
(\ref{fh1d=2}).\, The body of least resistance is a trapezium, the
tangent of slope of its lateral sides being equal to $u_+^0$ (see
Fig.~\ref{fig:2dh07u0p1068}). The minimal resistance equals
$$
\mathrm{R}(h) = {\mathcal R}_+(f_{h_+}^+) + {\mathcal
R}_-(f_{h_-}^-) = p_+(0) - B_+\, h + p_-(0).
$$

2) $u_+^0 \le h \le u_*$. Here one has $f_{h_+ =h}^+(t) = -h + h\,
t$, hence the optimal body is an isosceles triangle (see
Fig.~\ref{fig:2d-tri-h3}), and
$$
\mathrm R(h) = p_+(h) + p_-(0).
$$

In the cases 3) and 4)\, one has $\bar p_+'(h) > -B_-$, hence
$p_h'(h) = \bar p_+'(h) - \bar p_-'(0) > 0$. On the other hand,
$p_h'(u_+^0) = \bar p_+'(u_+^0) - \bar p_-'(h - u_+^0) \le - B_+ +
B_- < 0$. Moreover, using statement (a) of lemma \ref{l2}, one
finds that the function $\bar p_+'(u) = p_+'(u)$,\, $u \in
[u_+^0,\, h]$ is monotone increasing, hence $p_h'$ is also
monotone increasing on this interval. It follows that the function
$p_h$ has a unique minimum $z \in (u_+^0,\, h)$ and $f_{h_+
=z}^+(t) = - z + zt$.

3) $u_* < h < u_* + u_-^0$. One has $p_h'(u_*) = \bar p_+'(u_*) -
\bar p_-'(h - u_*) = -B_- + B_- = 0$, therefore $p_h$ reaches its
minimal value at $z = u_*$;\, thus, the optimal values of $h_+$
and $h_-$ are $h_+ = u_*$,\, $h_- = h - u_*$. The function $f_{h_-
= h - u_*}^-$ is given by (\ref{fh1d=2-}). Here the optimal body
is the union of a triangle and a trapezium, as shown on
Fig.~\ref{fig:2d-tri-trap-hp479hm121}. The tangent of slope of
lateral sides of the trapezium equals $-u_-^0$. ÍThe minimal
resistance equals
$$
\mathrm R(h) = p_+(u_*) + p_-(0) - B_-\, (h - u_*).
$$

4) $h \ge u_* + u_-^0$. One has $p_h'(h - u_-^0) = \bar p_+'(h -
u_-^0) + B_- \ge 0$,\, hence the minimum of $p_h$ is reached at a
point $z \in (u_+^0,\, h - u_-^0]$, and the optimal values $h_+ =
z$,\, $h_- = h - z$, as well as the minimal resistance, are
obtained from the relations
$$
\begin{array}{ll}
h_+ + h_- = h, \\
p_+'(h_+) = p_-'(h_-), \\
h_+ \ge u_+^0, \  h_- \ge u_-^0,\\
\mathrm R(h) = p_+(h_+) + \bar p_-(h_-).
\end{array}
$$
Here one should distinguish between two cases.

4a) If $h_- \in \RRR \setminus \OOO_{p_-}$ then $f_{h_-}^- (t) =
-h_- + h_- t$, and the optimal body is a union of two isosceles
triangles with common base, of heights $h_+$ and $h_-$ (see
Fig.~\ref{fig:2d-tritri-hp479hm304}).

4b) If $h_-$ belongs to some interval $\OOO_i = (u_i^-,\,
u_i^+)$,\, $i \ne 1$, then $f_{h_-}^-$ is given by (\ref{fhd=2a}),
and the optimal body is the union of two isosceles triangles and a
trapezium (see Fig.~\ref{fig:2dmixhp10-15hm1-23u1-072u2-2.24}).

Note that the case 4b) is realized for the values $h$ from an open
(maybe empty) set contained in $(u_* + u_-^0,\, +\infty)$. This
set is defined by the parameters $\s$ and $V$. The case 4a) is
realized for the values $h$ from the complement of this set in
$(u_* + u_-^0,\, +\infty)$, which is always non-empty.

\subsection{The problem in three and more dimensions}

Let $d \ge 3$. Using lemma \ref{l3b}, one concludes that there
exists a unique solution $f_h^\ve$ of the problem
\begin{equation*}\label{minimd=3}
\inf_{f\in\MMM(h)} \RRRR_\ve(f), \ \ \ \ \ \RRRR_\ve (f) =
\int_0^1 \, p_\ve (f'(t))\, dt^{d-1},
\end{equation*}
besides
\begin{equation*}\label{Rd=3}
\RRRR_\ve (f_h^\ve) = \bar p_\ve(U) +
\frac{Q_\ve(U)}{q_\ve(U)^{d-1}}\,,
\end{equation*}
where $U$ is defined (not necessarily uniquely) by the relation
\begin{equation*}\label{U}
U - \frac{Q_\ve(U)}{q_\ve(U)} = h;
\end{equation*}
here $q_\ve(U) = |\bar{p}_\ve'(U)|^{-1/(d-2)}$,\ $Q_\ve(U) =
\int_0^U q_\ve(u)\, du$.

Thus, the problem
$$
\inf_{h_+ + h_- = h} \left( {\mathcal
R}_+(f_{h_+}^+) + {\mathcal R}_-(f_{h_-}^-) \right)
$$
amounts to the following problem:
\begin{equation*}
\inf_{\mathrm{h}_+(u_+) + \mathrm{h}_-(u_-) = h} \left(
\mathrm{r}_+(u_+) + \mathrm{r}_-(u_-) \right),
\end{equation*}
where
\begin{equation*}
\mathrm{r}_+(u) = \bar p_+(u) + \frac{Q_+(u)}{q_+(u)^{d-1}}, \ \ \
\ \mathrm{r}_-(u) = \bar p_-(u) + \frac{Q_-(u)}{q_-(u)^{d-1}}
\end{equation*}
and
\begin{equation*}
\mathrm{h}_+(u) = u -  \frac{Q_+(u)}{q_+(u)}, \ \ \ \ \
\mathrm{h}_-(u) = u -  \frac{Q_-(u)}{q_-(u)}, \ \ u \ge 0.
\end{equation*}

The functions $\mathrm{r}_\ve$ and $\bar p_\ve$,\ $\ve \in \{ -,\,
+ \}$ are monotone non-increasing, and $\mathrm{h}_\ve$ is
monotone non-decreasing from 0 to $+\infty$ when $u \in
\mathbb{R}_+$, besides any interval of constancy of one of these
functions is at the same time the interval of constancy of two
others. For each $z \ge 0$ choose $u$ such that $\mathrm{h}_\ve
(u) = z$ and put $\mathrm{r}^{(\ve)}(z) := \mathrm{r}_{\ve}(u)$,\,
$\pi^{(\ve)} (z) := \bar p_\ve' (u)$. From the stated above it
follows that the functions $\mathrm{r}^{(\ve)}$ and $\pi^{(\ve)}$
well defined on $\mathbb{R}_+$ and are monotone decreasing. Denote
$$
\mathrm{r}_h (z) = \mathrm{r}^{(+)} (z) + \mathrm{r}^{(-)} (h -
z).
$$
After some algebra one obtains that the function $\mathrm{r}_h$ is
differentiable and
\begin{equation}\label{mais}
\mathrm{r}_h' (z) = (d - 1)\, (\bar p_+' (u_+) - \bar p_-' (u_-)),
\end{equation}
where the values $u_+$,\, $u_-$ are chosen from the relations
$\mathrm{h}_+ (u_+) = z$,\, $\mathrm{h}_- (u_-) = h - z$. Both
values in the right hand side of (\ref{mais}),\, $\bar p_+' (u_+)
= \pi^{(+)} (z)$ and $- \bar p_-' (u_-) = -\pi^{(-)} (h - z)$, are
monotone increasing functions of $z$, hence $\mathrm{r}_h' (z)$ is
also monotone increasing from $\mathrm{r}_h' (0) = (d - 1)\, (\bar
p_+'(0) - \bar p_-'(U_-))$ to $\mathrm{r}_h' (h) = (d - 1)\, (\bar
p_+'(U_+) - \bar p_-'(0))$, where $U_+$ and $U_-$ are defined from
the relations $\mathrm{h}_+(U_+) = h$,\, $\mathrm{h}_-(U_-) = h$.
Note that $\bar p_+'(0) = -B_+$ and $\bar p_-'(U_-) \ge -B_-
> -B_+$, therefore $\mathrm{r}_h' (0) < 0$.

Recall that $u_*$ is defined in the subsection 4.1.3 by
$\bar p_+'(u_*) = -B_-$. Designate
\begin{equation}\label{hstar}
h_* := \mathrm{h}_+(u_*) = u_* - B_-^{\frac{1}{d-2}}\, Q_+(u_*)
\end{equation}
and consider two cases.

1)~ $h \le h_*$.~ One has $\mathrm{h}_+(U_+) \le
\mathrm{h}_+(u_*)$, hence $U_+ \le u_*$, therefore $\mathrm{r}_h'
(h) = (d - 1)\, (\bar p_+'(U_+) + B_-) \le (d - 1)\, (\bar p_+'
(u_*) + B_-) = 0$. This implies that $\mathrm{r}_h' (z) < 0$ for
$z \in [0,\, h)$, hence the function $\mathrm{r}_h$ has a unique
minimum at the point $z = h$, which corresponds to the values $h_+
= h$,\, $h_- = 0$. The minimal resistance equals
$$
\mathrm R(h) = \bar p_+(u_+) + Q_+(u_+)\, q_+(u_+)^{-d+1} + \bar
p_-(0).
$$

2)~ $h > h_*$.~ One has $U_+ > u_*$, therefore $\mathrm{r}_h' (h)
= (d - 1)\, (\bar p_+'(U_+) + B_-) > 0$. On the other hand,
$\mathrm{r}_h' (0) < 0$. Hence, there exists a unique value $z \in
(0,\, h)$ such that $\mathrm{r}_h' (z) = 0$. Thus, the function
$\mathrm{r}_h$ has a unique minimum at $z$; the optimal values of
$h_+$,\, $h_-$ are $h_+ = z > 0$,\, $h_- = h - z > 0$. These
values and the related auxiliary values $u_-$,\, $u_+$ are
uniquely defined from the system of four equations
\begin{equation*}
\begin{array}{l}
h_+ = u_+ -  Q_+(u_+) / q_+(u_+),
\\
h_- = u_- - Q_-(u_-) / q_-(u_-),
\\
h_+ + h_- = h,
\\
\bar p_+'(u_+) = \bar p_-'(u_-),
\end{array}
\end{equation*}
and the minimal resistance equals
$$
\mathrm R(h) = \bar p_+(u_+) + Q_+(u_+)\, q_+(u_+)^{-d+1} + \bar
p_-(u_-) + Q_-(u_-)\, q_-(u_-)^{-d+1}.
$$

\subsection{The limit cases}

Consider heuristically the limit behavior of solutions as $V \to
+\infty$ and as $V \to 0$, with fixed $h$ and $\sigma$. We shall
denote the pressure and resistance functions by $p_\pm(u,V)$ and
$\mathrm{R}(h,V)$, thus explicitly indicating dependence of these
functions on the parameter $V$.

\subsubsection{$\mathrm{V \to +\infty}$}

Denote by $\tilde p_\ve(u,V) = V^{-2} p_\ve(u,V)$,\, $\ve \in \{
-,+ \}$ the reduced pressure, and by $\tilde{\mathrm{R}}(h,V) =
V^{-2} \mathrm{R}(h,V)$, the minimal reduced resistance. One has
$$
\begin{array}{l}
\tilde p_+(u,V) = 1/(1 + u^2)\, + o\,(1),\\
\tilde p_-(u,V) = o\,(1), \ \ \tilde p_-'(u,V) = o\,(1), \ \ \ V
\to +\infty;
\end{array}
$$
in other words, as $V \to +\infty$, the functions $\tilde
p_+(u,V)$ and $\tilde p_-(u,V)$ tend to $1/(1 + u^2)$ and to 0,
respectively. These limit functions determine pressure
distribution on the front part and on the rear part of body's
surface in Newton's classical problem.

Consider the cases $d = 2$ and $d = 3$ separately.

$\mathbf{d = 2}$ \ If $h < 1$ then for $V$ sufficiently large, the
figure of least resistance is a trapezium, and the inclination
angle of its lateral sides tends to $45^0$ as $V \to +\infty$.\,
If $h > 1$ then for $V$ sufficiently large, the figure of least
resistance is an isosceles triangle coincident with the solution
of two-dimensional Newton's problem.

$\mathbf{d = 3}$ \ For $V$ sufficiently large, the body of least
resistance is the first kind solution.
The front part of its surface is the the graph of a function defined
on a unit circle; as $V \to +\infty$, this function uniformly converges to
the function that describes Newton's classical solution with the same $h$.

The case $d > 3$ is similar to the three-dimensional one.

Finally, the limit value of minimal reduced resistance
$\tilde{\mathrm R}(h, \infty)$ is equal to the
resistance of Newton's optimal solution multiplied by the density of
particles' flow $\nu = \int_{\mathbb{R}^d} \s(|v|)\, dv$.

\subsubsection{$V \to 0^+$}

In this limit case one has
\begin{equation*}\label{asympt d=2}
p_\ve (u,V) = \ve b^{(d)} + V\, \frac{c^{(d)}}{\sqrt{1 + u^2}}\, +
o(V), \ \ \ \ \ \ve \in \{ -,+ \},
\end{equation*}
where
\begin{equation}\label{c1c2}
b^{(2)} = \frac{\pi}{2}\, \int_0^{+\infty} \sigma(r) r^3\, dr, \ \
\ \ \ c^{(2)} = 4 \int_0^{+\infty} \sigma(r) r^2\, dr,
\end{equation}
\begin{equation}\label{k1k2}
b^{(3)} = \frac{2\pi}{3}\, \int_0^{+\infty} \sigma(r) r^4\, dr, \
\ \ \ \ c^{(3)} = 2\pi \int_0^{+\infty} \sigma(r) r^3\, dr.
\end{equation}

This formula will be derived in Appendix B. One readily obtains
that $u_-^0$ and $u_+^0$ tend to the value $a := \sqrt{(1 + \sqrt
5)/2} \approx 1{.}272$, and $B_\pm = V \cdot a^{-5} + o(V)$.
Taking into account that $\bar p_+'(u) < \bar p_-'(u) < 0$, one
concludes that $u_*$ tends to the same value $a$, and $u_+^0 +
u_*$ tends to $2a$.

Let us describe the shape of optimal body and determine the
minimal reduced resistance $\hat{\mathrm{R}}(h,V) = V^{-1}
\mathrm{R}(h,V)$ in the limit $V \to 0^+$. We shall distinguish
between two cases:~ $d = 2$ and $d = 3$.

$\mathbf{d = 2}$ \hspace*{1.2mm} (a) $0 < h < a$: the optimal body
is a trapezium.

\hspace*{11.5mm} (b) $h = a$: an isosceles triangle.

\hspace*{11.5mm} (c) $a < h < 2a$: the union of a triangle and a
trapezium.

\hspace*{11.5mm} (d) $h \ge 2a$: a rhombus.

\begin{figure}
\subfigure[$h = 0.5$]{
  \begin{minipage}[b]{0.5\textwidth}
  \begin{center}
    \includegraphics[scale=0.3]{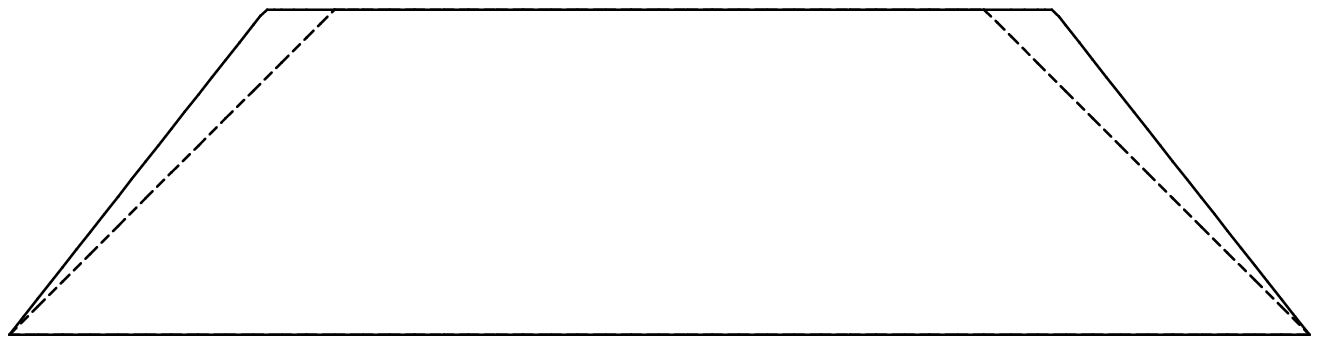}
  \end{center}
  \end{minipage}
} \subfigure[$h = 1.27$]{
  \begin{minipage}[b]{0.5\textwidth}
  \begin{center}
    \includegraphics[scale=0.3]{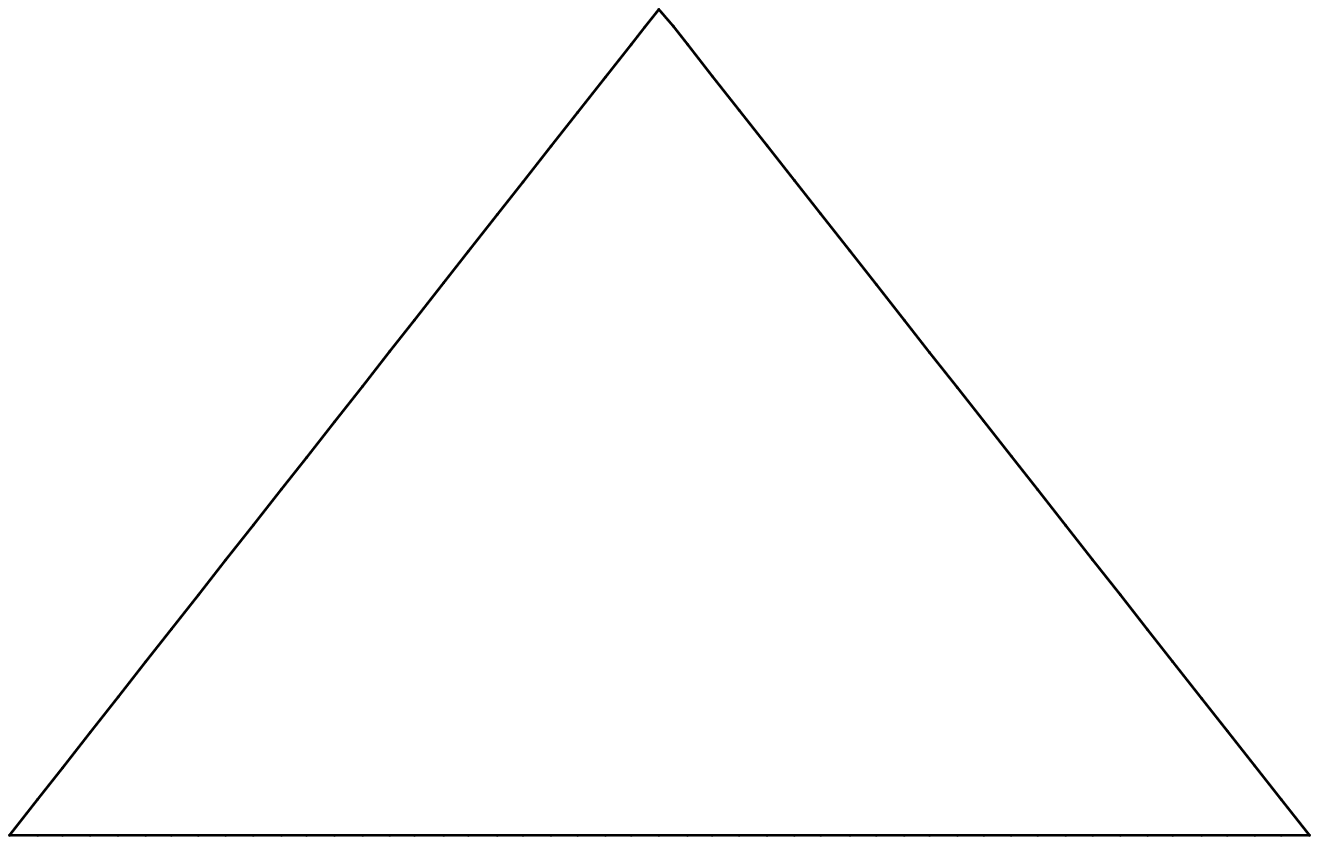}
  \end{center}
  \end{minipage}
} \subfigure[$h = 1.77$]{
  \begin{minipage}[b]{0.5\textwidth}
  \begin{center}
    \includegraphics[scale=0.3]{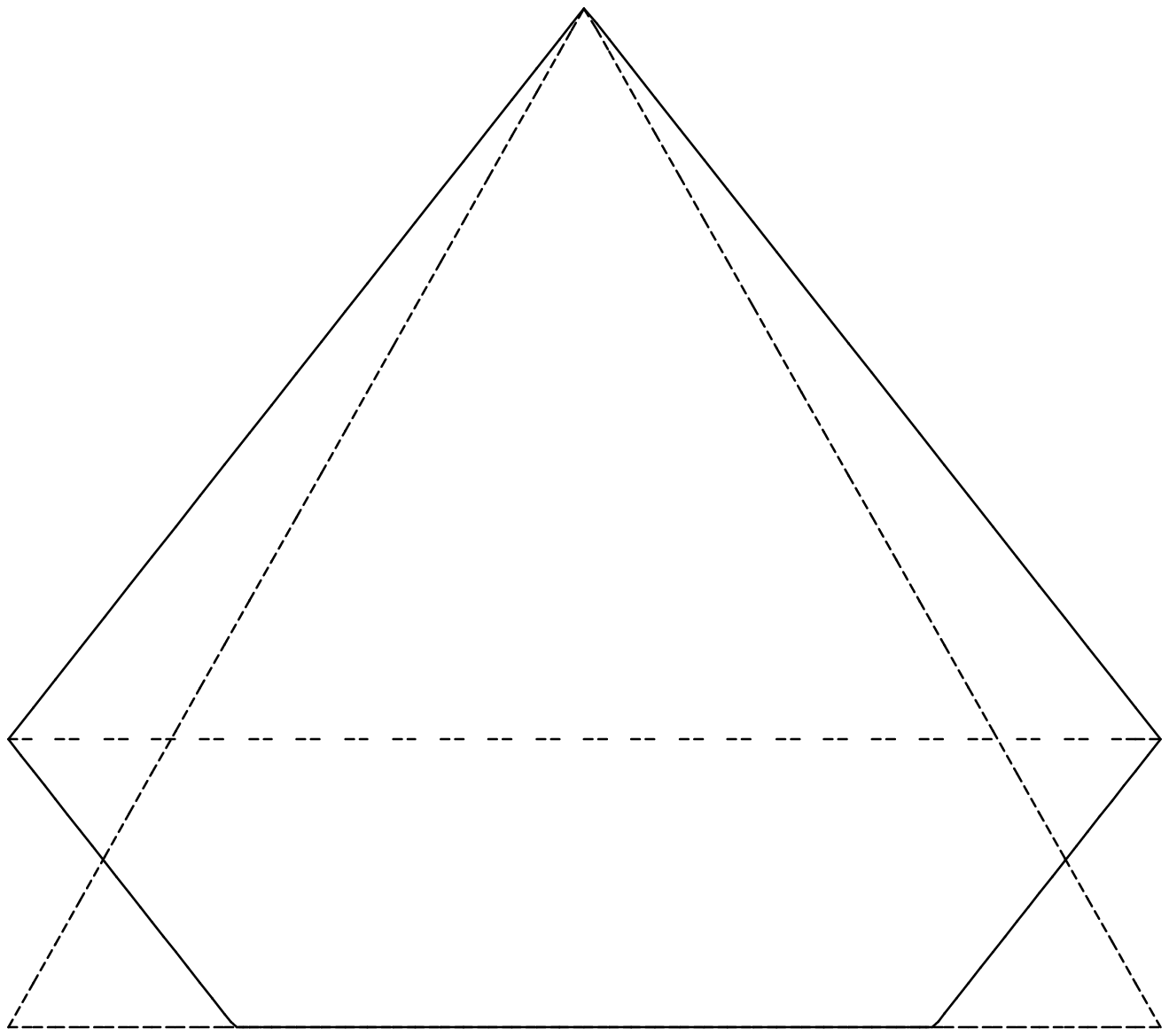}
  \end{center}
  \end{minipage}
} \subfigure[$h = 3$]{
  \begin{minipage}[b]{0.5\textwidth}
  \begin{center}
    \includegraphics[scale=0.3]{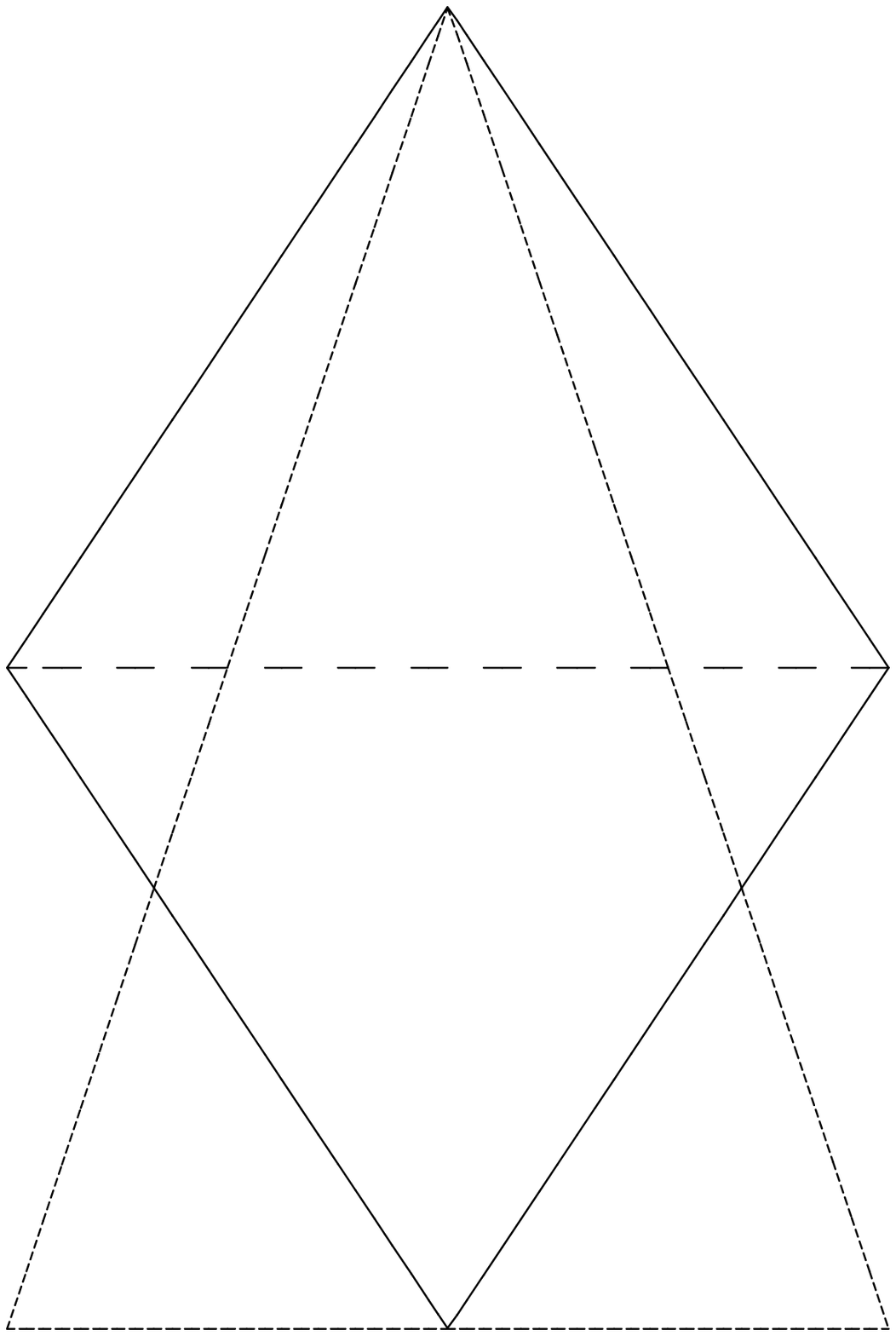}
  \end{center}
  \end{minipage}
} \caption{Two-dimensional case. Solutions in the limit $V
\rightarrow 0^+$ are shown by solid line. The corresponding
solutions of Newton's classical problem are shown by dashed line.
The case (b) is the unique one where these two solutions coincide.}
\label{fig:2d-cl}
\end{figure}

In the first three cases, the tangent of slope of lateral sides of
optimal figures equals $\arctan a \approx 51.8^0$, and in the last
case, exceeds this value. The examples of optimal figures are shown on
Fig.~\ref{fig:2d-cl}.

%
%
%
%

The minimal reduced resistance equals
\begin{equation}\label{Rasympt d=2}
\hat{\mathrm R}(h,V) = 2 c^{(2)} \bar p(h/2)\, +\, o(1), \ \ \ \ \
V \to 0^+,
\end{equation}
where
\begin{equation}\label{onemore}
\bar p(u) = \left\{
\begin{array}{ll}  1 - a^{-5}\, u,  & \text{ if } \ u \le a,\\
1 /\sqrt{1 + u^2}, & \text{ if } \ u \ge a.
\end{array}
\right.
\end{equation}

$\mathbf{d = 3}$ \ Let
\begin{equation*}
q(u) = \left\{
\begin{array}{ll}  ~a^5 & \text{if } u \le a\\
\frac{(1 + u^2)^{3/2}}{u} & \text{if } u \ge a,
\end{array}
\right.
\end{equation*}
\begin{equation*}
 Q(u) = \left\{
\begin{array}{ll}  ~a^5 u & \text{if } u \le a\\
\sqrt{1 + u^2}\ \frac{4 + u^2}{3} + \frac{2 + a^2}{3}  + \ln
\frac{\sqrt{1 + u^2} - 1}{u} - \ln \frac{a^2 - 1}{a} & \text{if }
u \ge a,
\end{array}
\right.
\end{equation*}
and let $U$ be a (unique) solution of (\ref{eql3b}). Define the
function $f_h$ like in the statement (b) of lemma \ref{l3b}. Then
the body of least resistance is
\begin{equation*}
\{ (x\text{'}, x_3) \in \mathbb{R}^3 :\, |x\text{'}| \le 1,\ |x_3|
\le -f_h(|x\text{'}|)\}\,,
\end{equation*}
where $x\text{'} = (x_1, x_2)$. Thus, the body is symmetric with
respect to the horizontal plane $\{ x_3 = 0 \}$. The front and
rear parts of its surface contain flat circular disks of equal
radius, and the angle of slope of lateral surface near these disks
equals $\arctan a \approx 51.8^0$.

\begin{figure}
\begin{center}
  \includegraphics[scale=0.3]{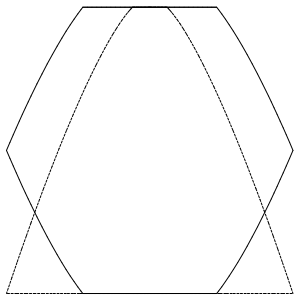}
\end{center}
\caption{Three-dimensional case, $h=1$. The solution in the limit
$V \to 0^+$ (solid line) and the solution of Newton's classical
problem (dashed line).} \label{fig:3d-cl-gNh1-gSh1}
\end{figure}

On Fig.~\ref{fig:3d-cl-gNh1-gSh1}, the projections of two optimal bodies of height
$h = 1$ on the plane $Ox_1x_3$ are shown. The bodies are: the body of least resistance
in the limit $V \to 0^+$ and the solution of Newton's classical problem.

%

The minimal reduced resistance equals
\begin{equation*}
\hat{\mathrm R}(h,V) = 2 c^{(3)} \left( \bar p(U) +
\frac{Q(U)}{q^2(U)}\right),
\end{equation*}
where $\bar p$ is defined in (\ref{onemore}), and $U$ is given by
(\ref{eql3b}).

It is interesting to note that in these limit cases, the shape
of optimal body does not depend on the distribution $\s$; moreover,
the reduced minimal resistance is proportional to $\nu$ in the limit
$V \to +\infty$, and is proportional to the factor $c^{(d)}$
given by (\ref{c1c2}),\,(\ref{k1k2}) in the limit $V \to 0^+$.
This factor can be interpreted as the sum of absolute values of
impulses of particles of the medium in unit volume, in the frame
of reference associated with the medium.

\section{Gaussian distribution of velocities: exact solutions}
\label{sec:GaussDV}

Suppose that the function $\rho = \rho_V$ is the density of
circular gaussian distribution with mean $-Ve_d$ and variance 1,
\textrm{i.e.},
\begin{equation}\label{gauss}
\rho_V(v) = \s(|v + V e_d|), \ \ \text{ where } \ \ \s(r) =
(2\pi)^{-d/2}\, e^{-r^2/2}.
\end{equation}
This function describes the particles' distribution over
velocities in a frame of reference moving in a homogeneous
monatomic ideal gas, where the velocity of motion equals $V$ and
the mean square velocity of molecules equals 1 (see example~\ref{e1}).

The function $\rho_V$ satisfies the condition A, therefore the
results obtained in the previous section can be applied in this case.
Here, the pressure functions $p_\pm(u, V)$ are calculated analytically
in the cases $d = 2$ and $d = 3$, and then, using numerical simulation,
the following results are obtained:

1) The parameter set $V$-$h$ is divided into several
subsets corresponding to different kinds of solutions.
This partition is shown on Figures\,\ref{fig:vh} and
\ref{fig:3DheightversusVel}.

2) The least resistance $\mathrm R(h, V)$ is calculated for
various values of $h$,\, $V$. The results are shown on
Figures\,\ref{resistance2d} and \ref{fig:RRh3d}.

3) For several values of parameters $h$ and $V$, the body of least
resistance is constructed. Two such bodies are shown on
Figures\,\ref{fig:sol3DNLSd3},\, \ref{fig:sol2D:fk:sk}\,(a)--(d)
and \ref{fig:NLSd3:wap}.

Here the value $V$ is allowed to vary, so the pressure and
resistance functions are designated by $p_\pm(u, V)$ and $\mathrm{R}(h,
V)$ instead of $p_\pm(u)$ and $\mathrm{R}(h)$.

Consider the cases $d = 2$ and $d = 3$.

\subsection{Two-dimensional case}

Fixing the sign ``+'' and passing to polar coordinates $v = (-r
\sin\vphi, -r \cos\vphi)$ in the formula (\ref{formula}), one
obtains
\begin{equation}\label{rhod=211}
p_+(u, V) = \int\!\!\int \frac{r^2 (\cos\vphi +
u\sin\vphi)_{\!+}^{\,\ 2}}{1 + u^2}\, \rho_+ (r, \vphi, V)\, r dr
d\vphi,
\end{equation}
where $z_+ := \max \{ 0,\, z \}$, and $\rho_+(r, \vphi, V)$ is the
gaussian density function $\rho_V$ (\ref{gauss}) written in the
introduced polar coordinates,
\begin{equation}\label{rhod=222}
\rho_+(r, \vphi, V) = \frac{1}{2\pi}\, e^{-\frac 12 (r^2 - 2Vr
\cos\vphi + V^2)}.
\end{equation}
Next, fixing the sign ``$-$'' and introducing polar coordinates in a
slightly different manner, $v = (-r \sin\vphi, r \cos\vphi)$, one
obtains
\begin{equation}\label{rhod=233}
p_-(u, V) = -\int\!\!\int \frac{r^2 (\cos\vphi +
u\sin\vphi)_{\!+}^{\,\ 2}}{1 + u^2}\, \rho_- (r, \vphi, V)\, r dr
d\vphi,
\end{equation}
Here $\rho_-(r, \vphi, V)$ is the same density function $\rho_V$
(\ref{gauss}) written in these coordinates,
\begin{equation}\label{rhod=244}
\rho_-(r, \vphi, V) = \frac{1}{2\pi}\, e^{-\frac 12 (r^2 + 2rV
\cos\vphi + V^2)}.
\end{equation}
Combining the formulas (\ref{rhod=211}), (\ref{rhod=222}),
(\ref{rhod=233}), and (\ref{rhod=244}), one comes to the more
general expression
\begin{equation*}
p_\ve(u, V) = \ve\, \frac{e^{-V^2/ 2}}{2\pi}\ \ \
\int\!\!\hspace{-8mm}\int\limits_{\cos\varphi + u\sin\varphi > 0}
\hspace{-4mm} \frac{\left(\cos\varphi + u\sin\varphi\right)^2}{1 +
u^2}\ e^{-\frac 12 r^2 + \ve 2rV \cos\varphi}\ r^3\, dr d\varphi
\,,
\end{equation*}
where $\ve \in \{ -,\, + \}$. Passing to the iterated integral and
integrating over $r$, one obtains
\begin{equation}\label{asympt1}
p_\ve(u, V) = \ve\, \frac{e^{-V^2/2}}{\pi} \hspace{-2mm}
\int\limits_{\cos\varphi + u\sin\varphi > 0} \hspace{-3mm}
\frac{(\cos\varphi + u\sin\varphi)^2}{1 + u^2}\ l(\ve V
\cos\varphi)\, d\varphi \, ,
\end{equation}
where
$$
l(z) = 1 + \frac{z^2}{2} + \frac{\sqrt\pi}{2\sqrt 2}\ e^{z^2/ 2}
\left(3z + z^3\right) \left( 1 + \text{erf} \left( {z}/{\sqrt 2}
\right) \right) ,
$$
and $\text{erf}(x) = \frac{2}{\sqrt{\pi}}\,\int _{0}^{x}\!{e^{-{t}^{2}}}{dt}$.
Changing the variable $\tau = \vphi - \arcsin(u/ \sqrt{1+u^2})$,
one finally comes to
\begin{equation}\label{asympt2}
p_\ve(u, V) = \ve \frac{e^{-V^2/ 2}}{\pi} \int_{-\pi/2}^{\pi/2}
\cos^2 \tau\ l \left( \ve V\, \frac{\cos\tau - u\sin\tau}{\sqrt{1
+ u^2}} \right) d\tau.
\end{equation}

Numerical simulations are done using Maple and verified by Matlab.

\begin{figure}
\begin{center}
\psline(3.25,2.2)(3,1.8)(3.5,1.8)(3.25,2.2) 
\psline(3.25,3.78)(3,3.28)(3.5,3.28)(3.25,3.78) 
\psline(3,3.28)(3.1,3.13)(3.4,3.13)(3.5,3.28) 
\psline(3.25,7.5)(3,6.5)(3.5,6.5)(3.25,7.5) 
\psline(3,6.5)(3.25,6)(3.5,6.5) 
\psline(3.1,1.15)(3,1.04)(3.5,1.04)(3.4,1.15)(3.1,1.15) 
\psfrag{v}{$V$} \psfrag{h}{$h$} \psfrag{0}{$0$} \psfrag{1}{$1$}
\psfrag{2}{$2$} \psfrag{3}{$3$} \psfrag{4}{$4$} \psfrag{5}{$5$}
\psfrag{10}{$10$} \psfrag{15}{$15$} \psfrag{20}{$20$}
\psfrag{25}{$25$}
\includegraphics[width=8cm,height=10cm]{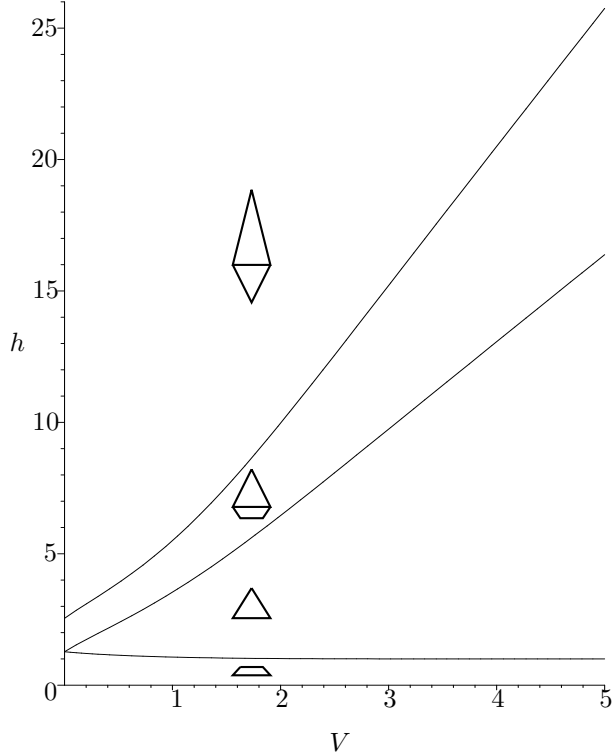}
\end{center}
\caption{Two-dimensional case. Four regions shown on the parameter
space correspond to four kinds of solutions.
 }
 \label{fig:vh}
\end{figure}
Graphs of the functions
$$
h = u_+^0(V),\ \ \ h = u_*(V),\ \ \ h = u_*(V) + u_-^0(V),
$$
are shown on Figure~\ref{fig:vh}, where the values $u_+^0$,\,
$u_-^0$ and $u_*$ (which are the functions of $V$) are defined in
subsections 4.1.1, 4.1.2 and 4.1.3, respectively. These graphs
separate the parameter space $\mathbb{R}_+^2$ into four regions
corresponding to the four different kinds of solutions. The lower
function $h = u_+^0(V)$ tends to 1 as $V \to \infty$. Further, at
$V = 0$, the lower, the middle, and the upper functions take the
values $a$,\, $a$,\, and $2a$, respectively, where $a = \sqrt{(1 +
\sqrt 5)/ 2} \approx 1.272$:~ $\lim_{V\to 0} u_+^0(V) = a$,\,
$\lim_{V\to 0} u_*(V) = a$,\, $\lim_{V\to 0} (u_*(V) + u_-^0(V)) =
2a$.

The solutions of fifth kind (union of two triangles and a
trapezium) were not found in numerical simulations. (These
solutions correspond to the case where the set
$\mathcal{O}_{p_-,V} = \{ u : \bar{p}_-(u,V) < p_-(u,V) \}$
contains at least two connected components.) We believe that in
the considered case corresponding to the gaussian distribution
$\rho_V$ this kind of solutions does not appear at all. (Notice
that, according to the statement (b) of lemma~\ref{l2}, this kind
of solutions does appear for some distributions corresponding to
mixtures of homogeneous gases.)

\begin{figure}
\subfigure[]{
  \label{fig:RRh1}
  \begin{minipage}[b]{0.5\textwidth}
  \begin{flushleft}
    \hspace*{-0.6cm}
    \psfrag{h}{$h$}
    \psfrag{RR}{$\hspace*{-7mm}\tilde{\mathrm R}(V,h)$}
    \psfrag{0}{$0$}
    \psfrag{0.5}{$0.5$}
    \psfrag{1}{$1$}
    \psfrag{1.5}{$1.5$}
    \psfrag{2}{$2$}
    \psfrag{2.5}{$2.5$}
    \psfrag{3}{$3$}
    \psfrag{3.5}{$3.5$}
    \psfrag{2}{$\hspace*{-1mm}2$}
    \psfrag{4}{$\hspace*{-1mm}4$}
    \psfrag{6}{$\hspace*{-1mm}6$}
    \psfrag{8}{$\hspace*{-1mm}8$}
    \psfrag{10}{$\hspace*{-1mm}10$}
    \psfrag{12}{$\hspace*{-1mm}12$}
    \psfrag{14}{$\hspace*{-1mm}14$}
    \psfrag{16}{$\hspace*{-1mm}16$}
    \psfrag{v=0.1}{$V = 0.1$}
    \psfrag{v=0.2}{$V = 0.2$}
    \psfrag{v=0.5}{$V = 0.5$}
    \psfrag{v=1}{$V = 1$}
    \psfrag{v=infinity}{$V = \infty$}
    \includegraphics[width=6.5cm,height=10cm]{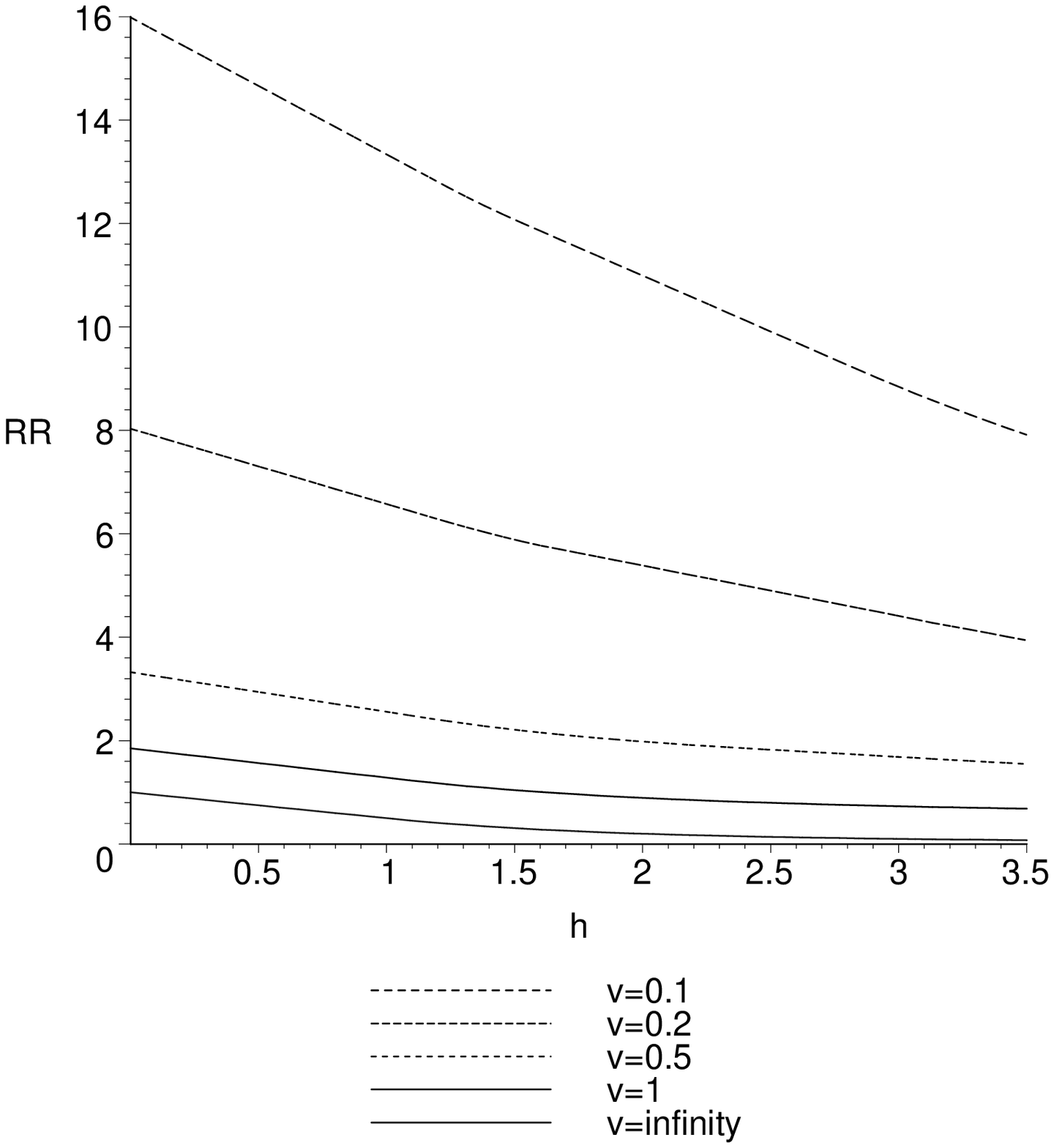}
   \end{flushleft}
  \end{minipage}}
\subfigure[]{
  \label{fig:RRh2}
  \begin{minipage}[b]{0.5\textwidth}
    \begin{flushright}
    \hspace*{0.6cm}
    \psfrag{h}{$h$}
    \psfrag{RR}{$\hspace*{-7mm}\tilde{\mathrm R}(V,h)$}
    \psfrag{0}{$0$}
    \psfrag{1}{$\hspace*{-1mm}1$}
    \psfrag{2}{$2$}
    \psfrag{3}{$3$}
    \psfrag{4}{$4$}
    \psfrag{5}{$5$}
    \psfrag{6}{$6$}
    \psfrag{0.2}{$\hspace*{-1mm}0.2$}
    \psfrag{0.4}{$\hspace*{-1mm}0.4$}
    \psfrag{0.6}{$\hspace*{-1mm}0.6$}
    \psfrag{0.8}{$\hspace*{-1mm}0.8$}
    \psfrag{1.2}{$\hspace*{-1mm}1.2$}
    \psfrag{v=2}{$V = 2$}
    \psfrag{v=3}{$V = 3$}
    \psfrag{v=5}{$V = 5$}
    \psfrag{v=infinity}{$V = \infty$}
    \includegraphics[width=6.5cm,height=10cm]{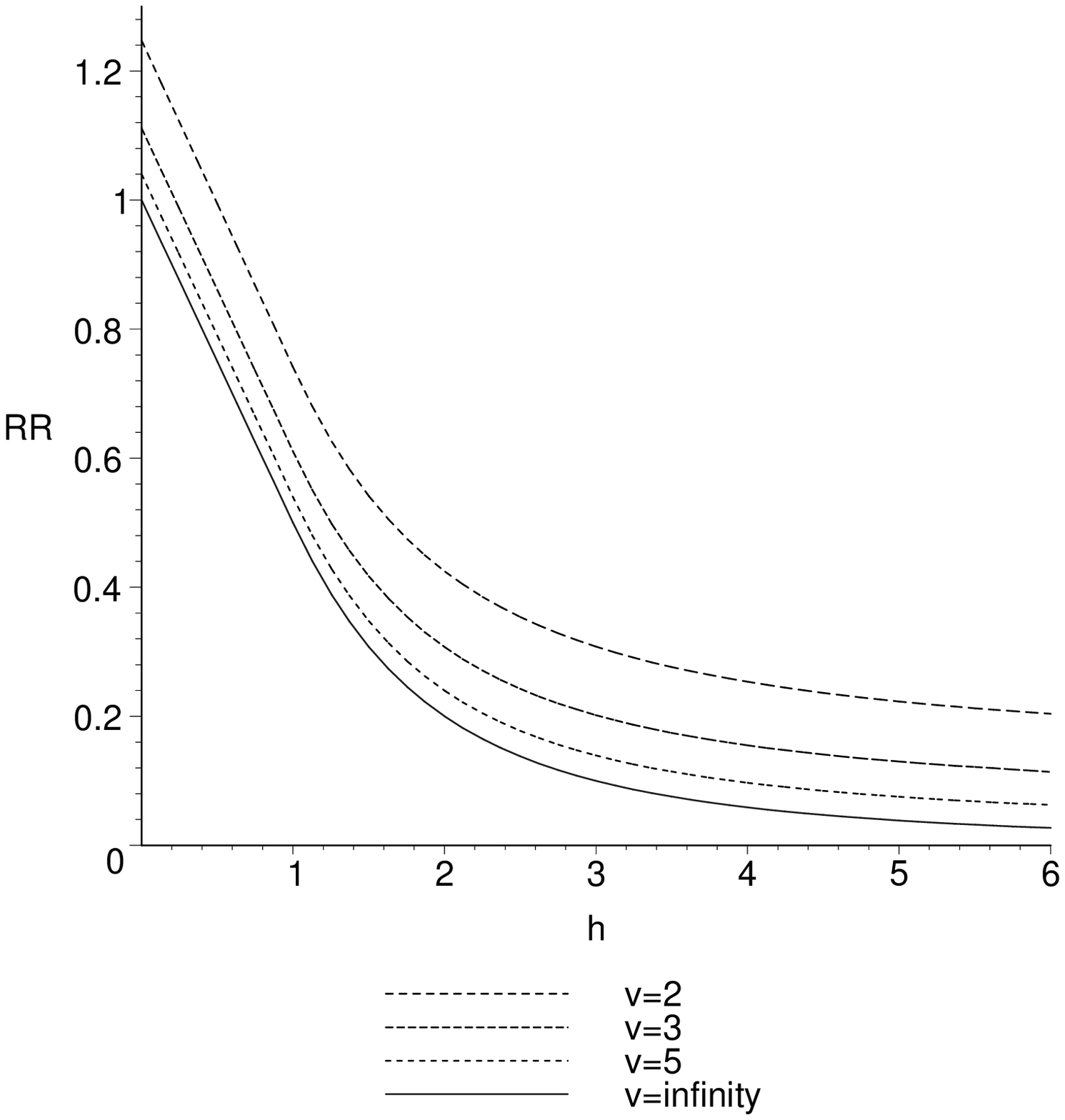}
    \end{flushright}
  \end{minipage}}
\caption{Two-dimensional case. Least reduced resistance
$\tilde{\mathrm R}(V,h)$ versus height $h$ of the body}
\label{resistance2d}
\end{figure}

Further, using the formulas from subsections 4.1.1 and 4.1.2, the
functions $f_+$ and $f_-$ are calculated, which allow one to
construct the optimal figures
(Fig.\,\ref{fig:sol2D:fk:sk}\,(a)--(d)), and using the formulas
from 4.1.3, the minimal resistance $\mathrm{R}(h,V)$ is
calculated. The graphs of reduced minimal resistance
$\tilde{\mathrm R}(h,V) = V^{-2} \mathrm R(h,V)$ versus $h$ are
shown on figures \ref{fig:RRh1} and \ref{fig:RRh2} for several
values of $V$.


\subsection{Three-dimensional case}

Fix the sign ``+''. In spherical coordinates $v = (-r \sin\vphi
\cos\theta,\, -r \sin\vphi \sin\theta,\, -r \cos\vphi)$,\, $r \ge
0$,\, $0 \le \vphi \le \pi$,\, $-\pi \le \theta \le \pi$ the
formula (\ref{formula}) takes the form
$$
p_+(u, V) = \int\!\!\!\int\!\!\!\int \frac{r^2 (\cos\vphi +
u\sin\vphi \cos\theta)_{\!+}^{\,\ 2}}{1 + u^2}\, \rho_+ (r, \vphi,
\theta, V)\, r^2 \sin\vphi\, dr d\vphi d\theta,
$$
where
$$
\rho_+(r, \vphi, \theta, V) = \frac{1}{(2\pi)^{3/2}}\, e^{-\frac
12 (r^2 - 2Vr \cos\vphi + V^2)}.
$$
Now, fix the sign ``$-$''. In spherical coordinates $v = (-r
\sin\vphi \cos\theta,\, -r \sin\vphi \sin\theta,\, r \cos\vphi)$,
one has
$$
p_-(u, V) = -\int\!\!\!\int\!\!\!\int \frac{r^2 (\cos\vphi +
u\sin\vphi \cos\theta)_{\!+}^{\,\ 2}}{1 + u^2}\, \rho_- (r, \vphi,
\theta, V)\, r^2 \sin\vphi\, dr d\vphi d\theta,
$$
where
$$
\rho_-(r, \vphi, \theta, V) = \frac{1}{(2\pi)^{3/2}}\, e^{-\frac
12 (r^2 + 2Vr \cos\vphi + V^2)}.
$$
Summarizing, one comes to the formula
$$
p_\ve(u, V) = \ve\, \frac{e^{-V^2/2}}{(2\pi)^{3/2}}\ \int
\!\int\!\!\!\hspace{-10mm}\int\limits_{\cos\vphi +
u\sin\vphi\cos\theta
> 0} \frac{(\cos\vphi + u\sin\vphi\ cos\theta)^2}{1 + u^2}\, \cdot
$$
\begin{equation}\label{ppmu}
\cdot\, e^{-\frac 12 r^2 + \ve V r \cos\vphi} r^4 \sin\vphi\, dr
d\vphi d\theta, \ \ \ \ \ve \in \{ -,\, + \}.
\end{equation}
We shall use two formulas, which are easy to verify. First,
$$
\int_0^{+\infty} e^{-\frac 12 r^2 + \ve V r \cos\vphi} r^4\, dr =
I(\ve V\cos\vphi),
$$
where
$$
I(z) = \sqrt{\pi/2}\ e^{z^2/2}\, (3 + 6z^2 + z^4) (1 +
\text{erf}\, (z/\sqrt 2)) + 5z + z^3.
$$
Second,
$$
\int_{\cos\vphi + u\sin\vphi\cos\theta > 0} \frac{(\cos\vphi +
u\sin\vphi\cos\theta)^2}{1 + u^2}\, d\theta\, =\, J(u,\cos\vphi),
$$
where
$$
J(u,\zeta) = \left\{ \begin{array}{ll} 0, \ \ \ & \text{ if } \
-1 \le \zeta \le -u/ \sqrt{1 + u^2}\\
J_1(u, \zeta), \ \ \ & \text{ if } \
|\zeta| < u/ \sqrt{1 + u^2}\\
J_2(u, \zeta), \ \ \ & \text{ if } \ u/ \sqrt{1 + u^2} \le \zeta
\le 1,
\end{array} \right.
$$
and
$$
J_1(u, \zeta) = \frac{1}{1 + u^2} \left[ \theta_0\, (2\zeta^2 +
u^2(1 - \zeta^2)) + 3\zeta \sqrt{u^2 - \zeta^2(1 + u^2)} \right],
$$
$$
J_2(u, \zeta) = \frac{\pi}{1 + u^2}\, [ 2\zeta^2 + u^2 (1 -
\zeta^2)], \ \ \ \ \ \theta_0 = \arccos \left(
 -\frac{\zeta}{u \sqrt{1 - \zeta^2}} \right).
$$
Taking into account these formulas and changing the variable
$\zeta = \cos\vphi$ in (\ref{ppmu}), one gets
$$
p_\ve(u, V) = \ve\, \frac{e^{-V^2/2}}{(2\pi)^{3/2}} \int_{-1}^1\,
I(\ve V\zeta)\, J(u,\zeta)\, d\zeta =
$$
$$
= \ve\, \frac{e^{-V^2/2}}{(2\pi)^{3/2}} \left( \int_{-u/\sqrt{1 +
u^2}}^{u/\sqrt{1 + u^2}} I(\ve V\zeta)\, J_1(u,\zeta)\, d\zeta +
\int_{u/\sqrt{1 + u^2}}^1 I(\ve V\zeta)\, J_2(u,\zeta)\, d\zeta
\right).
$$

\begin{figure}
\subfigure[$V = 1$, $h = 1.97$]{ 
  \label{fig:wa-solV1Up321}
  \begin{minipage}[b]{0.5\textwidth}
  \centering
  \includegraphics[scale=0.4]{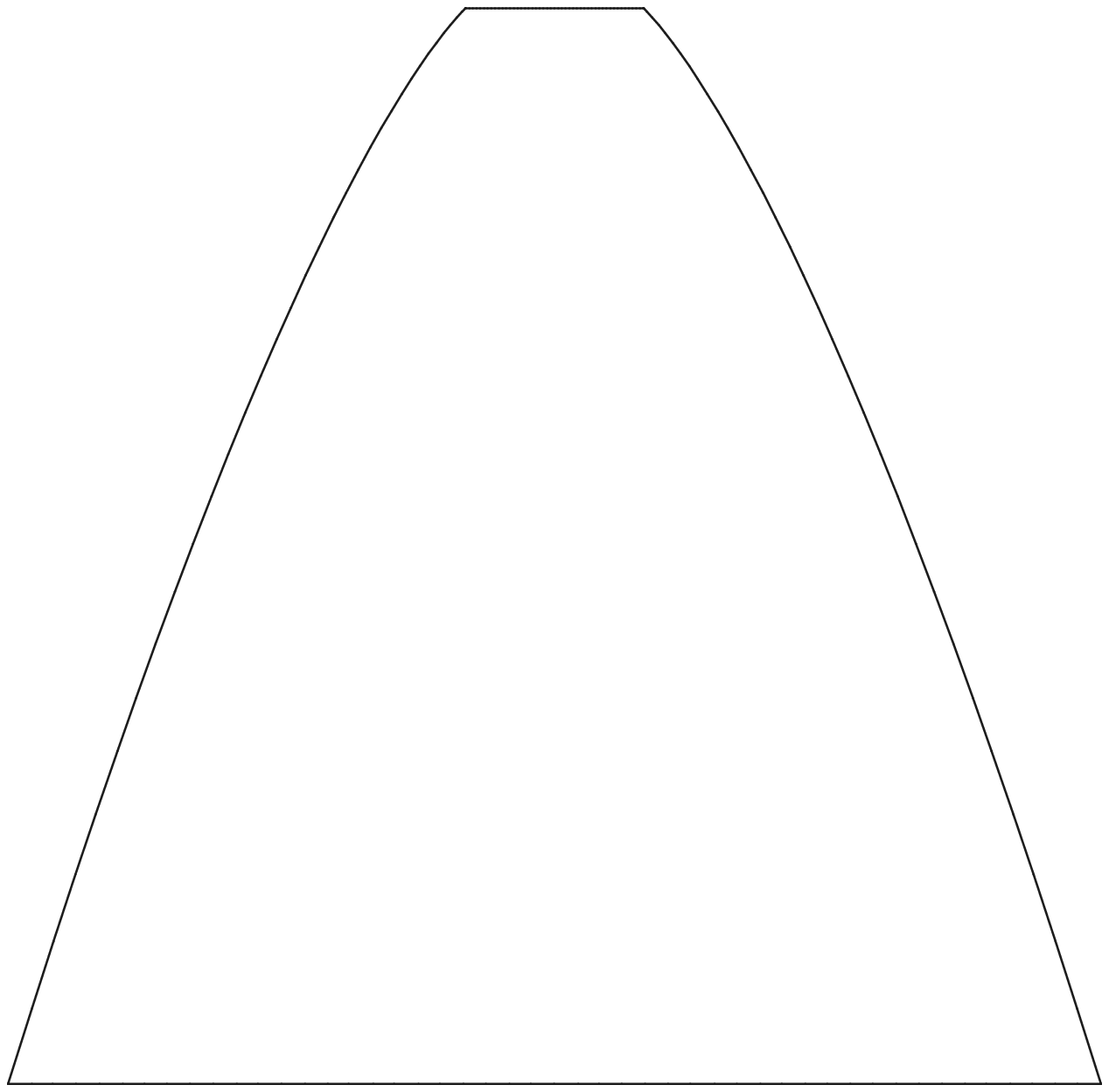}
  \end{minipage}}
\subfigure[$V = 1$, $h = 3.11$]{ 
  \label{fig:wa-solV1Up401}
  \begin{minipage}[b]{0.5\textwidth}
  \centering
  \includegraphics[scale=0.4]{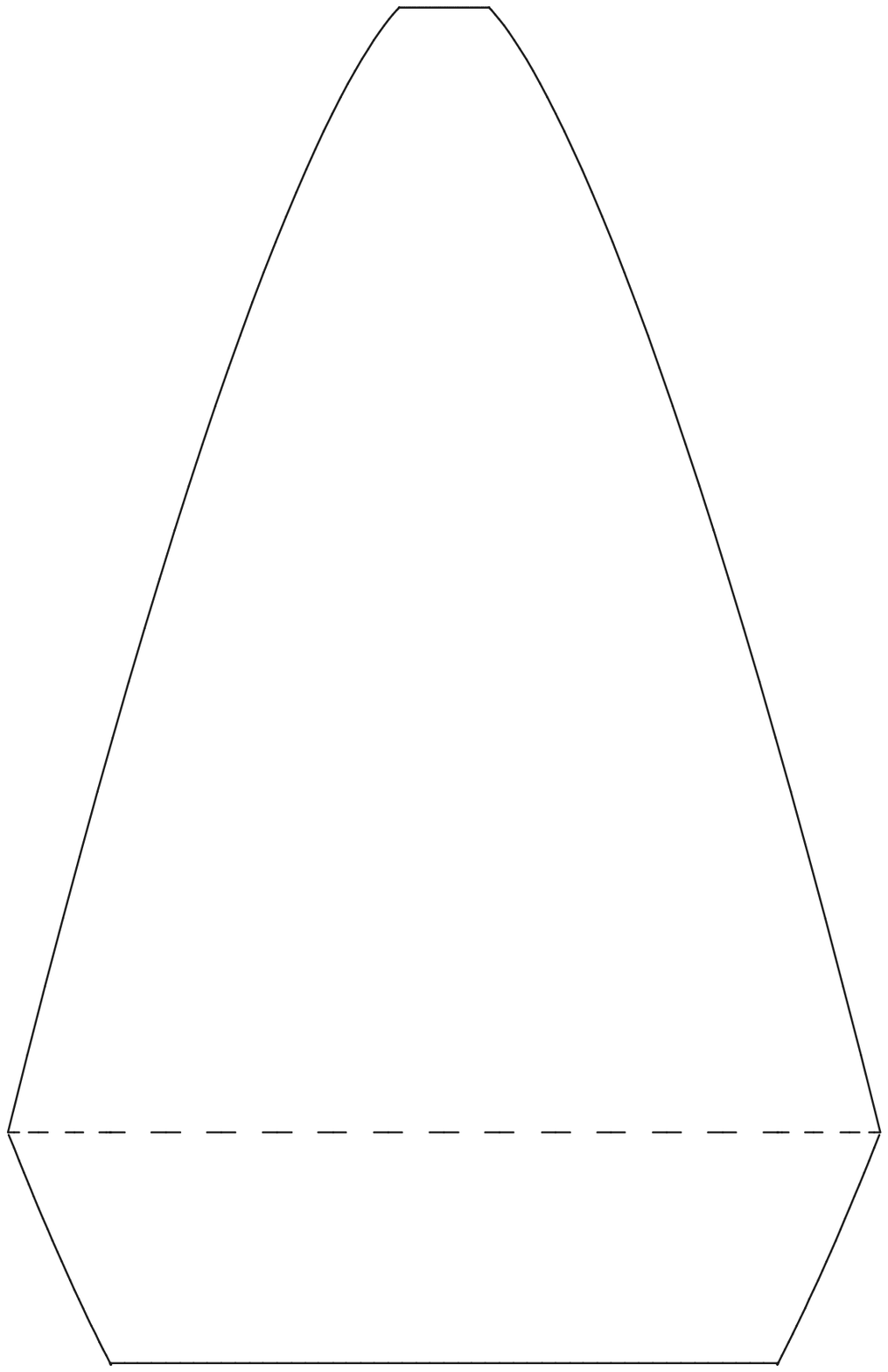}
  \end{minipage}}
\caption{The solutions in three-dimensional case, for the
distribution function $\rho_V$ (\ref{gauss}).}
\label{fig:NLSd3:wap}
\end{figure}

\begin{figure}
  \subfigure[$h_*(V)$]{
    \label{fig:gauss3Dvh}
    \begin{minipage}[b]{0.5\textwidth}
      \begin{center}
        \psfrag{V}{$V$}
        \psfrag{h}{$h$}
        \pscurve[linewidth=0.1pt]{-}(2.5,2)(2.58,2.25)(2.71,2.5)
        \psline[linewidth=0.1pt,linearc=0.0]{-}(2.71,2.5)(2.80,2.5)
        \pscurve[linewidth=0.1pt]{-}(2.80,2.5)(2.93,2.25)(3.01,2)
        \psline[linewidth=0.1pt,linearc=0.0]{-}(3.01,2)(2.5,2) 
        \pscurve[linewidth=0.1pt]{-}(2,4.5)(2.092,4.83)(2.204,5.08)
        \psline[linewidth=0.1pt,linearc=0.0]{-}(2.204,5.08)(2.25,5.08)
        \pscurve[linewidth=0.1pt]{-}(2.25,5.08)(2.363,4.83)(2.455,4.5)
        \psline[linewidth=0.1pt,linearc=0.0,linestyle=dashed,dash=1pt 1pt]{-}(2.455,4.5)(2,4.5)
        \psline[linewidth=0.1pt,linearc=0.0]{-}(2,4.5)(2.053,4.38)(2.402,4.38)(2.455,4.5)
        \includegraphics[scale=0.4]{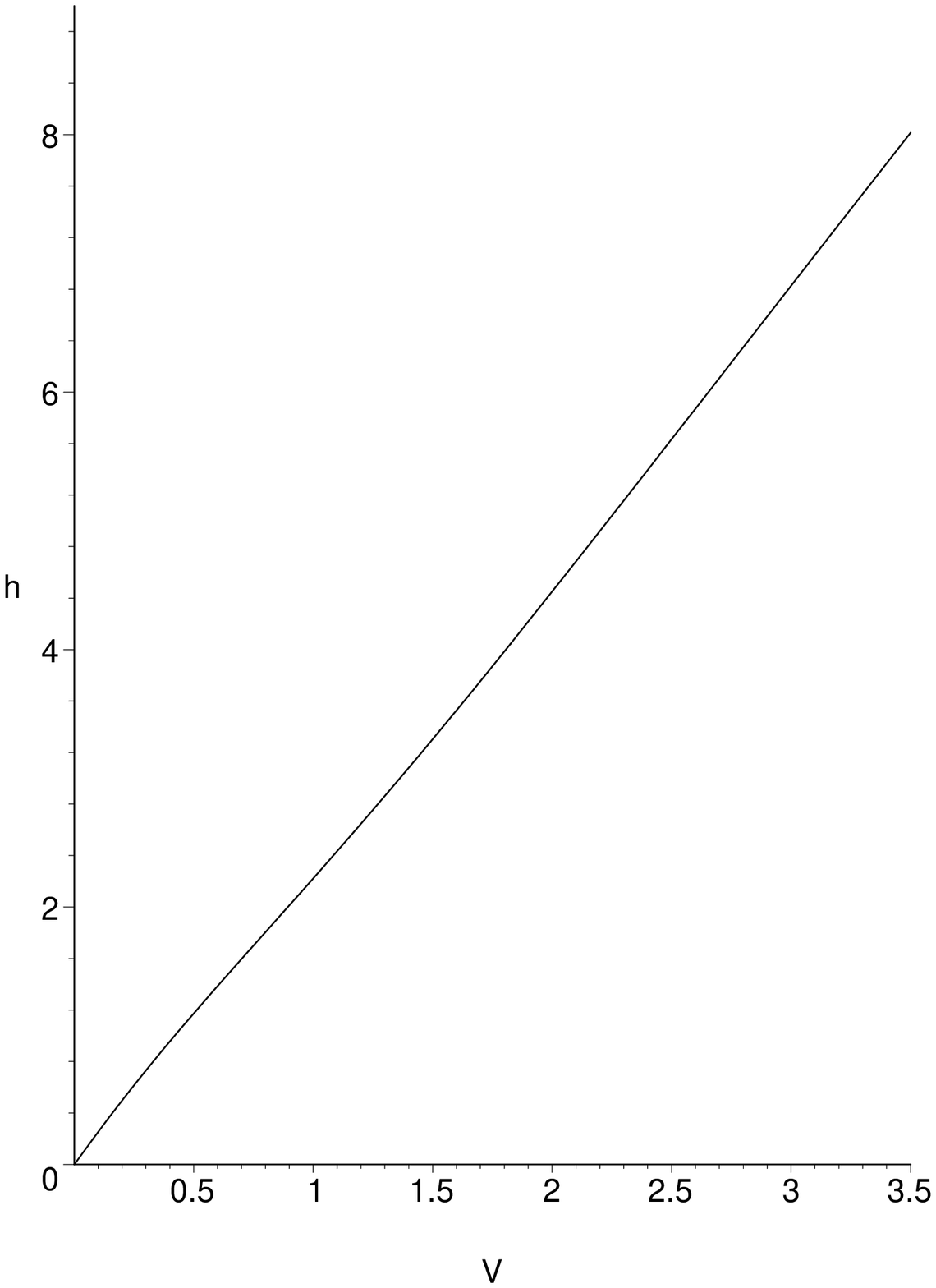}
      \end{center}
    \end{minipage}}
  \subfigure[${dh_*}/{dV}$]{
    \label{fig:gauss3Ddvdh}
    \begin{minipage}[b]{0.5\textwidth}
    \begin{center}
      \psfrag{V}{$V$}
      \psfrag{hprime}{$\hspace*{-3mm}h_*'(V)$}
      \includegraphics[scale=0.4]{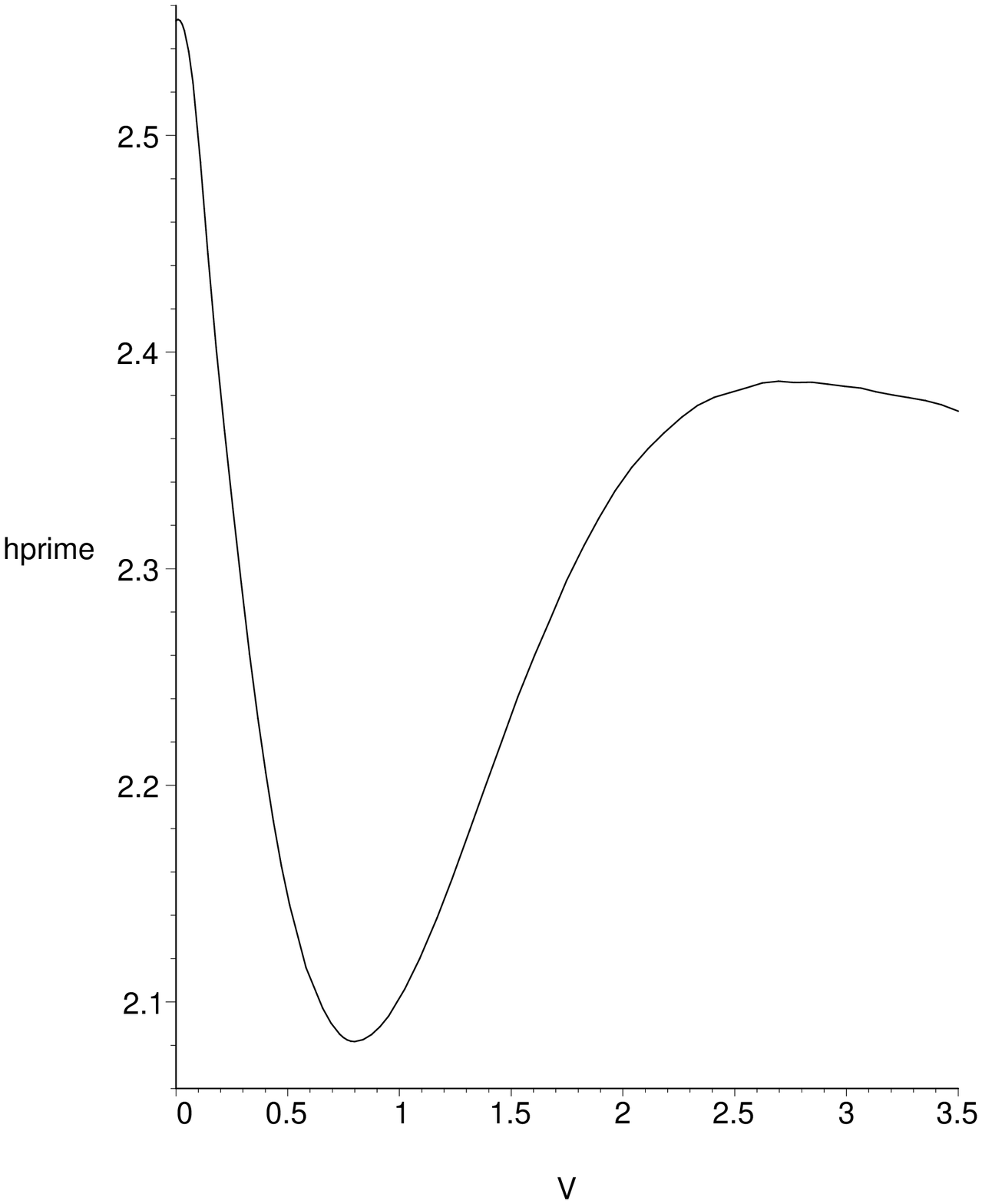}
    \end{center}
    \end{minipage}}
\caption{Three-dimensional case. The function $h_*(V)$ divides the
parameter plane $V$-$h$ into two subsets corresponding to the two
kinds of solutions.} \label{fig:3DheightversusVel}
\end{figure}

\begin{figure}
\begin{center}
\psfrag{h}{$h$}
\psfrag{R}{$\hspace*{-10mm}\tilde{\mathrm R}(h,V)$}
\psfrag{v=0.5}{$V = 0.5$}
\psfrag{v=1}{$V = 1$}
\includegraphics[width=6.5cm,height=10cm]{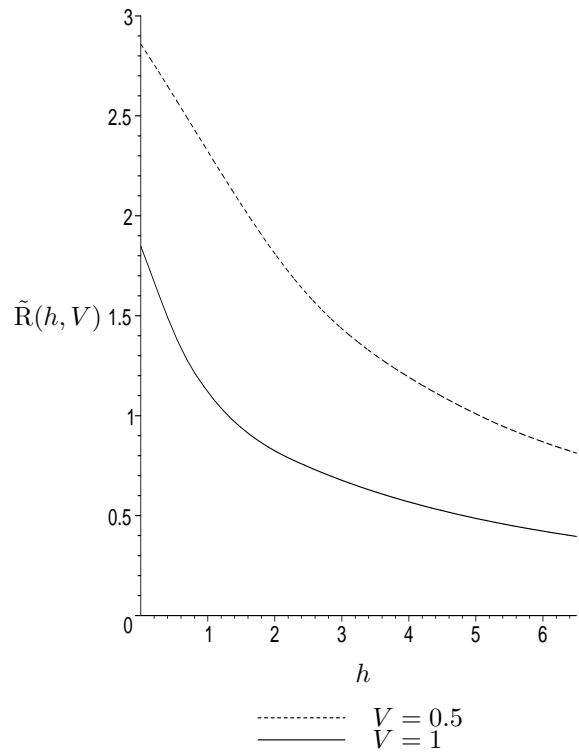}
\end{center}
\caption{Three-dimensional case. Minimal reduced resistance
$\tilde{\mathrm{R}}(V,h)$ versus height $h$ of the body.}
\label{fig:RRh3d}
\end{figure}

Next, one numerically calculates the function $h_*(V)$ according
to the formula (\ref{hstar}). This function is shown on
Fig.\,\ref{fig:gauss3Dvh}; it divides the parameter set
$\mathbb{R}_+^2$ into two subsets corresponding to two different
kinds of solutions. The function $h_*(V)$ looks like linear, but
is not; the graph of its derivative is shown on
Fig.\,\ref{fig:gauss3Ddvdh}. The function $\mathrm{R}(h,V)$ is
calculated according to the formulas given in subsection 4.2; the
graphs of $\tilde{\mathrm{R}}(h,V)$ versus $h$ are shown on
Fig.~\ref{fig:RRh3d} for several values of $V$. On
Fig.~\ref{fig:NLSd3:wap} the examples of solutions of the first
and the second kind are presented, for parameters indicated there.

\appendix

\section*{Appendix A}

\renewcommand{\theequation}{A.\arabic{equation}}
\setcounter{equation}{0}

\subsection*{Proof of lemma \ref{l1}}

Changing the variable $v = r\nu$,\, $r \ge 0$,\, $\nu \in S^{d-1}$
in the integral (\ref{formula}) and using condition A, one obtains
$$
p_\ve(u) = \ve \int_{S^{d-1}} d\HHH^{d-1}(\nu)\, \int_0^{+\infty}
r^2\ \frac{(\nu_1 u + \ve \nu_d)\!_-^{\,\ 2}}{1 + u^2}\
\rho(r\nu)\, r^{d-1} dr =
$$
\begin{equation*}
= \ve \int_{S^{d-1}} \frac{(\nu_1 u + \ve \nu_d)\!_-^{\,\ 2}}{1 +
u^2}\ \bar\rho(\nu)\ d\HHH^{d-1}(\nu),
\end{equation*}
where
\begin{equation*}
\bar\rho(\nu) := \int_0^{+\infty} r^2\, \rho(r\nu)\, r^{d-1} dr =
\int_0^{+\infty} r^2\, \sigma(\sqrt{r^2 + 2rV\, \nu_d + V^2})\,
r^{d-1} dr.
\end{equation*}
Substituting $u = \tan\vphi$, $\vphi \in [0,\, \pi/2]$, one
obtains
$$
p_\ve(\tan\vphi) = \ve \int_{S^{d-1}} (\nu_1 \sin\vphi + \ve \nu_d
\cos\vphi)\!_-^{\,\ 2}\, \bar\rho(\nu)\, d\HHH^{d-1}(\nu).
$$

Substitute ``+'' for $\ve$ and consider the rotation $T_\vphi$ that
sends the vector $(\sin\vphi,\, 0, \ldots, 0,\, \cos\vphi)$ to
$e_d$ and leaves the vectors $e_i$,\, $i = 2, \ldots, d-1 \,$
unchanged. For any $\nu \in \RRR^d$ one has $T_\vphi \nu =
(\cos\vphi\, \nu_1 - \sin\vphi\, \nu_d,\ \nu_2, \ldots,
\nu_{d-1},\ \sin\vphi\, \nu_1 + \cos\vphi\, \nu_d)$.~ Changing the
variable $T_\vphi \nu = \om$, one gets
\begin{equation}\label{square}
p_+(\tan\vphi) = \int_{S_-^{d-1}} \om_d^2\ \bar\rho (T_\vphi^{-1}
\om)\ d\HHH^{d-1}(\om),
\end{equation}
where $S_-^{d-1} := \{ \om \in S^{d-1} : \om_d < 0 \}$. Designate
\begin{equation}\label{triangle}
\varrho(z) := \int_0^\infty r^2\, \sigma(\sqrt{r^2 + 2rV z +
V^2})\, r^{d-1} dr, \ \ \ \ |z| \le 1;
\end{equation}
obviously, $\bar\rho(\nu) = \varrho(\nu_d)$. Using condition A,
one concludes that the function $\varrho$ is continuously
differentiable, and its derivative
$$
\varrho'(z) = \int_0^\infty r^2\, \frac{\sigma'(\sqrt{r^2 + 2rV z
+ V^2})} {\sqrt{r^2 + 2rV z + V^2}}\, rV\, r^{d-1} dr
$$
is negative and monotone non-decreasing;~ in particular,
\begin{equation}\label{circle}
\text{as } \ \ z > 0, \ \ \ \ \ \varrho'(z) > \varrho'(-z).
\end{equation}
Using that
$$
T_\vphi^{-1} \om = T_{-\vphi} \om = (\cos\vphi\, \om_1 +
\sin\vphi\, \om_d,\ \om_2, \ldots, \om_{d-1},\ -\sin\vphi\, \om_1
+ \cos\vphi\, \om_d),
$$
from (\ref{square}) and (\ref{triangle}) one obtains
\begin{equation}\label{p+tg}
p_+(\tan\vphi) = \int_{S_-^{d-1}} \om_d^2\ \varrho (-\sin\vphi\,
\om_1 + \cos\vphi\, \om_d)\ d\HHH^{d-1}(\om).
\end{equation}

Now, substitute ``$-$'' for $\ve$ and consider the orthogonal
reflection $U_\vphi$ with respect to the hyperplane $\{
\sin\frac{\vphi}{2}\,~ \om_1 = \cos\frac{\vphi}{2}\,~ \om_d \}$;~
for any $\nu \in \RRR^d$ one has
$$
U_\vphi \nu = (\cos\vphi\, \nu_1 + \sin\vphi\, \nu_d,\ \nu_2,
\ldots, \nu_{d-1},\ \sin\vphi\, \nu_1 - \cos\vphi\, \nu_d).
$$
Changing the variable $U_\vphi \nu = \om$, one gets
$$
p_-(\tan\vphi) = -\int_{S_-^{d-1}} \om_d^2\ \bar\rho (U_\vphi^{-1}
\om)\ d\HHH^{d-1}(\om).
$$
Using that $U_\vphi^{-1} = U_\vphi$, one obtains
\begin{equation}\label{p-tg}
p_-(\tan\vphi) = -\int_{S_-^{d-1}} \om_d^2\ \varrho (\sin\vphi\,
\om_1 - \cos\vphi\, \om_d)\ d\HHH^{d-1}(\om).
\end{equation}

The formulas (\ref{p+tg}) and (\ref{p-tg}) can be written in the
unified form
\begin{equation}\label{pve}
p_\ve(\tan\vphi) = \ve \int_{S_-^{d-1}} \om_d^2\ \varrho ( \ve
(-\sin\vphi\, \om_1 + \cos\vphi\, \om_d))\ d\HHH^{d-1}(\om).
\end{equation}
Substituting $\vphi = \pi/2$ in (\ref{pve}), one obtains
$$
\lim_{u\to +\infty} p_\ve(u) = p_\ve(+\infty) = \ve
\int_{S_-^{d-1}} \om_d^2\, \varrho (- \ve \om_1)\,
d\HHH^{d-1}(\om).
$$
Using that $S_-^{d-1}$ is invariant with respect to reflection
$\om_1 \mapsto -\om_1$, one concludes that $p_+(+\infty) = -
p_-(+\infty)$, so (a) is proved.

Further, using (\ref{pve}), one concludes that the function
$p_\ve$ is continuously differentiable, and
$$
p_\ve'(\tan\vphi) = \ve \cos^2 \vphi \cdot \int_{S_-^{d-1}}
\om_d^2\ \frac{\pl}{\pl\vphi}\, \varrho ( \ve(-\sin\vphi\, \om_1 +
\cos\vphi\, \om_d))\ d\HHH^{d-1}(\om) =
$$
\begin{equation}\label{cross}
= -\cos^2 \vphi \cdot \int_{S_-^{d-1}} \om_d^2\ (\cos\vphi\, \om_1
+ \sin\vphi\, \om_d) \, \varrho' (\ve (-\sin\vphi\, \om_1 +
\cos\vphi\, \om_d))\ d\HHH^{d-1}(\om).
\end{equation}
Substituting $\vphi = 0$ in (\ref{cross}), one obtains
$$
p_\ve'(0) = -\int_{S_-^{d-1}} \om_d^2\, \om_1\, \varrho' (\ve
\om_d)\ d\HHH^{d-1}(\om).
$$
Using that $S_-^{d-1}$ is invariant and the integrand is
antisymmetric with respect to reflection $\om_1 \mapsto -\om_1$,
one concludes that $p_\ve'(0) = 0$. Next, substituting $\vphi =
\pi/2$ in (\ref{cross}), one obtains
$$
\lim_{u\to+\infty} p_\ve'(u) = p_\ve'(+\infty) = 0 \cdot
\int_{S_-^{d-1}} \om_d^3\, \varrho' (-\ve \om_1)\ d\HHH^{d-1}(\om)
= 0.
$$
Thus, (b) is proved.

Further, one has
\begin{equation}\label{exp}
p_+'(\tan\vphi) - p_-'(\tan\vphi) = \cos^2 \vphi \cdot
\int_{S_-^{d-1}} \om_d^2\ \frac{\pl}{\pl\vphi}\, \Phi (\vphi,
\om_1, \om_d)\ d\HHH^{d-1}(\om),
\end{equation}
\begin{equation}\label{exp'}
p_+'(\tan\vphi) = \cos^2 \vphi \cdot \int_{S_-^{d-1}} \om_d^2\
\frac{\pl}{\pl\vphi}\, \Phi_+ (\vphi, \om_1, \om_d)\
d\HHH^{d-1}(\om),
\end{equation}
where
$$
\Phi (\vphi, \om_1, \om_d) = \varrho (-\sin\vphi~ \om_1 +
\cos\vphi~ \om_d) + \varrho (\sin\vphi~ \om_1 - \cos\vphi~ \om_d).
$$
$$
\Phi_+ (\vphi, \om_1, \om_d) = \varrho (-\sin\vphi~ \om_1 +
\cos\vphi~ \om_d).
$$
Designate
\begin{equation*}
I(c, \vphi) = \int_{\Gam_c}\, \om_d^2\ \frac{\pl}{\pl\vphi}\, \Phi
(\vphi, \om_1, \om_d)\ d\HHH^1(\om_1,\om_d),
\end{equation*}
\begin{equation*}
I_+ (c, \vphi) = \int_{\Gam_c}\, \om_d^2\ \frac{\pl}{\pl\vphi}\,
\Phi_+ (\vphi, \om_1, \om_d)\ d\HHH^1(\om_1,\om_d),
\end{equation*}
where $\Gam_c = \{ (\om_1,\om_d) : \om_1^2 + \om_d^2 = c^2,\ \om_d
< 0 \}$. Let us prove that
\begin{equation}\label{star}
\text{for any $c \in (0,\, 1)$ and } \vphi \in (0,\, \pi/2), \ \ \
I(c, \vphi) < 0 \ \text{ and } \ I_+(c, \vphi) < 0;
\end{equation}
then, integrating $I(1 - |\tilde\om|^2, \vphi)$ and $I_+(1 -
|\tilde\om|^2, \vphi)$ over $\tilde\om = (\om_2, \ldots,
\om_{d-1})$ and multiplying by $\cos^2 \vphi$, one will conclude
that the right hand sides of (\ref{exp}) and of (\ref{exp'})  are
negative, and so, (c) and the first inequality in (d) are true.

Parametrize the curve $\Gam_c$ by $\om_1 = c \cos\theta$,\, $\om_d
= -c \sin\theta$,\, $\theta \in [0,\, \pi]$, then
\begin{equation*}
I(c, \vphi) = \int_0^\pi c^2 \sin^2 \theta\,
\frac{\pl}{\pl\vphi}\, \big[ \varrho (c\sin (\vphi + \theta)) +
\varrho (-c\sin (\vphi + \theta)) \big]\, c\, d\theta = c^3 \III_1
+ c^3 \III_2,
\end{equation*}
where
\begin{equation}\label{i1}
\III_1 = \int_0^{\pi-2\vphi} \sin^2 \theta\,
\frac{\pl}{\pl\vphi}\, \big[ \ldots \big]\, d\theta,
\end{equation}
\begin{equation}\label{i2}
\III_2 = \int_{\pi-2\vphi}^\pi \sin^2 \theta\,
\frac{\pl}{\pl\vphi}\, \big[ \ldots \big]\, d\theta,
\end{equation}
and
\begin{equation*}
I_+(c, \vphi) = \int_0^\pi c^2 \sin^2 \theta\,
\frac{\pl}{\pl\vphi}\, \varrho (-c\sin (\vphi + \theta))\, c\,
d\theta = c^3 \III_1^+ + c^3 \III_2^+,
\end{equation*}
where
\begin{equation}\label{i1+}
\III_1^+ = \int_0^{\pi-2\vphi} \sin^2 \theta\,
\frac{\pl}{\pl\vphi}\, \varrho (-c\sin (\vphi + \theta))\,
d\theta,
\end{equation}
\begin{equation}\label{i2+}
\III_2^+ = \int_{\pi-2\vphi}^\pi \sin^2 \theta\,
\frac{\pl}{\pl\vphi}\, \varrho (-c\sin (\vphi + \theta))\,
d\theta.
\end{equation}
Changing the variable $\psi = \theta + \vphi - \pi/2$ in
(\ref{i1}), one obtains
$$
\III_1 = \int_{-\pi/2+\vphi}^{\pi/2-\vphi} \cos^2 (\vphi - \psi)\,
\frac{d}{d\psi}\, \big[ \varrho (c \cos \psi) + \varrho (-c \cos
\psi) \big]\, d\psi,
$$
and using the fact that the function $\frac{d}{d\psi} [\ldots]$ under the
sign of integral is odd, one gets
$$
\III_1 = \int_0^{\pi/2-\vphi} (\cos^2 (\vphi - \psi) - \cos^2
(\vphi + \psi))\, \frac{d}{d\psi}\, \big[ \varrho (c \cos \psi) +
\varrho (-c \cos \psi) \big]\, d\psi.
$$
One has $\cos^2 (\vphi - \psi) - \cos^2 (\vphi + \psi) = \sin
2\vphi\, \sin 2\psi > 0$. Taking into account (\ref{circle}), one
also has that $\frac{d}{d\psi}\, \big[ \ldots \big] = -c\sin\psi\
(\varrho'(c\cos\psi) - \varrho'(-c\cos\psi)) < 0$. Hence, $\III_1
< 0$.

Making the same change of variable in (\ref{i1+}), one gets
$$
\III_1^+ = \int_{-\pi/2+\vphi}^{\pi/2-\vphi} \cos^2 (\vphi -
\psi)\, \frac{d}{d\psi}\, \varrho (-c \cos \psi)\, d\psi =
$$
$$
= \int_0^{\pi/2-\vphi} (\cos^2 (\vphi - \psi) - \cos^2 (\vphi +
\psi))\, \frac{d}{d\psi}\, \varrho (-c \cos \psi)\, d\psi.
$$
One has $\cos^2 (\vphi - \psi) - \cos^2 (\vphi + \psi) > 0$, and
$\frac{d}{d\psi}\, \varrho (-c \cos \psi) = c\sin\psi\,
\varrho'(-c\cos\psi)) < 0$, thus $\III_1^+ < 0$.

On the other hand, changing the variable $\chi = \theta + \vphi -
\pi$ in (\ref{i2}), one obtains
$$
\III_2 = \int_{-\vphi}^{\vphi} \sin^2 (\vphi - \chi)\,
\frac{d}{d\chi}\, \big[ \varrho (c \sin \chi) + \varrho (-c \sin
\chi) \big]\, d\chi.
$$
The function $\frac{d}{d\chi} [\ldots]$ is odd, therefore
$$
\III_2 = \int_0^{\vphi} (\sin^2 (\vphi - \chi) - \sin^2 (\vphi +
\chi))\, \frac{d}{d\chi}\, \big[ \varrho (c \sin \chi) + \varrho
(-c \sin \chi) \big]\, d\chi.
$$
One has $\sin^2 (\vphi - \chi) - \sin^2 (\vphi + \chi) = -\sin
2\vphi\, \sin 2\chi < 0$, and $\frac{d}{d\chi}\, \big[ \ldots
\big] = c\cos\chi\, (\varrho'(c\sin\chi) - \varrho'(-c\sin\chi)) >
0$. Hence, $\III_2 < 0$.

Further, as $\pi - 2\vphi \le \theta \le \pi$, one has
$\frac{\pl}{\pl\vphi}\, \varrho (-c\sin (\vphi + \theta)) = -c\cos
(\vphi + \theta)\, \varrho' (-c\sin (\vphi + \theta)) < 0$, and
using (\ref{i2+}), one concludes that $\III_2^+ < 0$.

Thus, the inequalities in (\ref{star}) are proved, and so, (c) and
the first inequality in (d) are true.

Passing to the limit $\vphi \to \pi/2$ in (\ref{p-tg}), one gets
\begin{equation*}
\lim_{\vphi\to\pi/2} p_-(\tan\vphi) = p_-(+\infty) =
-\int_{S_-^{d-1}} \om_d^2\ \varrho (\om_1)\ d\HHH^{d-1}(\om).
\end{equation*}
Thus, to prove the second inequality in (d), one needs to verify
that
\begin{equation*}
\int_{S_-^{d-1}} \om_d^2\ \varrho (\om_1)\ d\HHH^{d-1}(\om)
> \int_{S_-^{d-1}} \om_d^2\ \varrho (\sin\vphi\, \om_1 - \cos\vphi\,
\om_d)\ d\HHH^{d-1}(\om).
\end{equation*}
Denote
\begin{equation*}
J(c, \vphi) = \int_{\Gam_c}\, \om_d^2\, \varrho (\sin\vphi~ \om_1
- \cos\vphi~ \om_d)\, d\HHH^1(\om_1,\om_d) =  c^3 \int_0^\pi
\sin^2 \theta\, \varrho (c\sin (\vphi + \theta))\, d\theta.
\end{equation*}

It suffices to prove that
\begin{equation}\label{star-d}
\text{for any $c \in (0,\, 1)$ and } \vphi \in (0,\, \pi/2), \ \ \
\ J(c, \pi/2) > J(c, \vphi);
\end{equation}
then by integrating $J(1 - |\tilde\om|^2, \pi/2)$ and $J(1 -
|\tilde\om|^2, \vphi)$ over $\tilde\om = (\om_2, \ldots,
\om_{d-1})$, the inequality (\ref{star-d}) will be established.

One has $J(c, \vphi) = J_1 + J_2$,\, $J(c, \pi/2) = J_1^* +
J_2^*$, where
$$
J_1 = \int_0^{\pi/2-\vphi} \sin^2 \theta\, \varrho(c \sin(\vphi +
\theta))\, d\theta, \ \ \ \ J_2 = \int_{\pi/2-\vphi}^\pi \sin^2
\theta\, \varrho(c \sin(\vphi + \theta))\, d\theta,
$$
$$
J_1^* = \int_0^{\pi/2+\vphi} \sin^2 \theta\, \varrho(c
\cos\theta)\, d\theta, \ \ \ \ J_2^* = \int_{\pi/2+\vphi}^\pi
\sin^2 \theta\, \varrho(c \cos\theta)\, d\theta.
$$
As $0 < \theta < \pi/2 - \vphi$, one has $-\cos\theta < 0 <
\sin(\vphi + \theta)$, hence $\varrho(-c\cos\theta) >
\varrho(c\sin(\vphi + \theta)$, thus
\begin{equation}\label{pochtikonec}
J_2^* = \int_0^{\pi/2-\vphi} \sin^2 \theta\,
\varrho(-c\cos\theta)\, d\theta > \int_0^{\pi/2-\vphi} \sin^2
\theta\, \varrho(c \sin(\vphi + \theta))\, d\theta = J_1.
\end{equation}
Further, one has
$$
2J_1^* = \int_0^{\pi/2+\vphi} \sin^2 \theta\, \varrho(c
\cos\theta)\, d\theta + \int_0^{\pi/2+\vphi} \sin^2 (\pi/2 + \vphi
- \theta)\, \varrho(c \cos(\pi/2 + \vphi - \theta)\, d\theta,
$$
$$
2J_2 = \int_0^{\pi/2+\vphi} \sin^2 \theta\, \varrho(c \sin(\theta
- \vphi))\, d\theta + \int_0^{\pi/2+\vphi} \sin^2 (\pi/2 + \vphi -
\theta)\, \varrho(c \sin(\pi/2 - \theta))\, d\theta,
$$
hence
\begin{equation}\label{konec}
2J_1^* - 2J_2 = \int_0^{\pi/2+\vphi} [\sin^2 \theta - \cos^2
(\theta - \vphi)] [\varrho(c \cos\theta) - \varrho(c \sin(\theta -
\vphi)]\, d\theta.
\end{equation}
Taking into account that the function $\varrho$ is monotone
decreasing and that $\sin^2 \theta - \cos^2 (\theta - \vphi) =
(\sin (\theta - \vphi) - \cos \theta) (\sin (\theta - \vphi) +
\cos \theta)$, one concludes that the integrand in (\ref{konec})
is positive, hence $J_1^* > J_2$. From here and from
(\ref{pochtikonec}) it follows that (\ref{star-d}) is true. Lemma
1 is completely proved. \ $\Box$

\subsection*{Proof of lemma \ref{l2}}

{\bf (a)} \ Parametrize the set $S_-^1$ according to $\nu_1 =
-\sin \theta$,\, $\nu_2 = -\cos \theta$,\, $\theta \in [-\pi/2,\,
\pi/2]$, then (\ref{pve}) takes the form
\begin{equation}\label{a*}
p_\ve(\tan\vphi) = \ve \int_{-\pi/2}^{\pi/2} \cos^2 \theta\
\varrho (-\ve \cos(\vphi + \theta))\ d\theta, \ \ \ \ \ve \in \{
-,\ + \}.
\end{equation}
Twice differentiating both parts of this equation with respect to
$\vphi$, one obtains
$$
\frac{p_\ve'' (\tan\vphi)}{\cos^3 \vphi} = -\ve
\int_{-\pi/2}^{\pi/2} \cos^2 \theta\ 2\sin\vphi\
\frac{\pl}{\pl\vphi}\, \varrho (-\ve \cos(\vphi + \theta))\
d\theta +
$$
$$
+ \ve \int_{-\pi/2}^{\pi/2} \cos^2 \theta\ \cos\vphi\
\frac{\pl^2}{\pl\vphi^2}\, \varrho (-\ve \cos(\vphi + \theta))\
d\theta.
$$
Integrating the second integral by parts and taking into account
that $\frac{\pl^k}{\pl\vphi^k}\, \varrho (-\ve \cos(\vphi +
\theta)) = \frac{\pl^k}{\pl\theta^k}\, \varrho (-\ve \cos(\vphi +
\theta))$,\, $k = 1,\ 2$, one gets
$$
\frac{p_\ve'' (\tan\vphi)}{\cos^3 \vphi} = \ve
\int_{-\pi/2}^{\pi/2} 2\cos\theta\, \sin(\theta - \vphi)\,
\frac{\pl}{\pl\theta}\, \varrho (-\ve \cos(\vphi + \theta))\,
d\theta.
$$
Integrating by parts once more and denoting $g_\ve (\vphi) :=
p_\ve'' (\tan\vphi)/ (2\cos^3 \vphi)$, one gets
\begin{equation}\label{astra}
g_\ve (\vphi) = -\ve \int_{-\pi/2}^{\pi/2} \varrho (-\ve
\cos(\vphi + \theta))\, \cos(2\theta - \vphi)\, d\theta.
\end{equation}

Fix the sign ``$+$'' and prove that

(I) \ for $0 < \vphi < \pi/6$,\, $g_+(\vphi) < 0$;

(II) \ for $\vphi \ge 0.3 \pi$,\, $g_+(\vphi) > 0$;

(III) \ for $\pi/6 \le \vphi \le 0.3 \pi$,\, $g_+'(\vphi) > 0$.

The relations (I), (II), and (III) imply that there exists
$\bar{u}_+ \in (1/\sqrt 3,\ \tan (0.3 \pi))$ such that $p_+''(u) <
0$ as $u \in (0,\, \bar u_+)$, and $p_+''(u)
> 0$ as $u \in (\bar u_+,\, +\infty)$.
  \vspace{1mm}

(I) \ Changing the variable $\psi = \theta - \vphi/2 + \pi/4$, one
gets
\begin{equation}\label{asterisk}
g_+ (\vphi) = -\int_{-\pi/4 -\vphi/2}^{3\pi/4 -\vphi/2} \varrho
(-\cos(3\vphi/2 - \pi/4 + \psi))\, \sin 2\psi\, d\psi
\end{equation}
$$
= \LLL_1 + \LLL_2, \ \ \ \ \text{where} \ \ \ \LLL_1 =
-\int_{-\pi/4 -\vphi/2}^{\pi/4 +\vphi/2} \ldots, \ \ \ \LLL_2 =
-\int_{\pi/4 + \vphi/2}^{3\pi/4 -\vphi/2} \ldots\,.
$$
One has
\begin{equation}\label{L1}
\LLL_1 = \int_0^{\pi/4 + \vphi/2} \left[ \varrho \left(-\cos
\left( 3\vphi/2 - \pi/4 - \psi \right) \right) - \varrho
\left(-\cos \left( 3\vphi/2 - \pi/4 + \psi \right) \right) \right]
\, \sin 2\psi \, d\psi.
\end{equation}
One has $0 \le 2\psi \le \pi/2 + \vphi \le \pi$, hence $\sin 2\psi
\ge 0$. Using that $0 < \vphi < \pi/6$, one obtains that $-\pi/2
\le 3\vphi/2 - \pi/4 - \psi < - |3\vphi/2 - \pi/4 + \psi|$,
therefore $-\cos (3\vphi/2 - \pi/4 - \psi) > -\cos (3\vphi/2 -
\pi/4 + \psi)$. The function $\varrho$ monotone decreases,
therefore $\varrho\, (-\cos (3\vphi/2 - \pi/4 - \psi)) < \varrho\,
(-\cos (3\vphi/2 - \pi/4 + \psi))$. Thus, the integrand in
(\ref{L1}) is negative, and so, $\LLL_1 < 0$.

Change the variable $\chi = \psi - \pi/2$ in the integral
$\LLL_2$. One obtains
$$
\LLL_2 = \int_{-\pi/4 + \vphi/2}^{\pi/4 - \vphi/2} \varrho
\left(-\cos \left( 3\vphi/2 + \pi/4 + \chi \right) \right)\, \sin
2\chi \, d\chi =
$$
\begin{equation}\label{L2}
= \int_0^{\pi/4 - \vphi/2} \left[ \varrho \left(-\cos \left(
3\vphi/2 + \pi/4 + \chi \right) \right) - \varrho \left(-\cos
\left( 3\vphi/2 + \pi/4 - \chi \right) \right) \right] \, \sin
2\chi \, d\chi.
\end{equation}
One has $0 \le 3\vphi/2 + \pi/4 - \chi \le 3\vphi/2 + \pi/4 + \chi
\le \pi$, hence $-\cos (3\vphi/2 + \pi/4 - \chi) \le -\cos
(3\vphi/2 + \pi/4 + \chi)$, and so, $\varrho (-\cos (3\vphi/2 +
\pi/4 - \chi)) \ge \varrho (-\cos (3\vphi/2 + \pi/4 + \chi))$.
Therefore the integrand in (\ref{L2}) is negative, and $\LLL_2 \le
0$. Thus, (I) is proved.
 \vspace{1mm}

(II) \ By (\ref{asterisk}), one has
$$
g_+ (\vphi) = \III_1 + \III_2 + \III_3 + \III_4,
$$
where
$$
\III_1 = -\int_{-\vphi}^{\vphi} (...), \ \ \III_2 = -\int_{-\pi/4
-\vphi/2}^{-\vphi} (...), \ \ \III_3 = -\int_{\vphi}^{\pi/4 +
\vphi/2} (...), \ \ \III_4 = -\int_{\pi/4 + \vphi/2}^{3\pi/4
-\vphi/2} (...),
$$
and $(...) = \varrho (-\cos(3\vphi/2 - \pi/4 + \psi))\, \sin
2\psi\, d\psi$.

One has
\begin{equation}\label{I1}
\III_1 = \int_0^{\vphi} [ \varrho (-\cos ( 3\vphi/2 - \pi/4 -
\psi) ) - \varrho (-\cos ( 3\vphi/2 - \pi/4 + \psi) ) ] \, \sin
2\psi \, d\psi;
\end{equation}
using that $\vphi \ge 0.3 \pi$, one easily verifies that as $0 <
\psi < \vphi$,\, $|3\vphi/2 - \pi/4 - \psi| < 3\vphi/2 - \pi/4 +
\psi \le \pi$, hence $-\cos(3\vphi/2 - \pi/4 - \psi) <
-\cos(3\vphi/2 - \pi/4 + \psi)$. Using that $\varrho$ monotone
decreases, one concludes that the integrand in (\ref{I1}) is
positive, thus $\III_1 > 0$.

Next, as $-\pi/4 - \vphi/2 \le \psi \le -\vphi$, one has $\sin
2\psi \le 0$ and $|3\vphi/2 - \pi/4 + \psi| \le 3\vphi/2 - \pi/4 +
\psi + 2\vphi \le \pi$, hence $-\cos(3\vphi/2 - \pi/4 + \psi) \le
-\cos(3\vphi/2 - \pi/4 + \psi + 2\vphi)$, and thus, $\varrho
(-\cos(3\vphi/2 - \pi/4 + \psi)) \le \varrho (-\cos(3\vphi/2 -
\pi/4 + \psi + 2\vphi))$. Therefore,
$$
\III_2 \ge -\int_{-\pi/4 -\vphi/2}^{-\vphi} \varrho
(-\cos(3\vphi/2 - \pi/4 + \psi + 2\vphi))\, \sin 2\psi\, d\psi =
$$
$$
= \int_{\vphi}^{\pi/4 +\vphi/2} \varrho (-\cos(3\vphi/2 - \pi/4 -
\chi + 2\vphi))\, \sin 2\chi\, d\chi\,,
$$
and
$$
\III_2 + \III_3 \ge \int_{\vphi}^{\pi/4 +\vphi/2} [ \varrho\,
(-\cos(3\vphi/2 - \pi/4 - \psi + 2\vphi)) -
$$
\begin{equation}\label{I2+I3}
- \varrho\, (-\cos(3\vphi/2 - \pi/4 + \psi)) ]\, \sin 2\psi\,
d\psi.
\end{equation}
On the other hand, one has
\begin{equation}\label{I4}
\III_4 = \int_{\pi/2}^{3\pi/4 -\vphi/2} [ \varrho\,
(-\cos(3\vphi/2 + 3\pi/4 - \psi)) - \varrho\, (-\cos(3\vphi/2 -
\pi/4 + \psi)) ]\, \sin 2\psi\, d\psi.
\end{equation}
Changing the variable $\theta = \psi - \vphi$ in (\ref{I2+I3}) and
$\theta = \psi - \pi/2$ in (\ref{I4}) and summing both parts of
these relations, one obtains
$$
\III_2 + \III_3 + \III_4 \ge \int_0^{\pi/4 -\vphi/2}
\Psi(\theta)\, d\theta,
$$
where
$$
\Psi(\theta) = [ \varrho\, (-\cos(5\vphi/2 - \pi/4 - \theta)) -
\varrho\, (-\cos(5\vphi/2 - \pi/4 + \theta)) ]\, \sin (2\theta +
2\vphi) -
$$
$$
- [ \varrho\, (-\cos(3\vphi/2 + \pi/4 - \theta)) - \varrho\,
(-\cos(3\vphi/2 + \pi/4 + \theta)) ]\, \sin 2\theta.
$$
Let us show that $\Psi(\theta) \ge 0$; it will follow that $\III_2
+ \III_3 + \III_4 \ge 0$, and thus, (II) will be proved.

One has $0 \le 2\theta \le \pi/2 - \vphi$,\, $\pi/2 - \vphi \le
2\vphi \le 2\theta + 2\vphi \le \pi/2 + \vphi$, hence
\begin{equation}\label{starone}
0 \le \sin 2\theta \le \sin (2\theta + 2\vphi).
\end{equation}
Denote $J_1 (\theta) = \varrho\, (-\cos(5\vphi/2 - \pi/4 -
\theta)) - \varrho\, (-\cos(5\vphi/2 - \pi/4 + \theta))$,\ $J_2
(\theta) = \varrho\, (-\cos(3\vphi/2 + \pi/4 - \theta)) -
\varrho\, (-\cos(3\vphi/2 + \pi/4 + \theta))$. One has
$$
J_1 (\theta) = -\int_{-\theta}^{\theta} \varrho' (-\cos(5\vphi/2 -
\pi/4 + \chi))\, \sin (5\vphi/2 - \pi/4 + \chi)\, d\chi\,,
$$
$$
J_2 (\theta) = -\int_{-\theta}^{\theta} \varrho' (-\cos(3\vphi/2 +
\pi/4 + \chi))\, \sin (3\vphi/2 + \pi/4 + \chi)\, d\chi\,.
$$
Using that $\vphi \ge 0.3 \pi$, one gets $\pi - 2\vphi \le 3\vphi
- \pi/2 \le 5\vphi/2 - \pi/4 + \chi \le 2\vphi$ and $2\vphi \le
3\vphi/2 + \pi/4 + \chi \le \vphi + \pi/2 \le \pi$, hence
\begin{equation}\label{startwo}
\sin (5\vphi/2 - \pi/4 + \chi) \ge \sin (3\vphi/2 + \pi/4 + \chi)
\ge 0
\end{equation}
and $-\cos (5\vphi/2 - \pi/4 + \chi) \le -\cos (3\vphi/2 + \pi/4 +
\chi)$. Taking into account that $\varrho'$ is negative and
monotone increasing, one gets
\begin{equation}\label{starthree}
\varrho' (-\cos (5\vphi/2 - \pi/4 + \chi)) \le \varrho' (-\cos
(3\vphi/2 + \pi/4 + \chi)) \le 0.
\end{equation}
From (\ref{startwo}) and (\ref{starthree}) it follows that $J_1
(\theta) \ge J_2 (\theta) \ge 0$, and taking into account
(\ref{starone}), one concludes that $\Psi(\theta) \ge 0$.

(III) \ One has
$$
g_+' (\vphi) = -\frac{d}{d\vphi}\ \int_{-\pi/4 -\vphi/2}^{3\pi/4
-\vphi/2} \varrho (-\cos(3\vphi/2 - \pi/4 + \psi))\, \sin 2\psi\,
d\psi =
$$
$$
= - \int_{-\pi/4 -\vphi/2}^{3\pi/4 -\vphi/2} \frac{\pl}{\pl\vphi}\
\varrho (-\cos(3\vphi/2 - \pi/4 + \psi))\ \sin 2\psi\, d\psi +
$$
$$
+ \frac 12\, \cos\vphi\, [\varrho(-\cos (\vphi - \pi/2)) -
\varrho(-\cos (\vphi + \pi/2))].
$$
One has $\cos\vphi > 0$,\, $-\cos (\vphi - \pi/2) < 0 < -\cos
(\vphi + \pi/2)$, hence $\varrho(-\cos (\vphi - \pi/2)) -
\varrho(-\cos (\vphi + \pi/2)) > 0$, and thus,
$$
g_+' (\vphi) > - \int_{-\pi/4 -\vphi/2}^{3\pi/4 -\vphi/2}
\frac{\pl}{\pl\vphi}\ \varrho (-\cos(3\vphi/2 - \pi/4 + \psi))\
\sin 2\psi\, d\psi = - (\KKK_1 + \KKK_2 + \KKK_3 + \KKK_4),
$$
where
$$
\KKK_1 = \int_{-\pi/4 -\vphi/2}^{\pi/2 -3\vphi} ..., \ \ \KKK_2 =
\int_{\pi/2 -3\vphi}^0 ..., \ \ \KKK_3 = \int_0^{\pi/4 +\vphi/2}
..., \ \ \KKK_4 = \int_{\pi/4 +\vphi/2}^{3\pi/4 -\vphi/2} ...\,.
$$
In order to prove (III), it suffices to show that $\KKK_1 \le
0$,\, $\KKK_2 \le 0$,\, $\KKK_3 \le 0$,\, and $\KKK_4 \le 0$.

One has
$$
\KKK_1 = \frac 32 \int_{-\pi/4 -\vphi/2}^{\pi/2 -3\vphi} \varrho'
(-\cos(3\vphi/2 - \pi/4 + \psi))\, \sin (3\vphi/2 - \pi/4 +
\psi)\, \sin 2\psi\, d\psi.
$$
Using that $\vphi \le 0.3 \pi$, one easily verifies that $-\pi/4
-\vphi/2 \le \pi/2 -3\vphi$. For $-\pi/4 -\vphi/2 \le \psi \le
\pi/2 -3\vphi$, one has $-\pi/2 \le 3\vphi/2 - \pi/4 + \psi \le
0$,\, $-\pi \le 2\psi \le 0$, hence $\sin (3\vphi/2 - \pi/4 +
\psi) \le 0$,\, $\sin 2\psi \le 0$, and taking into account that
$\varrho' \le 0$, one concludes that $\KKK_1 \le 0$.

Changing the variable $\psi = 3\vphi/2 - \pi/4 + \psi$, one has
$$
\KKK_2 = \frac 32\ \int_{\pi/4 -3\vphi/2}^{3\vphi/2 -\pi/4}
\frac{d}{d\chi}\, \varrho(-\cos\chi)\, \cos (2\chi - 3\vphi)\,
d\chi =
$$
$$
= \frac 32\ \int_0^{3\vphi/2 -\pi/4} \frac{d}{d\chi}\,
\varrho(-\cos\chi)\, [\cos (2\chi - 3\vphi) - \cos (2\chi +
3\vphi)]\, d\chi.
$$
One has $\frac{d}{d\chi}\, \varrho(-\cos\chi) \le 0$,\, $\cos
(2\chi - 3\vphi) - \cos (2\chi + 3\vphi) = 2\sin 2\chi \sin 3\vphi
\ge 0$, hence $\KKK_2 \le 0$.

Further, one has
$$
\KKK_3 = \frac 32\ \int_0^{\pi/4 +\vphi/2} \varrho'(-\cos
(3\vphi/2 - \pi/4 + \psi))\, \sin (3\vphi/2 - \pi/4 + \psi)\, \sin
2\psi\, d\psi.
$$
It is easy to verify that the integrand is negative, hence $\KKK_3
\le 0$.

Changing the variable $\theta = \psi - \pi/2$ in the integral
$\KKK_4$, one obtains
$$
\KKK_4 = - \int_{-\pi/4 +\vphi/2}^{\pi/4 -\vphi/2}
\frac{\pl}{\pl\vphi}\ \varrho (-\cos(3\vphi/2 + \pi/4 + \theta))\
\sin 2\theta\, d\theta =
$$
$$
= \frac 32\ \int_0^{\pi/4 -\vphi/2} [ \varrho' (-\cos(3\vphi/2 +
\pi/4 - \theta))\, \sin (3\vphi/2 + \pi/4 - \theta) -
$$
\begin{equation*}
\label{twostars}
- \varrho' (-\cos(3\vphi/2 + \pi/4 + \theta))\, \sin (3\vphi/2 +
\pi/4 + \theta) ]\, \sin 2\theta\, d\theta.
\end{equation*}
Using that $\vphi \ge \pi/6$, one easily verifies that $(3\vphi/2
+ \pi/4 + \theta) - \pi/2 \ge |(3\vphi/2 + \pi/4 - \theta) -
\pi/2|$, hence $\sin (3\vphi/2 + \pi/4 + \theta) \le \sin
(3\vphi/2 + \pi/4 - \theta)$ and $-\cos (3\vphi/2 + \pi/4 +
\theta) \ge -\cos (3\vphi/2 + \pi/4 - \theta)$, therefore $0 \le
-\varrho' (-\cos(3\vphi/2 + \pi/4 + \theta)) \le -\varrho'
(-\cos(3\vphi/2 + \pi/4 - \theta))$. This implies that $\KKK_4 \le
0$.

{\bf (b)} \ Let us slightly change the notation just introduced:
we shall write down
$$
\varrho(z,V) = \int_0^\infty r^3\, \sigma(\sqrt{r^2 + 2rV z +
V^2})\, dr,
$$
$$
p_-(\tan\vphi,V) = - \int_{-\pi/2}^{\pi/2} \cos^2 \theta\ \varrho
(\cos(\vphi + \theta),V)\ d\theta,
$$
$$
g_-(\vphi,V) = \frac{d^2}{du^2}\rfloor_{u=\tan\vphi} p_- (u,V)/
(2\cos^3 \vphi) = \int_{-\pi/2}^{\pi/2} \varrho (\cos(\vphi +
\theta),V)\, \cos(2\theta - \vphi)\, d\theta,
$$
thus explicitly indicating dependence of these functions on $V$.
Using these formulas, we see that the function $\s^{\al,\bt}(r) =
\s(r) + \al\, \s(\bt r)$ generates the function
$\varrho^{\al,\bt}(z) = \varrho(z,V) + \frac{\al}{\bt^4}
\varrho(z, \bt V)$ and the pressure function $p_-^{\al,\bt}(u) =
p_-(u,V) + \frac{\al}{\bt^4} p_-(u, \bt V)$.

\begin{figure}
\begin{center}
\psfrag{u1}{$u_1$}
\psfrag{u2}{$u_2$}
\psfrag{pm}{$p_-$}
\includegraphics[scale=0.6]{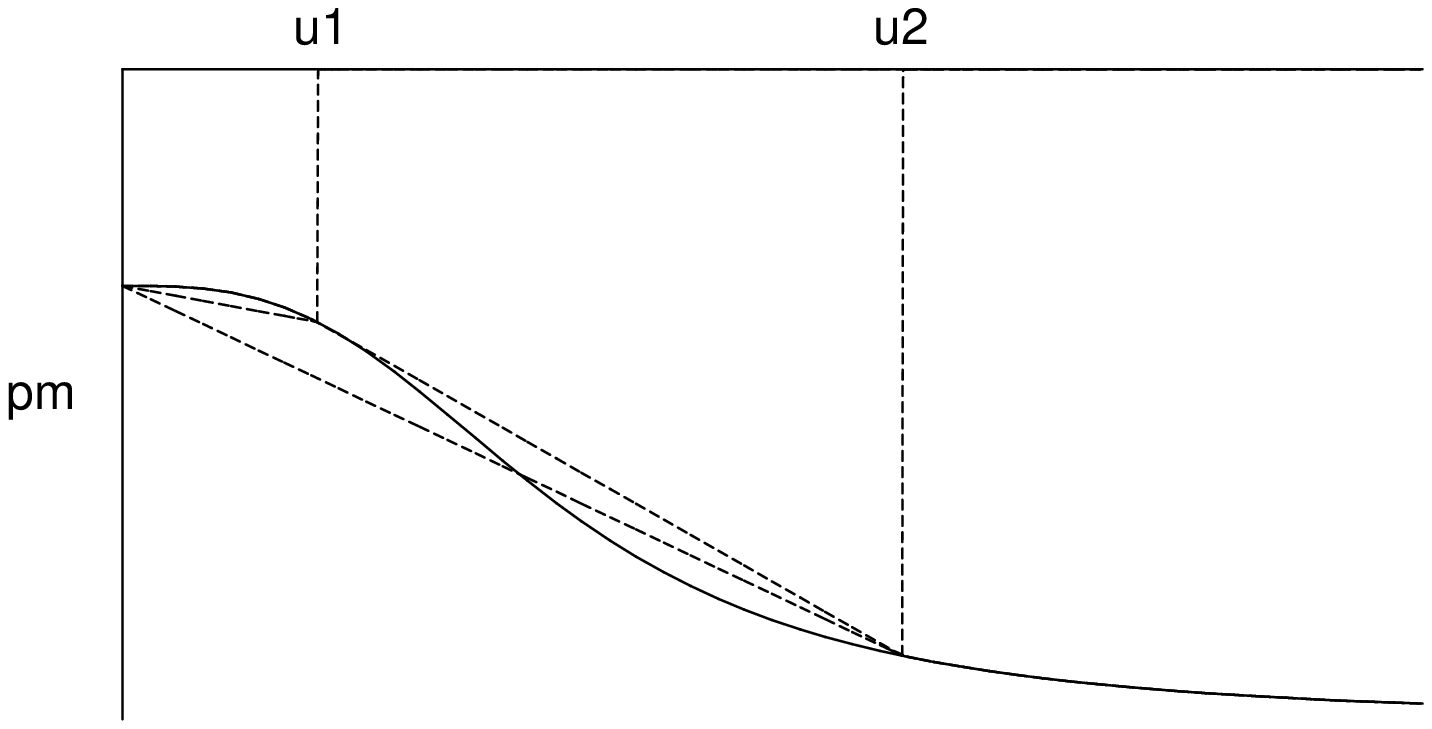}
\end{center}
\caption{}
\label{fig:AppAfig1}
\end{figure}

Suppose that the set $\OOO^0 := \{ u :\, p_-(u,V) > \bar p_-(u,V)
\}$ coincides with an interval $(0,\, u_-^0)$; otherwise, the
hypothesis of lemma \ref{l2}\,(b) is valid for $\al = 0$ and
arbitrary $\bt > 0$. Consider two arbitrary values $u_1 \in (0,\,
u_-^0)$ and $u_2 > u_-^0$, and designate
$$
\Delta_1(p) = \frac{p_-(u_1,V) - p_-(0,V)}{u_1}\,, \ \ \ \ \
\Delta_2(p) = \frac{p_-(u_2,V) - p_-(0,V)}{u_2}\,,
$$
$$
\Delta_{12}(p) = \frac{p_-(u_2,V) - p_-(u_1,V)}{u_2 - u_1}\,;
$$
one has $\Delta_1(p) > \Delta_2(p) > \Delta_{12}(p)$ and
$p_-'(u_2,V) > \Delta_2(p)$ (see Fig.~\ref{fig:AppAfig1}).

Note that the function $g_-(\vphi,V)$ is continuous with respect
to $\vphi$ on $[0,\, \pi/2)$, and has a limit as $\vphi \to \pi/2
- 0$, hence the value $\mathrm{g} := \sup_{\vphi\in [0, \pi/2)}
|g_-(\vphi,V)|$ is finite. Denote
$$
\om = \min \left\{ \frac{\Delta_1 p - \Delta_{12} p}{2}, \
\frac{p_-'(u_2,V) - \Delta_2 p}{2} \right\}.
$$
It is easy to see that if the functions $\check{p}(u) :=
\frac{\al}{\bt^4} p_-(u, \bt V)$ and $\check{g}(\varphi) := $\\
$\frac{\al}{\bt^4}\, \frac{1}{\cos^3 \vphi}\,
\frac{d^2}{du^2}\rfloor_{u = \tan\vphi} p_-(u, \bt V)$ satisfy the
inequalities
\begin{equation}\label{1plus}
|\check{p}'(u)| < \om \ \text{ for } \ 0 \le u \le u_2, \ \ \ \
\check{p}'(u) > -\om \ \text{ for } \ u > u_2,
\end{equation}
and
\begin{equation}\label{2plus}
\check{g}(\vphi) < -\mathrm{g} \ \ \text{ for some } \ \ \vphi >
\arctan u_2,
\end{equation}
then the function $p_-^{\al,\bt}(u) = p_-(u) + \check{p}(u)$ has
the following property: the set $\OOO^{\al,\bt} = \{ u :\,
p_-^{\al,\bt}(u) > \bar p_-^{\al,\bt}(u) \}$ has at least two
connected components; one of them is contained in $(0,\, u_2)$,
and the second one contains $\tan\vphi$ and is contained in
$(u_2,\, +\infty)$.

Denote
\begin{equation}\label{PP}
P^\bt (\vphi) = \max \{ \sup_{0\le u\le\tan\vphi} |p_-'(u,\bt V)|,
\ \sup_{u>\tan\vphi} (-p_-'(u,\bt V)) \}.
\end{equation}
The rest of the text (till the end of Appendix A) is devoted to
the proof of the fact that for any $\ve > 0$ there exist $\phi$,\,
$\vphi$, and $\bt$ such that
\begin{equation}\label{less}
\pi/2 - \ve \le \phi \le \pi/2 \ \ \text{ and } \ \ 0 <
\frac{P^\bt (\vphi)}{-g_-(\vphi, \bt V)} < \ve;
\end{equation}
then, letting $\phi = \tan u_2$,\, $\ve = \om/ \mathrm{g}$, one
can find $\al$ such that (\ref{1plus}) and (\ref{2plus}) are
valid, so the statement of lemma \ref{l2}\,(b) for given $\al$ and
$\bt$ is fulfilled.

\begin{figure}
\begin{center}
\psfrag{U}{\footnotesize{\hspace*{-0.2mm}$\UUU$}}
\psfrag{S}{\footnotesize{\hspace*{-0.4mm}$\Sigma$}}
\psfrag{mbv}{\tiny{\hspace*{-0.8mm}$-\beta V$}}
\psfrag{mbvmr/a}{\tiny{$-\beta V-\frac{R}{\alpha}$}}
\psfrag{mbvmr}{\tiny{\hspace*{-2mm}$-\beta V-R$}}
\psfrag{5fi/4+3pi/8}{\tiny{$\frac{5\varphi}{4}+\frac{3\pi}{8}$}}
\psfrag{fi+pi/2}{\tiny{$\varphi+\frac{\pi}{2}$}}
\pspolygon[linewidth=0.0pt,fillcolor=lightgray,fillstyle=solid](2.49,4.80)(2.71,4.50)%
(2.82,4.25)(2.89,4.00)(2.92,3.75)(2.93,3.45)
(3.26,3.42)(3.27,3.73)(3.24,3.98)(3.19,4.20)(3.18,4.27)%
(3.06,4.55)(2.89,4.85)(2.75,5.02)
\includegraphics[scale=0.5]{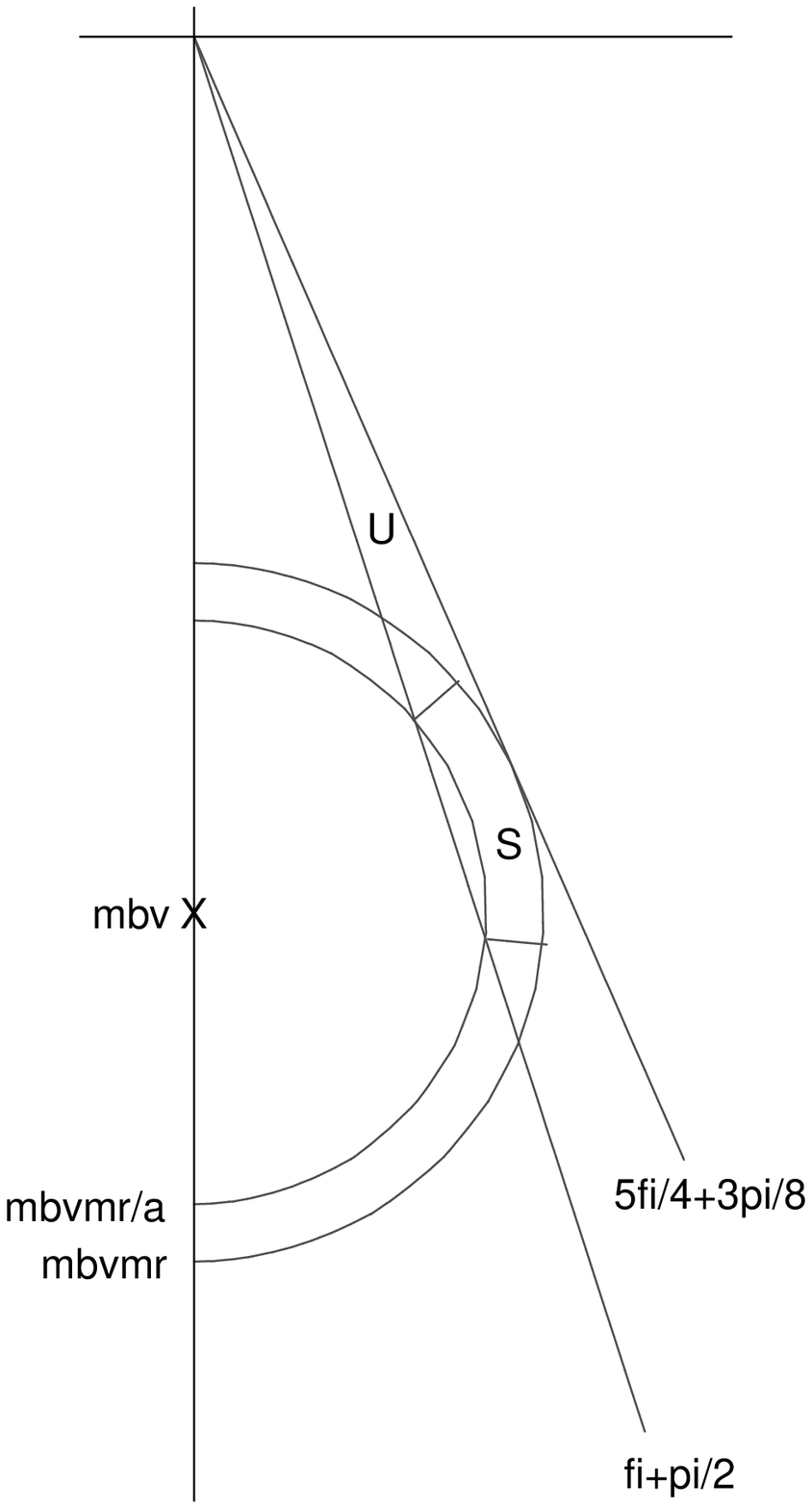}
\end{center}
\caption{}
\label{fig:AppAfig2}
\end{figure}

Let us carry out some auxiliary calculation. For $\vphi \in
(\pi/10,\, \pi/2)$,\, $\bt > 0$ define $R(\vphi,\bt)$ by
\begin{equation}\label{*0}
\frac{R}{\bt V} = \sin \left(\frac{5\pi}{8} -
\frac{5\vphi}{4}\right).
\end{equation}
Then pick out constants $\al > 1$,\, $c > 0$ in such a way that
there exists a circular sector $\Sigma$ of angle $c$, outer radius
$R$, and inner radius $R/\al$, with the center $(0, -\bt V)$, that
belongs to the angle $\UUU := \{ v = (v_1, v_2) :\, v_1 > 0, \
-\tan(\frac{5\pi}{8} - \frac{5\vphi}{4}) \le \frac{v_1}{v_2} \le
-\tan(\frac{\pi}{2} - \vphi) \}$ (see Fig.~\ref{fig:AppAfig2}).
One can choose, in particular, arbitrary $c \in (0,\ \pi/5)$
and $\al = \frac{\cos (c/2)}{\cos (\pi/10)}$.

Next, using the definition of $\varrho$ and passing to the
variables $v_1 = r \sin\chi$,\, $v_2 = r \cos\chi$, one obtains
\begin{equation}\label{*1}
\int_{\frac{5\vphi}{4} + \frac{3\pi}{8}}^{\vphi + \frac{\pi}{2}}
\varrho(\cos\chi, \bt V)\, d\chi = \int \!\!\! \int_\UUU (v_1^2 +
v_2^2)\, \s(|v + \bt V e_2|)\, dv_1 dv_2.
\end{equation}
Taking into account that for $v \in \Sigma$ one has $v_1^2 + v_2^2
\ge (\bt V - R)^2$,\, $\bt V = R/ \sin(\frac{5\pi}{8} -
\frac{5\vphi}{4})$, and $\Sigma \subset \UUU$, and passing to
polar coordinates, one gets
$$
\int \!\!\! \int_\UUU (v_1^2 + v_2^2)\, \s(|v + \bt V e_2|)\, dv_1
dv_2 \ge
$$
$$
\ge R^2 \left( \frac{1}{\sin^2 (\frac{5\pi}{8} -
\frac{5\vphi}{4})} - 1 \right)\, \int \!\!\! \int_\Sigma \s(|v +
\bt V e_2|)\, dv_1 dv_2 =
$$
$$
= R^2 \left( \frac{1}{\sin^2 (\frac{5\pi}{8} - \frac{5\vphi}{4})}
- 1 \right) \cdot c\, \int_{R/\al}^R \s(r)\, r\, dr \ge
$$
\begin{equation}\label{*2}
\ge \left( \frac{1}{\sin^2 (\frac{5\pi}{8} - \frac{5\vphi}{4})} -
1 \right) \cdot c\, \int_{R/\al}^R r^3\, \s(r)\, dr.
\end{equation}
From (\ref{*1}) and (\ref{*2}) one concludes that
\begin{equation}\label{*3}
\int_{\frac{5\vphi}{4} + \frac{3\pi}{8}}^{\vphi + \frac{\pi}{2}}
\varrho(\cos\chi, \bt V)\, d\chi \ge  \left( \frac{1}{\sin^2
(\frac{5\pi}{8} - \frac{5\vphi}{4})} - 1 \right) \cdot c\,
\int_{R/\al}^R r^3\, \s(r)\, dr.
\end{equation}
Besides, one has
$$
\int_{\vphi - \frac{\pi}{2}}^{\frac{3\vphi}{2} + \frac{\pi}{4}}
\varrho(\cos\chi, \bt V)\, d\chi \le 2 \int_0^{\frac{5\vphi}{4} +
\frac{3\pi}{8}} \varrho(\cos\chi, \bt V)\, d\chi \le
$$
\begin{equation}\label{*4}
\le \int \!\!\! \int_\VVV (v_1^2 + v_2^2)\, \s(|v + \bt V e_2|)\,
dv_1 dv_2,
\end{equation}
where $\VVV = \{ v = (v_1, v_2) :\, |v + \bt V e_2| \ge R \}$ is
the complement of the circle of radius $R$ with the center $(0,
-\bt V)$. Taking into account that for $v \in \VVV$,\, $\bt V \le
|v + \bt V e_2|/ \sin(\frac{5\pi}{8} - \frac{5\vphi}{4})$, one has
$|v| \le \bt V + |v + \bt V e_2| \le |v + \bt V e_2|\, (1 + 1/
\sin(\frac{5\pi}{8} - \frac{5\vphi}{4}))$, hence
$$
\int \!\!\! \int_\VVV (v_1^2 + v_2^2)\, \s(|v + \bt V e_2|)\, dv_1
dv_2 \le
$$
$$
\le \left(\frac{1}{\sin(\frac{5\pi}{8} - \frac{5\vphi}{4})} + 1
\right)^2 \int \!\!\! \int_\VVV |v + \bt V e_2|^2\, \s(|v + \bt V
e_2|)\, dv_1 dv_2 =
$$
\begin{equation}\label{*5}
= \left(\frac{1}{\sin(\frac{5\pi}{8} - \frac{5\vphi}{4})} + 1
\right)^2 \cdot 2\pi \int_R^\infty r^2\, \s(r)\, r dr.
\end{equation}
From (\ref{*4}) and (\ref{*5}) one obtains
\begin{equation}\label{*6}
\int_{\vphi - \frac{\pi}{2}}^{\frac{3\vphi}{2} + \frac{\pi}{4}}
\varrho(\cos\chi, \bt V)\, d\chi \le
\left(\frac{1}{\sin(\frac{5\pi}{8} - \frac{5\vphi}{4})} + 1
\right)^2 \cdot 2\pi \int_R^\infty r^3\, \s(r)\, dr.
\end{equation}

Further, by hypothesis (b) of lemma \ref{l2}, for $n > 0$ and for
$r$ sufficiently large the function $\gamma(r) := r^{n+3} \s(r)$
monotonically decreases, hence for $R$ sufficiently large one has
$$
\int_R^\infty r^{-n} \gamma(r)\, dr \le \gamma(R)\,
\frac{R^{-n+1}}{n-1}, \ \ \ \int_{R/\al}^R r^{-n} \gamma(r)\, dr
\ge \gamma(R) (\al^{n-1} - 1)\, \frac{R^{-n+1}}{n-1},
$$
thus
$$
\int_R^\infty r^3 \s(r)\, dr \le \frac{1}{\al^{n-1} - 1}
\int_{R/\al}^R r^3 \s(r)\, dr.
$$
By virtue of arbitrariness of $n$, one concludes that
\begin{equation}\label{*7}
\lim_{R\to +\infty} \frac{\int_R^\infty r^3 \s(r)\,
dr}{\int_{R/\al}^R r^3 \s(r)\, dr} = 0.
\end{equation}
From (\ref{*7}), (\ref{*3}), (\ref{*6}), and (\ref{*0}) it follows
that for any $\vphi \in (\pi/10,\, \pi/2)$
\begin{equation}\label{*8}
\lim_{\beta\to +\infty} \frac{\int_{\vphi -
\frac{\pi}{2}}^{\frac{3\vphi}{2} + \frac{\pi}{4}}
\varrho(\cos\chi, \bt V)\, d\chi}{\int_{\frac{5\vphi}{4} +
\frac{3\pi}{8}}^{\vphi + \frac{\pi}{2}} \varrho(\cos\chi, \bt V)\,
d\chi} = 0
\end{equation}
and
\begin{equation}\label{*8.1}
\int_{\frac{5\vphi}{4} + \frac{3\pi}{8}}^{\vphi + \frac{\pi}{2}}
\varrho(\cos\chi, \bt V)\, d\chi \ge c_1(\vphi)\, \int_{R/\al}^R
r^3\, \s(r)\, dr,
\end{equation}
where $c_1(\vphi) = c\, (1/ \sin^2 (\frac{5\pi}{8} -
\frac{5\vphi}{4}) - 1)$.

Substituting ``$-$'' for $\ve$ and changing the variable $\chi =
\vphi + \theta$  in (\ref{astra}), one gets
\begin{equation}\label{*8.2}
g_- (\vphi, \bt V) = \int_{\vphi-\pi/2}^{\vphi+\pi/2} \varrho
(\cos\chi, \bt V)\, \cos(2\chi - 3\vphi)\, d\chi.
\end{equation}
Denote $\UUU_1 := \{ v = (v_1, v_2) :\, \frac{v_2}{|v_1|} \ge
-\tan\vphi \}$,\ $\VVV_1 := \{ v = (v_1, v_2) :\, |v + \bt V e_2|
\ge \bt V \cos\vphi \}$;\, $\VVV_1$ is the complement of the
circle of radius $\bt V \cos\vphi$ with the center $(0, -\bt V)$.
One has $\UUU_1 \subset \VVV_1$, and for $v \in \VVV_1$\ $|v| \le
|v + \bt V e_2| + \bt V \le |v + \bt V e_2|\, (1 + 1/ \cos\vphi)$,
hence
$$
|g_- (\vphi, \bt V)| \le \int_{-\vphi-\pi/2}^{\vphi+\pi/2} \varrho
(\cos\chi, \bt V)\, d\chi = \int \!\!\! \int_{\UUU_1} (v_1^2 +
v_2^2)\, \s(|v + \bt V e_2|)\, dv_1 dv_2 \le
$$
$$
\le \left(1 + \frac{1}{\cos\vphi} \right)^2 \int \!\!\!
\int_{\VVV_1} |v + \bt V e_2|^2\, \s(|v + \bt V e_2|)\, dv_1 dv_2
=
$$
\begin{equation}\label{*9}
= \left(1 + \frac{1}{\cos\vphi} \right)^2 \cdot 2\pi \int_{\bt
V\cos\vphi}^\infty r^2\, \s(r)\, r dr.
\end{equation}
Taking into account that $\cos(2\chi - 3\vphi) \le 0$ as $\chi \in
[\frac{3\vphi}{2} + \frac{\pi}{4},\, \frac{5\vphi}{4} +
\frac{3\pi}{8}]$,\ $\cos(2\chi - 3\vphi) \le \cos \left(
\frac{3\pi}{4} - \frac{\vphi}{2} \right) = -\cos \left(
\frac{\pi}{4} + \frac{\vphi}{2} \right) < 0$ as $\chi \in
[\frac{5\vphi}{4} + \frac{3\pi}{8},\, \vphi + \frac{\pi}{2}]$, and
using (\ref{*8.2}), one obtains
\begin{equation}\label{*10}
g_- (\vphi, \bt V) \le
\int_{\vphi-\frac{\pi}{2}}^{\frac{3\vphi}{2} + \frac{\pi}{4}}
\varrho (\cos\chi, \bt V)\, d\chi - \cos \left( \frac{\pi}{4} +
\frac{\vphi}{2} \right) \cdot \int_{\frac{5\vphi}{4} +
\frac{3\pi}{8}}^{\vphi+\pi/2} \varrho (\cos\chi, \bt V)\, d\chi.
\end{equation}
From (\ref{*3}), (\ref{*8}), (\ref{*8.1}), and (\ref{*10}) it
follows that for any $\vphi$ and for $\bt$ sufficiently large
\begin{equation}\label{*11}
g_- (\vphi, \bt V) \le -c_2(\vphi)\, \int_{R/\al}^R r^3\, \s(r)\,
dr,
\end{equation}
where $R = R(\vphi,\bt)$ and $c_2(\vphi) = \frac 12\, \cos \left(
\frac{\pi}{4} + \frac{\vphi}{2} \right) \cdot c_1(\vphi)$.

Next, consider $\phi = \phi(\vphi) := \frac{5\vphi}{4} -
\frac{\pi}{8}$; one has $\bt V \cos\phi = \bt V \sin \left(
\frac{5\vphi}{4} + \frac{3\pi}{8} \right) = R$; hence, by
(\ref{*9}), for any $\psi \in [0,\, \phi]$
$$
|g_- (\psi, \bt V)| \le \left(1 + \frac{1}{\cos\psi} \right)^2\,
2\pi \int_{\bt V\cos\psi}^\infty r^3\, \s(r)\, dr \le
$$
\begin{equation}\label{*12}
\le \left(1 + \frac{1}{\cos\phi} \right)^2 \, 2\pi \int_R^\infty
r^3\, \s(r)\, dr.
\end{equation}
Further, one has
$$
p_-'(u, \bt V) = \int_0^u p_-''(\nu, \bt V)\, d\nu = 2\,
\int_0^{\arctan u} g_-(\psi, \bt V)\, \cos\psi\, d\psi,
$$
hence, according to (\ref{*12}),
$$
\sup_{0 \le u \le \tan\phi} |p_-'(u, \bt V)| \le 2\, \int_0^{\phi}
|g_-(\psi, \bt V)|\, d\psi \le
$$
\begin{equation}\label{*12.1}
\le 2\phi\, \left(1 + \frac{1}{\cos\phi} \right)^2\, 2\pi
\int_R^\infty r^3\, \s(r)\, dr
\end{equation}
and
$$
\sup_{u > \tan\phi} \left( -p_-'(u, \bt V) \right) =
-\int_0^{\tan\phi} p_-''(\nu, \bt V)\, d\nu + \sup_{u > \tan\phi}
\left( -\int_{\tan\phi}^u p_-''(\nu, \bt V)\, d\nu \right) \le
$$
\begin{equation}\label{*13}
\le 2\phi\, \left(1 + \frac{1}{\cos\phi} \right)^2\, 2\pi
\int_R^\infty r^3\, \s(r)\, dr + 2(\pi/2 - \phi) \cdot
\sup_{\psi\in [\phi, \pi/2]} (-g_-(\psi, \bt V)).
\end{equation}
Using (\ref{PP}), (\ref{*12.1}), (\ref{*13}), and (\ref{*11}), one
gets that for $\bt$ sufficiently large
\begin{equation}\label{*14}
0 < \frac{P^\bt(\phi)}{\sup_{\psi\in [\phi, \pi/2]} (-g_-(\psi,
\bt V))} \le c_3(\vphi)\, \frac{\int_R^{\infty} r^3\, \s(r)\,
dr}{\int_{R/\al}^R r^3\, \s(r)\, dr} + 2(\pi/2 - \phi),
\end{equation}
where $R = R(\vphi, \bt)$ and $c_3(\vphi) :=
\frac{4\pi}{c_2(\vphi)}\, \phi \left( 1 + \frac{1}{\cos\phi}
\right)^2$. For arbitrary $\ve > 0$, choose $\vphi$ such that $0 <
2(\pi/2 - \phi(\vphi)) < \ve/2$, and then, taking account of
(\ref{*7}), choose $\bt$ such that the first term in the right hand side
of (\ref{*14}) is less that $\ve/2$. The relation (\ref{less}) is
established. \ \ \ $\Box$

\section*{Appendix B}

\renewcommand{\theequation}{B.\arabic{equation}}
\setcounter{equation}{0}

{\large $\mathbf{d = 2}$} \ Here the formula (\ref{formula}) takes
the form
\begin{equation*}
p_\ve(u,V) = \ve \int\!\!\int \frac{(v_1 u + \ve v_2)\!_-^{\,\
2}}{1 + u^2}\ \rho(v)\, dv_1 dv_2, \ \ \ \ \ve \in \{ -,\ + \}.
\end{equation*}
Passing to the polar coordinates $v = (-r \sin\vphi, -\ve r
\cos\vphi)$ and using that $\rho(v) = \s(|v + V e_2|) = \s(r) -
\ve V \cos\vphi \, \s'(r) + o(V)$,\, $V \to 0^+$, one obtains
\begin{equation*}
p_\ve(u,V) = \ve \int_0^{2\pi} \!\! \int_0^\infty \frac{r^2
(\sin\vphi\, u + \cos\vphi)\!_+^{\,\ 2}}{1 + u^2}\ (\s(r) - \ve V
\cos\vphi \, \s'(r) + o(V))\, r\, dr\, d\vphi =
\end{equation*}
\begin{equation}\label{b2}
= \ve \int_0^\infty \s(r)\, r^3\, dr \cdot I^{(2)} - V
\int_0^\infty \s'(r)\, r^3\, dr \cdot J^{(2)} + o(V),
\end{equation}
where
$$
I^{(2)} = \int_0^{2\pi} (\cos(\vphi - \vphi_0))\!_+^{\,\ 2}\,
d\vphi, \ \ \ \ \ J^{(2)} = \int_0^{2\pi} (\cos(\vphi -
\vphi_0))\!_+^{\,\ 2} \cos\vphi\, d\vphi,
$$
$x_+ := \max \{ x,\, 0 \}$,\, $\vphi_0 := \arccos \frac{1}{\sqrt{1
+ u^2}}$. Changing the variable $\psi = \vphi - \vphi_0$ in these
integrals, one obtains
$$
I^{(2)} = \int_{-\pi/2}^{\pi/2} \cos^2 \psi\, d\psi = \pi/2, \ \ \
\ \ J^{(2)} = \cos\vphi_0 \int_{-\pi/2}^{\pi/2} \cos^3 \psi\,
d\psi = \frac{4}{3 \sqrt{1 + u^2}}.
$$
Substituting the obtained values in (\ref{b2}) and using that
$-\int_0^\infty \s'(r)\, r^3\, dr = 3 \int_0^\infty \s(r)\, r^2\,
dr$, one comes to the formula (\ref{asympt d=2}) with coefficients
(\ref{c1c2}).
 \vspace{3mm}

{\large $\mathbf{d = 3}$} \ Formula (\ref{formula}) takes the form
\begin{equation*}
p_\ve(u,V) = \ve \int\!\!\int\!\!\int \frac{(v_1 u + \ve
v_3)\!_-^{\,\ 2}}{1 + u^2}\ \rho(v)\, dv_1 dv_2 dv_3, \ \ \ \ \ve
\in \{ -,\ + \}.
\end{equation*}
Passing to the spherical coordinates $v = (-r \sin\vphi
\cos\theta, -r \sin\vphi \sin\theta, -\ve r \cos\vphi)$, one
obtains
\begin{equation*}
p_\ve(u,V) = \ve \int_0^\pi \!\! \int_0^{2\pi} \!\! \int_0^\infty
\frac{r^2 (\sin\vphi \cos\theta\, u + \cos\vphi)\!_+^{\,\ 2}}{1 +
u^2}\ (\s(r) - \ve V \cos\vphi\, \s'(r) +
\end{equation*}
\begin{equation}\label{b3}
+ o(V))\, r^2\, dr\, d\theta\, \sin\vphi d\vphi = \ve
\int_0^\infty \s(r)\, r^4\, dr \cdot I^{(3)} - V \int_0^\infty
\s'(r)\, r^4\, dr \cdot J^{(3)} + o(V),
\end{equation}
where
$$
I^{(3)} = \int_0^{2\pi} \!\! \int_0^\pi (\cos(\vphi -
\phi_0(\theta)))\!_+^{\,\ 2}\, \frac{1 + u^2 \cos^2 \theta}{1 +
u^2}\, \sin\vphi d\vphi\, d\theta,
$$
$$
J^{(3)} = \int_0^{2\pi} \!\! \int_0^\pi (\cos(\vphi -
\phi_0(\theta)))\!_+^{\,\ 2}\, \frac{1 + u^2 \cos^2 \theta}{1 +
u^2}\, \cos\vphi \sin\vphi d\vphi\, d\theta,
$$
$\phi_0(\theta) = \arccos \frac{1}{\sqrt{1 + u^2 \cos^2 \theta}}$.
Changing the variable $\psi = \vphi - \phi_0(\theta)$, one obtains
$$
I^{(3)} = \int_0^{2\pi} d\theta\ \frac{1 + u^2 \cos^2 \theta}{1 +
u^2}\ \sin\phi_0(\theta) \int_{-\phi_0(\theta)}^{\pi/2} \cos^3
\psi\, d\psi =
$$
$$
= \int_0^{2\pi} \frac{1 + u^2 \cos^2 \theta}{1 + u^2}\ \frac{(1 +
\sin\phi_0(\theta))^2}{3}\, d\theta =
$$
$$
= \int_0^{2\pi} \frac{(\sqrt{1 + u^2 \cos^2 \theta} +
u\cos\theta)^2}{3(1 + u^2)}\, d\theta = \frac{2\pi}{3},
$$
$$
J^{(3)} = \int_0^{2\pi} d\theta\ \frac{1 + u^2 \cos^2 \theta}{1 +
u^2} \int_{-\phi_0(\theta)}^{\pi/2} \cos^2 \psi\ \frac{\sin(2\psi
+ 2\phi_0(\theta))}{2}\ d\psi =
$$
$$
= \int_0^{2\pi} d\theta\ \frac{1 + u^2 \cos^2 \theta}{2(1 + u^2)}\
\left[ \cos 2\phi_0(\theta) \int_{-\phi_0(\theta)}^{\pi/2} \cos^2
\psi\, \sin 2\psi\, d\psi + \right.
$$
$$
\left. + \sin 2\phi_0(\theta) \int_{-\phi_0(\theta)}^{\pi/2}
\cos^2 \psi\, \cos 2\psi\, d\psi \right] = \int_0^{2\pi} \frac{1 +
u\, \phi_0(\theta) \cos\theta}{4(1 + u^2)}\ d\theta =
$$
$$
= \frac{1}{4(1 + u^2)} \int_0^{2\pi} [1 + (u \cos\theta)\, \arctan
(u \cos\theta)]\, d\theta = \frac{\pi}{2\sqrt{1 + u^2}}.
$$
Substituting these values in (\ref{b3}) and using that
$-\int_0^\infty \s'(r)\, r^4\, dr = 4 \int_0^\infty \s(r)\, r^3\,
dr$, one gets the formula (\ref{asympt d=2}) with coefficients
(\ref{k1k2}).


\section*{Acknowledgements}

This work was partially supported by the R\&D Unit CEOC
(Centre for Research in Optimization and Control).


\end{document}